\documentclass[12pt,letterpaper]{article}
\usepackage{amsmath}
\usepackage{amsthm}
\usepackage{amsfonts}
\usepackage{amssymb}
\usepackage{hyperref}
\usepackage{graphicx}
\usepackage{caption2}
\input xy
\xyoption{all}

\newtheorem{thm}{Theorem}[section]

\newtheorem{prop}[thm]{Proposition}
\newtheorem{lemma}[thm]{Lemma}

\theoremstyle{remark}
\newtheorem{rem}[thm]{Remark}
\newtheorem{ex}[thm]{Example}

\theoremstyle{definition}
\newtheorem{dfn}[thm]{Definition}

\newcommand{\R}{\mathbb R}
\newcommand{\C}{\mathbb C}
\newcommand{\N}{\mathbb N}
\newcommand{\Z}{\mathbb Z}

\renewcommand{\L}{\mathcal L}

\newcommand{\M}{\mathcal M}
\newcommand{\J}{\mathcal J}
\newcommand{\E}{\mathcal E}
\newcommand{\PP}{\mathcal P}
\newcommand{\s}{\mathfrak s}
\newcommand{\p}{\mathfrak p}
\newcommand{\D}{\underline D}
\newcommand{\W}{\mathcal W}
\newcommand{\K}{\mathcal K}
\newcommand{\G}{\mathcal G}

\renewcommand{\k}{{\bf k}}
\renewcommand{\d}{{\bf d}}
\renewcommand{\u}{{\bf u}}
\DeclareMathOperator{\ind}{index}
\DeclareMathOperator{\id}{Id} \DeclareMathOperator{\im}{Im}
\DeclareMathOperator{\coker}{coker}
\DeclareMathOperator{\sign}{sign}
\DeclareMathOperator{\rank}{rank}

\long\def\forget#1\forgotten{} %

\author{Jake P. Solomon \thanks{Massachusetts Institute of Technology;
email: jake@math.mit.edu.}}

\title{Intersection Theory on the Moduli Space of Holomorphic Curves with
Lagrangian Boundary Conditions}

\begin{document}
\maketitle

\begin{abstract}
We define a new family of open Gromov-Witten type invariants based on
intersection theory on the moduli space of pseudoholomorphic
curves of arbitrary genus with boundary in a Lagrangian
submanifold. We assume the Lagrangian submanifold arises as the
fixed points of an anti-symplectic involution and has dimension $2$ or $3.$
In the strongly semi-positive genus $0$ case, the new invariants coincide with
Welschinger's invariant counts of real pseudoholomorphic curves.
Furthermore, we calculate the new invariant for the real quintic threefold in
genus $0$ and degree $1$ to be $30.$
\end{abstract}

\tableofcontents

\section{Introduction}\label{sec:intro}
\subsection{The main idea}\label{ssec:main}
In 1985, Gromov initiated the study of pseudoholomorphic curves in
symplectic geometry with his seminal paper \cite{G}. Motivated by
applications of Gromov's techniques in string theory, Witten
developed a systematic way of organizing pseudoholomorphic curve
information, later known as Gromov-Witten invariants \cite{W1,W2}.
Over the following decade, mathematicians including Ruan-Tian \cite{RT},
McDuff-Salamon \cite{MS}, Li-Tian \cite{LT} and Fukaya-Ono \cite{FO},
successfully established a rigorous foundation for Gromov-Witten invariants.
Concurrently, Kontsevich \cite{KM,K} initiated research that
eventually succeeded in calculating the Gromov-Witten invariants in many
situations.

We briefly recall the definition of Gromov-Witten invariants.
Note that throughout the introductory portion of this paper, we assume all
moduli spaces have expected dimension in order to simplify the exposition. Let
$(X,\omega)$ be a symplectic manifold and denote by $\J_\omega$ the set of
$\omega$-tame almost complex structures
on $X.$ Fix a generic $J \in \J_\omega.$  For $d \in H_2(X),$ let
$\overline{\M}_{g,n}(X,d)$ denote the Gromov-compactification of
the moduli space of $J$-holomorphic maps from a surface of genus
$g$ to $X$ representing $d$ together with a choice of $n$ marked
points on the domain. There exist canonical evaluation maps
\begin{equation*}
ev_i: \overline{\M}_{g,n}(X,d) \rightarrow X.
\end{equation*}
Furthermore, denoting by $\overline{\M}_{g,n} =
\overline{\M}_{g,n}(\hbox{pt},0)$ the Deligne-Mumford
compactification of the moduli space of genus $g$ curves with $n$
marked points, there exists a canonical projection
\begin{equation*}
\pi : \overline{\M}_{g,n}(X,d) \rightarrow \overline{\M}_{g,n}.
\end{equation*}
Let $A_i \in H^*(X)$ and $B \in H^*(\overline{\M}_{g,n}).$ Choose
differential forms $\alpha_i \in \Omega^*(X)$ and $\beta \in
\Omega^*(\overline{\M}_{g,n})$ such that $[\alpha_i] = A_i.$ The
genus $g$ Gromov-Witten invariant of $X$ for cohomology classes
$A_i,\,B,$ takes the form of the integral
\begin{equation*}
\int_{\overline{\M}_{g,n}(X,d)} ev_1^*(\alpha_1)\wedge \ldots
\wedge ev_n^*(\alpha_n) \wedge \pi^*(\beta).
\end{equation*}
It follows from Stokes's Theorem and the fact that
$\overline{\M}_{g,n}(X,d)$ is a closed orbifold that this integral
does not depend on the choice of the forms $\alpha_i,\, \beta,$ or
the choice of $J \in \J_\omega.$ Hence, it is an invariant of the
deformation class of $\omega$ parametrized by the cohomology
classes $A_i,\,B.$ Roughly speaking, the Gromov-Witten invariants
of $X$ count the number of $J$-holomorphic maps from a Riemann
surface of fixed genus $g$ to $X$ representing $d$ and
intersecting fixed generic representatives of $PD(A_i).$ The class
$B$ can be used to fix the conformal structure on the domain
Riemann surface or the relative position of the marked points.

Now, let $L \subset X$ be a Lagrangian submanifold, and let
$(\Sigma,\partial \Sigma)$ be a Riemann surface with boundary. For
some time, physicists \cite{AKV,LMV,OV,W4} have predicted the
existence of ``open'' Gromov-Witten invariants counting
pseudoholomorphic maps $(\Sigma,\partial \Sigma) \rightarrow
(X,L)$ satisfying certain incidence conditions. These invariants
would naturally generalize classical Gromov-Witten invariants to
include maps from Riemann surfaces with boundary. Note, however,
that according to \cite{AKV} it may be necessary to specify some
additional structure on $L$ in order to uniquely determine the
invariants. Katz and Liu took a first step toward the mathematical
definition of open invariants given the additional
structure of an $S^1$ action on the pair $(X,L)$ \cite{KL,Liu}.
However, the existence of such an $S^1$ action is a rather
restrictive condition.

Before entering a more detailed discussion, let us briefly
establish some necessary notation. In the following, we denote by
$\Sigma$ a Riemann surface with boundary with fixed conformal
structure. This avoids the issue of degenerations of $\Sigma,$
which the author plans to treat in another paper in the near
future. For $d \in H_2(X,L),$ let $\M_{k,l}(L,\Sigma,d)$ be the
moduli space of configurations of $k$ distinct marked points in
$\partial \Sigma,$ $l$ distinct marked points in $\Sigma$ and
$J$-holomorphic maps $u: (\Sigma, \partial \Sigma) \rightarrow
(X,L)$ such that $u_*([\Sigma,\partial \Sigma]) = d.$ In this
moduli space, points which are equivalent by automorphisms of
$\Sigma$ are identified. We denote by
$\overline{\M}_{k,l}(L,\Sigma,d)$ the Gromov compactification of
$\M_{k,l}(L,\Sigma,d).$ Finally, we denote by
\begin{eqnarray*}
evb_i : \overline{\M}_{k,l}(L,\Sigma,d) \rightarrow L, \qquad i = 1 \ldots k, \\
evi_j: \overline{\M}_{k,l}(L,\Sigma,d) \rightarrow X, \qquad j = 1
\ldots l,
\end{eqnarray*}
the canonical evaluation maps at the marked points.

From a mathematical perspective, two main difficulties have
obstructed progress on open invariants: orientation and
bubbling in codimension one. Indeed, Fukaya et al.\ \cite{FOOO}
showed that $\M_{k,l}(L,\Sigma,d)$ need not be orientable. In the
same paper, they proved orientability if $L$ is orientable and
``relatively spin.'' However, in many interesting examples, i.e.
$(X,L) = (\C P^2, \R P^2),$ $L$ is not orientable and neither is
$\M_{k,l}(L,\Sigma,d).$ In Theorem \ref{thm:or} we show that even
if $L$ is not orientable, under reasonable assumptions, the
orientation bundle of $L^k$ pulls-back to the orientation bundle
of $\M_{k,l}(L,\Sigma,d)$ under the map $\prod_i ev_i.$ This
allows us to pull-back differential forms with values in the
orientation bundle of $L,$ wedge and integrate.

Considerably more troublesome is the problem of bubbling in
codimension one. Put differently, $\overline{\M}_{k,l}(L,\Sigma,d)$ is an
orbifold with corners. Intuitively, one should think of a manifold with many
boundary components. The boundary consists of codimension one strata of the
Gromov-compactification. This stands in contrast to the moduli
space associated to a closed surface $\overline{\M}_{g,n}(X,d),$
which has no boundary since all strata of the Gromov
compactification have codimension two or more. By analogy to the
classical Gromov Witten invariants, we would like to define invariants
parametrized by cohomology classes $A_i \in H^*(L)$ and $C_i \in
H^*(X).$ Choose $\alpha_i \in \Omega^*(L)$ with $[\alpha_i]= A_i$
and $\gamma_j \in \Omega^*(X)$ with $[\gamma_j] = C_i.$ The
desired invariant should take the form
\begin{equation}\label{eq:int}
\int_{\overline{\M_{k,l}}(L,\Sigma,d)} evb_1^*(\alpha_1)\wedge
\ldots \wedge evb_k^*(\alpha_k) \wedge evi_1^*(\gamma_1) \wedge
\ldots \wedge evi_l^*(\gamma_l).
\end{equation}
However, trouble arises in trying to prove independence of the
choices of $\alpha_i,\, \gamma_j.$ For example, suppose
$\alpha_1'$ also satisfies $[\alpha_1'] = A_1.$ Then $\alpha_1 -
\alpha_1' = d\delta$ for some $\delta$ and hence
\begin{equation*}
\begin{array}{l}
evb_1^*(\alpha_1) \wedge \ldots \wedge evi_l^*(\gamma_l) -
evb_1^*(\alpha_1') \wedge \ldots \wedge evi_l^*(\gamma_l) = \\
\hspace{5 cm}= d\left ( evb_1^*(\delta) \wedge evb_2^*(\alpha_2)
\wedge \ldots \wedge evi_l^*(\gamma_l) \right).
\end{array}
\end{equation*}
We would like to integrate the right-hand side of the above
equation over $\overline{\M_{k,l}}(L,\Sigma,d)$ to obtain zero by
Stokes's theorem. However, contributions from the integral of
$ev_1^*(\delta) \wedge ev_2^*(\alpha_2) \wedge \ldots \wedge
ev_l^*(\gamma_l)$ over the boundary of
$\overline{\M_{k,l}}(L,\Sigma,d)$ may spoil this vanishing. So,
the integral (\ref{eq:int}) may depend on the choice of
$\alpha_i,\,\gamma_j.$

Now, let us assume there exists an anti-symplectic involution
\begin{equation*}
\phi : X \rightarrow X, \qquad \phi^* \omega = -\omega,
\end{equation*}
such that $L = Fix(\phi).$ We limit our discussion to the special
case that $A_i \in H^{\dim L}(L),$ $C_j \in H^{\dim X}(X)$ and
$\dim X \leq 6.$ If $L$ is not orientable, we assume $\dim L \leq
4.$ Consequently, we can actually prove independence of
(\ref{eq:int}) from the choice of $\alpha_i,\, \gamma_j.$

We proceed to explain the idea of the proof. The extra structure
$\phi$ enters the definition of the invariants through the almost
complex structure. Indeed, we define
\begin{equation*}
\J_{\omega,\phi} : = \{ J \in \J_\omega | \phi^* J = - J \}.
\end{equation*}
In the following, we fix a generic $J \in \J_{\omega,\phi}.$ Let
$\overline{\M_{k,l}}(L,\Sigma,d)^{(1)}$ denote the union of the
codimension one strata of $\overline{\M_{k,l}}(L,\Sigma,d),$ that
is, strata consisting of two-component stable-maps. Think of
$\overline{\M_{k,l}}(L,\Sigma,d)^{(1)}$ as the boundary of
$\overline{\M_{k,l}}(L,\Sigma,d).$ Recall that it may have many
connected components. We identify a subset of these components,
the union of which we refer to as
$\overline{\M_{k,l}}(L,\Sigma,d)^{(1a)},$ satisfying the following
properties:
\begin{itemize}
\item There exists an orientation reversing involution
$\widetilde\phi_2$ of $\overline{\M_{k,l}}(L,\Sigma,d)^{(1a)}$
that does not preserve any single connected component. Hence the
quotient
\begin{equation*}
\widehat{\M_{k,l}}(L,\Sigma,d) : =
\overline{\M_{k,l}}(L,\Sigma,d)/\widetilde\phi_2(x) \sim x
\end{equation*}
carries a natural orientation. Here, we use the $\phi$ invariance
of $J.$

\item The forms $evb_i^*(\alpha_i),\, evi_j^*(\gamma_j)$ and $evb_1^*(\delta)$
descend naturally to
\[
\widehat{\M_{k,l}}(L,\Sigma,d)
\]
under the assumption that the $\gamma_j$ are $\phi$ invariant.

\item The differential form $evb_1^*(\delta) \wedge
evb_2^*(\alpha_2) \wedge \ldots \wedge evi_l^*(\gamma_l)$ has
support away from the boundary of
$\widehat{\M_{k,l}}(L,\Sigma,d).$
\end{itemize}
The independence of the integral (\ref{eq:int}) of the choices of
$\beta_i$ and $\gamma_j$ follows immediately from Stokes's
theorem: Indeed, we may replace the domain of integration in
(\ref{eq:int}) with $\widehat{\M_{k,l}}(L,\Sigma,d).$ Since
$evb_1^*(\delta) \wedge evb_2^*(\alpha_2) \wedge \ldots \wedge
evi_l^*(\gamma_l)$ vanishes on the boundary of
$\widehat{\M_{k,l}}(L,\Sigma,d),$
\begin{equation*}
\int_{\widehat{\M_{k,l}}(L,\Sigma,d)} d\left( evb_1^*(\delta)
\wedge evb_2^*(\alpha_2) \wedge \ldots \wedge evi_l^*(\gamma_l)
\right) = 0
\end{equation*}
as desired. Independence of the choice of $J \in \J_{\omega,\phi}$
follows by a similar argument. Note however, that a priori the
invariant so obtained should depend on the choice of $\phi.$ We
may interpret the choice of $\phi$ as the extra parameter involved
in defining open invariants predicted by \cite{AKV}.

\subsection{The definition}\label{sec:def}
In the following, we denote by $(X, \omega)$ a symplectic manifold
of dimension $2n$ and by $L \subset X$ a Lagrangian submanifold.
Let $\J_\omega$ denote the set of $\omega$-tame almost complex
structures on $X,$ and let $J\in \J_\omega.$  Let $\PP$ denote the
set of $J$-anti-linear inhomogeneous perturbation terms generalizing
those introduced by Ruan and Tian in \cite{RT}, and let $\nu \in \PP.$ See
Section \ref{sec:def2} for more details. Fix a Riemann surface with boundary
$(\Sigma,\partial \Sigma),$ let $\M_\Sigma$
denote the moduli space of conformal structures on
$(\Sigma,\partial \Sigma),$ and fix $j \in \M_\Sigma.$ Suppose
$\partial \Sigma = \coprod_{a=1}^m (\partial \Sigma)_a,$ where
$(\partial \Sigma)_a \simeq S^1.$ Let
\begin{equation*}
\d = (d,d_1,\ldots, d_m) \in H_2(X,L)\oplus H_1(L)^{\oplus m},
\end{equation*}
let $\k = (k_1,\ldots,k_m) \in \N^m$ and let $l \in \N.$ By
$\M_{\k,l}(L,\Sigma,\d),$ we denote the moduli space of
$(j,J,\nu)$-holomorphic maps $u:(\Sigma,\partial \Sigma)
\rightarrow (X,L)$ with $k_a$ marked points on $(\partial
\Sigma)_a$ and $l$ marked points on $\Sigma$ such that
$u_*([\Sigma,\partial \Sigma]) = d$ and $u|_{(\partial \Sigma)_a
*}([(\partial \Sigma)_a])= d_a.$ Let
$\overline{\M_{\k,l}}(L,\Sigma,\d)$ denote its Gromov
compactification. There exist natural evaluation maps
\begin{align*}
evb_{ai} &: \overline{\M_{\k,l}}(L,\Sigma,\d) \rightarrow L, \qquad i = 1 \ldots
k_a, \; a = 1 \ldots m, \\
evi_j &: \overline{\M_{\k,l}}(L,\Sigma,\d) \rightarrow X, \qquad
j = 1 \ldots l.
\end{align*}

We now digress for a moment to discuss the notion of a relatively
$Pin^\pm$ Lagrangian submanifold. Let $V\rightarrow B$ be a vector
bundle. Define the characteristic classes $p^\pm(V) \in
H^2(B,\Z/2\Z)$ by
\begin{equation*}
p^+(V) = w_2(V), \qquad p^-(V) = w_2(V) + w_1(V)^2.
\end{equation*}
According to \cite{KT}, $p^\pm(V)$ is the obstruction to the
existence of a $Pin^\pm$ structure on $V.$ See \cite{KT} for a
detailed discussion of the definition of the groups $Pin^\pm$ and
the notion of $Pin^\pm$ structures.

Now, suppose $(X,\omega)$ is a symplectic manifold and $L \subset
X$ is a Lagrangian submanifold. Note that we do not assume $L$ is
the fixed points of an anti-symplectic involution yet. We say that
$L$ is relatively $Pin^\pm$ if
\begin{equation*}
p^\pm(TL)  \in \im\left(i^*: H^2(X) \rightarrow H^2(L)\right).
\end{equation*}
and $Pin^\pm$ if $p^\pm(TL) = 0.$ If $L$ is $Pin^\pm,$ we define
a $Pin^\pm$ structure for $L$ to be a $Pin^\pm$ structure for
$TL.$ If $L$ is relatively $Pin^\pm,$ a relative $Pin^\pm$
structure for $L$ consists of the choice of a triangulation for
the pair $(X,L),$ an oriented vector bundle $V$ over the three skeleton of
$X$ such that $w_2(V) = p^\pm(TL)$ and a $Pin^\pm$ structure on
$TL|_{L^{(3)}} \oplus V|_{L}.$ Note that the definition of a relative
$Pin^\pm$ structure given here directly generalizes the notion of a relative
$Spin$ structure given in \cite{FOOO}.

By the Wu relations \cite{MiS}, if $n \leq 3$ then $p^-(TL) = 0,$ so
that $L$ is always $Pin^-.$ It follows that for the applications considered in
this paper, we need only consider honest $Pin^\pm$ structures.
However, we state Theorem \ref{thm:or} in full generality, since
that requires little extra effort.
\begin{thm}\label{thm:or}
Assume $L$ is relatively $Pin^\pm$ and fix a relative $Pin^\pm$
structure on $(X,L).$ If $L$ is not orientable, assume $k_a \cong w_1(d_a) + 1
\mod 2.$ If $L$ is orientable, fix an orientation.
The relative $Pin^\pm$ structure on $(X,L)$ and the orientation of $L$
if $L$ is orientable canonically determine an isomorphism
\begin{equation*}
\det(T\overline{\M}_{\k,l}(L,\Sigma,\d))
\stackrel{\sim}{\longrightarrow} \bigotimes_{a,i} evb_{ai}^*
\det(TL).
\end{equation*}
\end{thm}
\begin{rem}
This theorem was proved in \cite{FOOO} in the case that
$L$ is orientable.
\end{rem}

Under the assumptions of Theorem \ref{thm:or}, we define an
invariant as follows. Let $H^*(L,\det(TL))$ denote the
cohomology of $L$ with coefficients in the flat line bundle
$\det(TL)$. Poincare duality will hold whether or not $L$ is
orientable. Let $\Omega^*(L,\det(TL))$ denote differential forms on $L$ with
values in $\det(TL),$ and let $\Omega^*(X)$ denote ordinary differential forms
on $X.$ For $a = 1,\ldots m,$ and $i = 1, \ldots, k_a,$ let
$\alpha_{ai}\in \Omega^n(L,\det(TL))$ represent the Poincare
dual of a point in $H^n(L,\det(TL)).$ Furthermore, for $j = 1,
\ldots, l,$ let $\gamma_j \in \Omega^{2n}(X)$ represent the
Poincare dual of twice the point-class. We define
\begin{equation*}
N_{\Sigma,\d,\k,l}: =
\int_{\overline{\M}_{\k,l}(L,\Sigma,\d)}evb_{11}^*\alpha_{11}\wedge\ldots
\wedge evb_{mk_m}^*\alpha_{mk_m}\wedge evi_1^*\gamma_1 \wedge
\ldots \wedge evi_m^*\gamma_m.
\end{equation*}
This integral makes sense because by Theorem \ref{thm:or}, the
integrand is a differential form taking values in the orientation
line bundle of the moduli space over which it is to be integrated.
Let $\mu: H_2(X,L) \rightarrow \Z$ denote the Maslov index as
defined in \cite{CG}. Denote by $g$ the genus of the closed
Riemann surface $\Sigma \cup_{\partial\Sigma} \overline\Sigma$
obtained by doubling $\Sigma.$ Furthermore, we employ the
shorthand $|\k| = k_1 +\ldots + k_m.$ We note that by calculating
the expected dimension of $\overline{\M}_{\k,l}(L,\Sigma,\d),$ it
follows that unless
\begin{equation}\label{eq:dim}
(n-1)(|\k| + 2l) = n (1-g) + \mu(d) - \dim Aut(\Sigma)
\end{equation}
the above integral must vanish.

Now, suppose there exists an anti-symplectic involution $\phi: X
\rightarrow X,$ such that $L = Fix(\phi).$ We define
$\J_{\omega,\phi}$ to be the set of $J \in \J_\omega$ such that
$\phi^* J = - J.$ Define
\begin{equation*}
\Omega_\phi^*(X) : = \{\gamma \in \Omega^*(X)|\phi^*\gamma =
\gamma\}.
\end{equation*}
Furthermore, define $\tilde h = h \circ r$ where $h: \pi_2(X) \rightarrow
H_2(X)$ is the Hurewicz homomorphism and $r: H_2(X) \rightarrow H_2(X,L)$ is the
natural homomorphism.

Assume that $\dim X \leq 6,$ and if $L$ is not orientable assume $\dim X \leq
4.$ If $\dim X = 4,$ assume that $k_a \cong w_1(d_a) + 1 \mod 2.$ Note that
these assumptions imply the hypothesis of Theorem \ref{thm:or}. If $\Sigma =
D^2$ and $k = 0$ assume that
\begin{equation}\label{eq:aod}
d \notin \im\left(\tilde h : \pi_2(X) \rightarrow H_2(X,L)\right).
\end{equation}
This is necessary to avoid a certain type of bubbling that requires taking into
consideration real curves with empty real part. The author intends to treat such bubbling
in another paper in the near future.
\begin{thm}\label{thm:inv}
The integers $N_{\Sigma,\d,\k,l}$ do not depend on the choice of
$J \in \J_{\omega,\phi},\newline \: \nu \in \PP, \: j \in
\M_\Sigma,$ the choice of $\alpha_{ai} \in \Omega^n(L,\det(TL))$
or the choice of $\gamma_j \in \Omega^n_\phi(X).$ That is, the
numbers $N_{\Sigma,\d,\k,l}$ are invariants of the triple
$(X,\omega,\phi).$
\end{thm}
\begin{rem}
The condition that $k_a \cong w_1(d_a) + 1 \mod 2$ when $L$ is not
orientable is redundant if $g = 0,$ as it can easily
be derived from the dimension condition $(\ref{eq:dim}).$
\end{rem}
\begin{rem}
The definition of the integers $N_{\Sigma,\d,\k,l}$ does not use
$\phi$ or the condition that $\dim X \leq 6$ in an essential way.
The author believes that there exist far more general conditions
under which similarly defined integers are invariant.
\end{rem}
We now present an example of a non-trivial calculation of these
invariants. See also Section \ref{ssec:re} where we develop the
relationship with Welschinger invariants, for which many
interesting calculations have already been carried out \cite{IKS}.
\begin{ex}\label{ex:q}
Let $(X,L)$ be the pair consisting of the quintic threefold and
its real part. That is, $X := \{ \sum_{i = 0}^4 z_i^5= 0 \}
\subset \C P^4$ equipped with the symplectic form coming from the
restriction of the Fubini-Study metric, and $L := X \cap \R P^4.$
Let $\ell \in H_2(X,L)$ denote the generator with positive
symplectic area. It is not hard to see that $\ell$ satisfies condition
\eqref{eq:aod}. We calculate $N_{D^2,\ell,0,0} = 30.$ This may be
interpreted as the number of oriented lines in the real quintic.
It is interesting to compare this with the classical computation
of $2875$ lines in the complex quintic.
\end{ex}

\subsection{The relationship with real algebraic
geometry}\label{ssec:re}

Real algebraic geometry provides a rich source of examples of
symplectic manifolds admitting anti-symplectic involutions.
Indeed, given any smooth real projective algebraic variety, we can
take $X$ to be its complexification, $\omega$ to be the pull-back
of the Fubini-Study metric and $\phi$ to be complex conjugation.
For this reason, it makes sense to call triples $(X,\omega,\phi)$
real symplectic manifolds. Fix an $\omega$-compatible almost
complex structure $J$ such that $\phi^*J = -J.$ Let $\Sigma$ be a
Riemann surface with an anti-holomorphic involution $c: \Sigma
\rightarrow \Sigma,$ and let $\nu$ be a $c-\phi$ equivariant inhomogeneous
perturbation. We can define real $(J,\nu)$-holomorphic curves to be
$(J,\nu)$-holomorphic maps $u: \Sigma \rightarrow X$ such that $\phi \circ u
\circ c = u \circ a$ for some $a \in Aut(\Sigma,\nu).$ Note that a given Riemann
surface may have several different anti-holomorphic involutions. So, when we
need to specify that a curve is real with respect to a particular
anti-holomorphic involution, we use the terminology $c$-real.

Now, suppose $(\Sigma',\partial\Sigma')$ is a Riemann surface with
boundary such that $\Sigma \simeq \Sigma' \bigcup_{\partial
\Sigma'} \overline{\Sigma'}$ and $c$ acts by exchanging $\Sigma'$
and $\overline{\Sigma'}.$ Any $c$-real $(J,\nu)$-holomorphic curve $u:
\Sigma \rightarrow X,$ must satisfy
\begin{equation*}
u^{-1}(Fix(\phi)) = Fix(c) \simeq \partial \Sigma'.
\end{equation*}
Since, $Fix(c)$ divides $\Sigma$ into $\Sigma'$ and
$\overline{\Sigma'},$ restricting $u$ to either $\Sigma'$ or
$\overline{\Sigma'}$ gives a $(J,\nu)$-holomorphic curve with boundary
in the Lagrangian submanifold $L = Fix(\phi).$  Conversely, given
a $(J,\nu)$-holomorphic map $u' : (\Sigma',\partial \Sigma') \rightarrow
(X, L),$ we can construct a $c$-real $(J,\nu)$-holomorphic map $\Sigma
\rightarrow X$ by gluing $u'$ and $\phi \circ u': (\overline
{\Sigma'}, \overline{\partial \Sigma'}) \rightarrow (X,L)$ along
their common boundary $\partial \Sigma'$ by the Schwarz reflection
principle.

Let us denote by $\overline{\M_n}(X,\Sigma,d)$ the
Gromov-compactification of the space of $(J,\nu)$-holomorphic maps
$\Sigma \rightarrow X$ with $n$ marked points and let
$\overline{\R_c\M_n}(X,d)$ denote its $c$-real locus. Let $r:
H_2(X) \rightarrow H_2(X,L).$ We have just shown there exists a
canonical map
\begin{equation*}
\mathop{\mathop{\coprod_{\k, |k| = n}}_{\d',\, 2d' = r(d),}}_{
\sum_a d_a' = \partial d'}
\overline{\M_{\k,0}}(L,\Sigma',\d')\rightarrow
\overline{\R_c\M_n}(X,\Sigma, d).
\end{equation*}
If $(\Sigma',\nu|_{\Sigma'})$ is not biholomorphic to
$(\overline{\Sigma'},\nu|_{\overline\Sigma'}),$ this
map is $1:1$ on the open stratum. If $(\Sigma',\nu|_{\Sigma'})$ is biholomorphic
to $(\overline{\Sigma'},\nu|_{\overline\Sigma'}),$ then restricted to the open
stratum, this map is a $2:1$ covering map. As an immediate consequence, we have
the following proposition:
\begin{prop}\label{prop:lb} If $(\Sigma',\nu|_{\Sigma'})$ is not biholomorphic
to $(\overline{\Sigma'},\nu|_{\overline\Sigma'}),$ the number of $c$-real
$(J,\nu)$-holomorphic maps $\Sigma \rightarrow X$ intersecting $n$ generic real
points of $X$ is bounded below by
\begin{equation}\label{eq:lb}
\mathop{\mathop{\sum_{\k, |k| = n}}_{\d',\, 2d' = j_*d,}}_{ \sum_a
d_a' = \partial d'} N_{\Sigma',\d',\k,l}.
\end{equation}
If $(\Sigma',\nu|_{\Sigma'})$ is biholomorphic
to $(\overline{\Sigma'},\nu|_{\overline\Sigma'}),$ then we
should take one half of \eqref{eq:lb} as a lower bound instead.
\end{prop}

In \cite{We1,We2,We3}, for strongly semi-positive real symplectic
manifolds, Welschinger defined invariants counting real rational
$J$-holomorphic curves intersecting a generic $\phi$-invariant
collection of marked points. Unlike in the usual definition of
Gromov-Witten invariants, which depends on intersection theory on
the moduli space of $J$-holomorphic curves, Welschinger defined
his curve count by assigning signs to individual curves based on
certain geometric-topological criteria. However, it turns out that
Welschinger's invariants admit the following intersection
theoretic interpretation:
\begin{thm}\label{thm:eq} Let $X$ be a strongly semi-positive real symplectic
manifold satisfying the assumptions of Theorem \ref{thm:inv}. Then
the numbers $N_{D^2,d,k,l}$ are twice the corresponding
Welschinger invariant.
\end{thm}

\subsection{An overview of the paper}
In Section \ref{sec:opcr}, we define the notion of a Cauchy-Riemann Pin
boundary value problem, and prove that the determinant of the associated
Fredholm operator admits a natural orientation. Furthermore, we examine
how this orientation is effected by certain changes in the underlying boundary
value problem and determine when it is preserved by the conjugation morphism. In
Section \ref{sec:or}, we apply the results of Section \ref{sec:opcr} to the
determinant bundle of the $\bar\partial$ operator over the moduli space of
$W^{1,p}$ maps $(\Sigma,\partial\Sigma) \rightarrow (X,L).$ Section
\ref{sec:def2} explains the precise definition of the moduli space of stable
maps we consider, concludes the proof of Theorem \ref{thm:or}, and rigorously
defines the invariants. Section \ref{sec:conj} calculates the sign of the
involution $\widetilde\phi_2$ of Section \ref{ssec:main} on the boundary of the
moduli space. In Section \ref{sec:inv}, assuming that the Maslov index of all
bubbles is strictly positive, we analyze the boundary of the moduli space and
use the results of Section \ref{sec:conj} to prove Theorem \ref{thm:inv}.
Section \ref{sec:km} employs the Kuranishi structure developed in \cite{FO,FOOO}
to treat the case when the Maslov index of bubbles may be $0.$ For the reader
who is not familiar with Kuranishi structures, we summarize the relevant
definitions and ideas in Appendix \ref{sec:aks}. Finally, in Section
\ref{sec:calc}, we introduce the notion of a short exact sequence of
Cauchy-Riemann boundary value problems. We determine how the associated
isomorphisms of determinant lines behave with respect to the orientation of
Section \ref{sec:or} and apply the results to prove Theorem \ref{thm:eq} and
Example \ref{ex:q}.

The author circulated an announcement of the results of this paper and began
giving talks about it in late 2005. Late in the preparation of the full
manuscript, the author learned of work of Cho \cite{CH} that addresses Theorem
\ref{thm:inv} in the strongly semi-positive genus $0$ case when $L$ is
orientable.

\subsection{Acknowledgements}
    I would like to thank my teacher Gang Tian for suggesting to me the
problem of open Gromov-Witten theory and for constant encouragement, support
and helpful advice in my work. In addition, I would like to thank P. Biran, C.
Manolescu, T. Mrowka, L. Polterovich, E. Shustin, and J. Walcher for many
stimulating conversations.

\section{Orienting Cauchy-Riemann operators}\label{sec:opcr}
In this section we analyze the choices necessary to orient the determinant of a
real-linear Cauchy-Riemann operator. In the following, we use the symbol
$\Gamma$ to denote an appropriate Banach space completion of the smooth
sections of a vector bundle. The exact choice of completion will not be
important. If $V \rightarrow B$ is a vector bundle, we denote by $\mathfrak
F(V)$ the principal $O(n)$ bundle with fiber at $x\in B$ given by the set of
orthonormal frames in $V_x.$ We call $\mathfrak F(V)$ the frame-bundle of $V.$
\begin{dfn}
A \emph{$Pin^\pm$ structure} $\p = (P,p)$ on a vector bundle $V\rightarrow B$
consists of principal $Pin^\pm$ bundle $P \rightarrow B$ and a
$Pin^\pm(n)$-$O(n)$ equivariant bundle map
\begin{equation*}
p: P \rightarrow \mathfrak F(V).
\end{equation*}
A morphism of vector bundles with $Pin$ structure $\phi: V\rightarrow V'$ is
said to preserve $Pin$ structure if there exists a lifting $\tilde\phi,$
\begin{equation*}
\xymatrix{
P \ar[r]^{\tilde\phi}\ar[d]^p & P' \ar[d]^{p'} \\
\mathfrak F(V) \ar[r]^\phi & \mathfrak F(V').
}
\end{equation*}
\end{dfn}
\begin{dfn}
A \emph{Cauchy-Riemann $Pin$ boundary value problem}
\[\D = (\Sigma,E,F,\p,D)
\]
consists of
\begin{itemize}
\item
A Riemann surface $\Sigma$ with boundary $\partial \Sigma = \coprod_{a=1}^m
(\partial\Sigma)_a, (\partial\Sigma)_a \simeq S^1.$
\item
A complex vector bundle $E \rightarrow \Sigma.$
\item
A totally real sub-bundle over the boundary
\begin{equation*}
\xymatrix{
 F \ar[r] \ar[d] & E \ar[d] \\
 \partial\Sigma \ar[r] & \Sigma .
}
\end{equation*}
\item A $Pin^\pm$ structure $\p$ on $F.$
\item An orientation of $F|_{(\partial\Sigma)_a}$ for each $a$ such that
$F|_{(\partial\Sigma)_a}$ is orientable.
\item A differential operator
\begin{equation*}
D : \Gamma\left(\vphantom{A^b}(\Sigma,\partial), (E,F)\right) \rightarrow
\Gamma\left(\Sigma,\Omega^{0,1}(E)\right),
\end{equation*}
satisfying, for $\xi \in \Gamma\left(\vphantom{A^b}(\Sigma,\partial),
(E,F)\right)$ and $f \in C^\infty(\Sigma,\R),$
\begin{equation*}
D(f\xi) = f D\xi + (\bar\partial f) \xi.
\end{equation*}
Such a $D$ is known as a \emph{real-linear Cauchy-Riemann operator}.
\end{itemize}
When it does not cause confusion, we will refer to such a collection by the
operator alone, i.e. $D,$ leaving the domain and range implicit.
\end{dfn}
\begin{dfn}
A \emph{morphism of Cauchy-Riemann $Pin$ boundary value problems}
$\underline\phi: \underline D \rightarrow \underline D'$ consists of
\begin{itemize}
\item A biholomorphism $f : \Sigma \rightarrow \Sigma'.$
\item A morphism of bundles $\phi : E \rightarrow E'$ covering $f$ such that
$\phi|_{\partial\Sigma}$ takes $F$ to $F'$ and $\phi\circ D = D' \circ \phi.$
\end{itemize}
Such a morphism is called an \emph{isomorphism} if $\phi$ is an isomorphism
of vector bundles preserving $Pin$ structure and preserving orientation if
$F,F',$ are orientable.
When it causes no confusion, we may refer to such a morphism by the
bundle-morphism component alone, i.e. $\phi.$
\end{dfn}
\begin{dfn}
We define the \emph{determinant line} of a Fredholm operator $D$ by
\begin{equation*}
\det(D) := \Lambda^{\mbox{\tiny max}}(\ker D)\otimes
\Lambda^{\mbox{\tiny max}}(\coker D)^*.
\end{equation*}
If $D$ is a family of Fredholm operators, then we denote by $\det(D)$ the
corresponding line bundle with the natural topology, as explained in, for
example, \cite[appendix A.2]{MS}.
\end{dfn}

We now briefly recall the definition of the Maslov index $\mu(E,F)$ of the
vector bundle pair $(E,F)$ appearing in the definition of a Cauchy-Riemann
boundary value problem in the case that $\partial\Sigma \neq \emptyset.$ Indeed,
if $\partial\Sigma \neq \emptyset,$ we may trivialize $E$ over $\Sigma.$ Writing
\begin{equation*}
\partial\Sigma = \coprod_{a=1}^m (\partial\Sigma)_a, \qquad (\partial \Sigma)_a
\simeq S^1,
\end{equation*}
the restriction of $F$ to each boundary component $(\partial\Sigma)_a$ defines a
loop of totally real subspaces of $\C^n.$ The Maslov index $\mu_a$ of such a
loop was defined in \cite{Ar}. We define
\begin{equation*}
\mu(E,F) = \sum_{a=1}^m \mu_a.
\end{equation*}
It is not hard to see that although $\mu_a$ may depend on the choice of
trivialization of $E,$ the sum $\mu(E,F)$ does not. On the other hand, $\mu_a$
is well defined $\pmod 2,$ and coincides with the first Steifel-Whitney class
$w_1\left(F|_{(\partial\Sigma)_a}\right).$
To properly understand the results of this paper, it is useful
to know the following topological classification of vector-bundle pairs $(E,F).$
\begin{lemma}
Two vector bundle pairs $(E,F)$ and $(E',F')$ of the same dimension admit an
isomorphism
\begin{equation*}
\xymatrix{
E \ar[r]^{\sim}_\phi & E' \\
F \ar[u]\ar[r]^\sim_{\phi|_{\partial \Sigma}} & F' \ar[u]
}
\end{equation*}
if an only if
\begin{equation*}
\mu(E,F) = \mu(E',F'), \qquad w_1(F) =
w_1(F').
\end{equation*}
\end{lemma}

Now, we describe a canonical orientation for the determinant line of a number
of special examples of Cauchy-Riemann operators. Let $\tau \rightarrow \C P^1$
denote the tautological bundle. Let $c': \C P^1 \rightarrow \C P^1$
and $\widetilde c': \tau \rightarrow \tau$ denote the automorphisms induced by
complex conjugation. The fixed points of $c'$ are simply $\R P^1,$
and they divide $\C P^1$ into two copies of $D^2.$ We define $\tau_\R$ to be the
fixed points of $\widetilde c'$ on $\tau|_{\R P^1}.$  Furthermore,
let $\underline \C^n\rightarrow D^2$ denote the trivial bundle and
let $\underline \R^n \hookrightarrow \underline \C^n$ denote the canonical
trivial real sub-bundle. We define basic Cauchy-Riemann $Pin$ boundary value
problems
\begin{align*}
\D(-1,n) &: = (D^2,\tau|_{D^2}\oplus \underline\C^{n-1},\tau_\R \oplus
\underline \R^{n-1},D_{-1,n},\p_{-1}), \\
\D(0,n) &: = (D^2,\underline\C^n, \underline\R^n, D_{0,n},\p_0).
\end{align*}
Here, the $Pin^\pm$ structure $\p_0$ is canonically induced by the splitting
into line bundles. On the other hand, $\p_{-1}$ is not canonical \emph{per se},
but we fix one choice and remain with it for the rest of the paper.
$D_{-1,n}$ and $D_{0,n}$ are taken to be the standard Cauchy-Riemann operators
on these bundles. Since $D_{-1,n}$ and $D_{0,n}$ are surjective, we have
\begin{align*}
\det(D_{-1,n}) &= \Lambda^{\mbox{\tiny max}}(\ker(D_{-1,n})) =
\Lambda^{\mbox{\tiny max}}(\R^{n-1}), \\
\det(D_{0,n}) &= \Lambda^{\mbox{\tiny max}}(\ker(D_{0,n})) =
\Lambda^{\mbox{\tiny
max}}(\R^n).
\end{align*}
So, $\det(D_{-1,n})$ and $\det(D_{0,n})$ admit canonical orientations.

Finally, we will need the following special automorphisms $Q(i)$
of the bundle pair $(\underline\C^2,\underline\R^2)\rightarrow(D^2,\partial
D^2).$ We define the restriction of $Q(i)$ to $\underline\R^2 \rightarrow
\partial D^2$ to be given by a loop in $SO(2)$ with homotopy class $i \in
\pi_1(O(2)) \simeq \Z.$ We define $Q(i)$ to be an arbitrary extension of this
automorphism over the inside of $D^2$ using the fact that the inclusion $SO(2)
\hookrightarrow U(2)$ induces the trivial map on the fundamental group.
\begin{lemma}\label{lem:aut}
If $n \geq 3,$ automorphisms of the bundle pair
\[
(\underline\C^n,\underline\R^n)\rightarrow (D^2,\partial D^2)
\]
preserving $Pin$ structure and orientation are all homotopic to the identity.
If $n = 2,$ automorphisms preserving $Pin$ structure and orientation are
homotopic to $Q(2i),\, i\in \Z.$
For all $n,$ there are two homotopy classes of automorphisms of
\[
(\tau\oplus\underline\C^{n-1},\tau_\R\oplus\underline\R^{n-1}) \rightarrow
(D^2,\partial D^2)
\]
preserving $Pin$ structure. One is homotopic to the identity and the other is
homotopic to $-\id_\tau\oplus \id_{\underline\C^n}.$
\end{lemma}
\begin{rem}\label{rem:n=1}
The assertions of Lemma \ref{lem:aut} are clearly true when $n = 1$
without any reference to $Pin$ structure. This is not surprising because all
automorphisms of a real line bundle preserve $Pin$ structure.
\end{rem}
\begin{proof}[Proof of Lemma \ref{lem:aut}]
For the trivial bundle pair $(\underline\C^n,\underline\R^n),$ homotopy classes
of
automorphisms preserving orientation are classified by
$\pi_2(U(n),SO(n)),$ which is easily calculated from the homotopy
long exact sequence of the pair. In the case $n \geq 3,$ we have
\[
\xymatrix@-1pc{ \pi_2(U(n))\ar[rr]\ar@{=}[d]&&\pi_2(U(n),SO(n))
\ar[rr]^(.55)\sim && \pi_1(SO(n))\ar[rr]^0 \ar@{=}[d] && \pi_1(U(n)). \\
0 &&&& \Z/2\Z }
\]
The automorphisms preserving $Pin$ structure map to $0 \in
\pi_1(SO(n))$ so they are all homotopic to the identity. In the case $n = 2,$
we have the same exact sequence, but $\pi_1(SO(2)) \simeq \Z$ and hence
$\pi_2(U(2),SO(2)) \simeq \Z.$ The automorphisms preserving $Pin$ structure map
to the subgroup $2\Z \subset \Z \simeq \pi_1(SO(2)).$ On the other hand, by
definition, exactly one of the automorphisms $Q(2i)$ maps to each element of the
subgroup $2\Z,$ implying the claim.

In the case of a non-trivial boundary condition, we will have to
make a more explicit argument. Most of the work will be devoted to
showing that the claim of the lemma is true for homotopy classes
of automorphisms of the boundary condition alone. To verify this,
we construct a convenient model for the boundary condition
$\tau_\R \oplus \underline\R^{n-1}.$ Indeed, let $r_1 \in O(n)$ be the
reflection that acts on $\R^n$ by
\[
r(x_1,x_2,\ldots,x_n) = (-x_1,x_2,\ldots,x_n).
\]
We identify
\[
\xymatrix{ \tau_\R\oplus \underline\R^{n-1} \ar[r]^(.3)\sim \ar[d] & \R^n
\times [0,1] / {(x,0) \sim
(r_1(x),1)}\ar[d] \\
\partial D^2 \ar[r] & [0,1]/0\sim 1.
}
\]
Let $\pi : Pin(n) \rightarrow O(n)$ denote the covering map.
Letting $e_i$ denote the standard basis vectors in $\R^n,$ and
thinking of $Pin(n)$ as the group generated by the unit vectors in
the Clifford algebra, we have $\pi(e_1) = r_1.$  So, we may define
a $Pin$ structure on $\tau_\R\oplus \underline\R^{n-1}$ by Diagram \ref{eq:mdl}.
\begin{figure}[ht]
\centering
\includegraphics{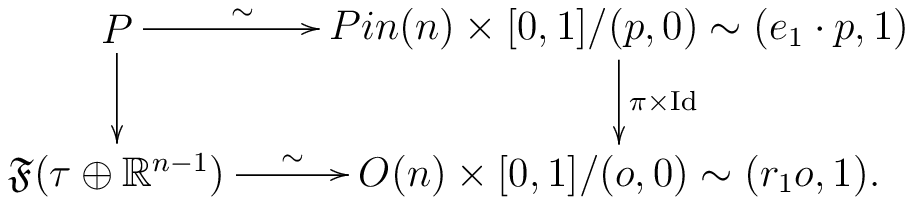}
\caption{}
\label{eq:mdl}
\end{figure}
Here, $\cdot$ denotes Clifford multiplication.  It follows that an automorphism
of $P$ is given by a map
\[
a: [0,1]\rightarrow Pin(n)
\]
such that
\begin{equation}\label{eq:pin}
e_1\cdot a(0) = a(1)\cdot e_1.
\end{equation}
In particular, we see that $\id_{\tau_\R} \oplus \id_{\underline\R^{n-1}}$ lifts
to the
identity automorphism of $P,$ and $-\id_{\tau_\R} \oplus
\id_{\underline\R^{n-1}}$ lifts
to the automorphism of $P$ given by $a(t) = e_1.$

We claim that up to homotopy, these are the only two possibilities. First we
consider the case $n \geq 3.$ Indeed, noting that $Pin(n)$ has two components,
one containing $\id_{Pin(n)}$ and the other containing $e_1,$ it suffices to
show that if automorphisms $a$ and $a'$ map to the same component of $Pin(n)$
then they are homotopic through automorphisms. Indeed, connect $a(0)$ to $a'(0)$
by an arbitrary path $b$. Connect $a(1)$ to $a'(1)$ by $\pm e_1 \cdot b \cdot
e_1,$ where the sign depends on whether we work in $Pin^+$ or $Pin^-.$ The
resulting loop is null homotopic by the simply-connectedness of each component
of $Pin(n).$ Reparameterizing a null-homotopy, we obtain a family of
automorphisms connecting $a$ and $a'$ as desired.

We turn to the case $n = 2.$ Since $Pin(2)$ is not simply-connected, we must be
more carefully. Indeed, topologically, $Pin(2) \simeq S^1 \coprod S^1.$  One
component consists of spinors of the form $cos(\theta) + sin(\theta) e_1 \cdot
e_2$ and the other consist of spinors of the form $cos(\theta) e_1 +
\sin(\theta) e_2.$ So, if we think of $S^1$ as the complex numbers of
unit length, conjugation by $e_1$ acts by complex conjugation on each component
of $Pin(2).$ So, an automorphism of $P$ is given by a path in one of the
two copies of $S^1$ with complex conjugate endpoints. It suffices to show that
any such path is homotopic through similar paths to the constant path at
$\pm 1$ or $\pm e_1.$ Indeed, $\pi(\pm 1) = \id_{\tau_\R \oplus \underline\R}$
and
$\pi(\pm e_1) = -\id_{\tau_\R}\oplus \id_{\underline\R}.$ So, consider the
covering
map $\R \rightarrow S^1$ given by $x \rightsquigarrow e^{2\pi ix}.$ Given a path
in $S^1$ with conjugate endpoints, we may lift it to path
\[
x: [0,1] \rightarrow \R
\]
such that $x(0) \cong -x(1) \pmod 1.$ Since either
\[
\frac{x(0) + x(1)}{2} \cong 0 \pmod 1 \qquad \text{or} \qquad \frac{x(0) +
x(1)}{2} \cong \frac{1}{2} \pmod 1,
\]
linear interpolation between $x$ and $\frac{x(0) + x(1)}{2}$ yields a homotopy
$x_t$ such that
\[
x_t(0) \cong -x_t(1) \pmod 1,\qquad x_0(s) = x(s),\qquad  x_1(s) \cong
0\;\text{or}\;\frac{1}{2} \pmod 1.
\]
Then $\pi(x_t)$ yields the desired homotopy of $a.$

Now we extend our conclusion to automorphisms of the pair. Trivializing $\tau
\oplus \underline\C^{n-1}$ over $D^2,$ we may identify an automorphism of $\tau
\oplus
\underline\C^{n-1}$ with a map
\[
A:D^2 \rightarrow U(n).
\]
We are interested in $A$ that preserve $\tau_\R \oplus \underline\R^{n-1}$
over $\partial D^2,$ such that the induced automorphism of $\tau_\R
\oplus \underline\R^{n-1}$ preserves $Pin$ structure. By the preceding
calculation,
two examples are given by $\id_\tau \oplus \id_{\underline\C^{n-1}}$ and
$-\id_\tau
\oplus \id_{\underline\C^{n-1}}.$ We claim that up to homotopy these are the
only two
examples. Indeed, given $A,$ choose a lift of the induced automorphism on $\tau
_\R \oplus \underline\R^{n-1}$ to $P$ and denote it by $a.$ As just proved, $a$
is
homotopic to either $\id_{Pin(n)}$ or $e_1.$ Denote the homotopy by
\[
\tilde B : [0,1] \rightarrow Aut(P),\qquad \tilde B(0) = a, \qquad \tilde B(1) =
\id_{Pin(n)} \mbox{ or } e_1.
\]
We will construct a homotopy from $A$ to $A': D^2 \rightarrow U(n),$
where $A'$ corresponds to the automorphism
\[\id_\tau \oplus
\id_{\underline\C^{n-1}} \mbox{ or } -\id_\tau \oplus \id_{\underline\C^{n-1}}
\]
as $\tilde B(1) = \id_{Pin(n)}$ or $\tilde B(1) = e_1$ respectively.
Indeed, composing with the given covering map, $\tilde B$ defines a path
\[
B : [0,1] \rightarrow Aut(\tau_\R \oplus \underline\R^{n-1}).
\]
Denote by
\[
i : Aut(\tau_\R \oplus \underline \R^{n-1}) \hookrightarrow Aut(\tau\oplus
\underline\C^{n-1}|_{\partial D^2})
\]
the inclusion given by complexification. Restricting the previously mentioned
trivialization of $\tau\oplus \underline\C^{n-1}$ to $\partial D^2,$ we may
identify the
path $i\circ B$ with a map
\[
\hat B : [0,1] \times \partial D^2 \rightarrow U(n).
\]
Capping off this cylinder with the disk $A$ at one end and the disk $A'$ at the
other end, we obtain a map $S^2 \rightarrow U(n),$ which is well known to be
null-homotopic. Reparameterizing a null-homotopy gives the required homotopy
through automorphisms from $A$ to $A'.$ All these automorphisms preserve the
boundary condition and its $Pin$ structure by the construction of $\hat B.$
\end{proof}

\begin{prop}\label{prop:orl}
The determinant line of a real-linear Cauchy-Riemann $Pin$
boundary value problem $\underline D$ admits a canonical
orientation. If $\underline \phi : \underline D \rightarrow
\underline D'$ is an isomorphism, then the induced morphism
\begin{equation*}
\phi : \det(D) \rightarrow \det(D')
\end{equation*}
preserves the canonical orientation. Furthermore, the canonical orientation
varies continuously in a family of Cauchy-Riemann operators. That is, it
defines a single component of the determinant line bundle over that family.
\end{prop}
\begin{proof}
Near each boundary component $(\partial\Sigma)_a$ choose a closed curve
$\gamma_a$ homotopic to $(\partial\Sigma)_a.$ Degenerate
$\Sigma$ by contracting the curves $\gamma_a$ to points to obtain a nodal
surface $\hat\Sigma.$ $\hat \Sigma$ consists of one closed component
$\tilde\Sigma,$ a disk $\Delta_a$ corresponding to each boundary component
$(\partial\Sigma)_a$ and a nodal point $\hat\gamma_a$ corresponding to each
curve $\gamma_a.$ There exists a continuous map $\pi: \Sigma \rightarrow
\hat\Sigma$ which is a smooth diffeomorphism away from the nodal
points $\hat \gamma_a.$ So, we may define
\[
\hat F = (\pi|_{(\partial\Sigma)_a})^{-1*}F.
\]
At the same time, degenerate $E$ to a vector bundle $\hat E \rightarrow \hat
\Sigma$ such that
\begin{equation}\label{eq:bis}
(\hat E|_{\Delta_a},\hat F|_{\partial \Delta_a}) \simeq
\begin{cases}
(\tau\oplus\underline\C^{n-1},\tau_\R\oplus\underline\R^{n-1}) & \text{if
$w_1(F|_{(\partial\Sigma)_a}) = 1$} \\
(\underline\C^n,\underline\R^n) & \text{if $w_1(F|_{(\partial\Sigma)_a}) = 0.$}
\end{cases}
\end{equation}
We choose the isomorphism \eqref{eq:bis} to preserve orientation in the
orientable case. Moreover, if $n \geq 2$ we choose isomorphism \eqref{eq:bis} to
preserve $Pin$ structure. If $n  = 1,$ we cannot always choose isomorphism
\eqref{eq:bis} to preserve $Pin$ structure. We compensate for this difference
at the next stage of the construction.

Equip $\hat E|_{\Delta_a}$ with the Cauchy-Riemann operator $D_a$ induced by
the isomorphism \eqref{eq:bis} from $D_{-1,n}$ (resp. $D_{0,n}$). The
isomorphism \eqref{eq:bis} induces an orientation on $\det(D_a)$ from the
canonical orientation of $\det(D_{-1,n})$ (resp. $\det(D_{0,n})$). If
$n=1,$ and isomorphism \eqref{eq:bis} does not preserve $Pin$ structure,
reverse the induced orientation on $\det(D_a).$ Choose a
Cauchy-Riemann operator $\tilde D$ on $\hat E|_{\tilde\Sigma}.$ Equip
$\det(\tilde D)$ with the canonical complex orientation. Define an operator
\[
d_{\hat\gamma_a} : \Gamma (\hat E|_{\Delta_a},\hat F|_{\partial \Delta_a})
\oplus \Gamma(\hat E|_{\tilde \Sigma}) \rightarrow E_{\hat \gamma_a}
\]
by
\[
d_{\hat\gamma_a}(\xi,\eta) = \xi(\hat\gamma_a) - \eta(\hat\gamma_a), \qquad \xi
\in \Gamma (\hat E|_{\Delta_a},\hat F|_{\partial \Delta_a}), \qquad \eta \in
\Gamma(\hat E|_{\tilde \Sigma}).
\]
Gluing the $D_a$ with $\tilde D$ at $\hat\gamma_a$ we obtain a Cauchy-Riemann
operator $\#_a D_a \# \tilde D$ on $E$ along with an isomorphism of virtual
vector spaces
\[
\ind\left(\#_a D_a \# \tilde D\right) \simeq \ind\left( \bigoplus_a D_a \oplus
\tilde D \oplus \bigoplus_a d_{\hat\gamma_a}\right),
\]
or equivalently, an isomorphism
\begin{equation}\label{eq:orl}
\det(\#_a D_a \# \tilde D) \simeq \bigotimes_a \det(D_a) \otimes \det(\tilde
D)\otimes \bigotimes_a \det(E_{\hat\gamma_a})^*.
\end{equation}
Since the space of Cauchy-Riemann operators on $E$ is contractible, choosing a
one-parameter family $\mathbf D_t$ with $\mathbf D_0 = D$ and
$\mathbf D_1 = \#_a D_a \# \tilde D$ and
trivializing the line bundle $\det(\mathbf D_t)$ over the family induces an
orientation on $\det(D).$

We claim that the orientation induced on $\det(D)$ is independent of the choice
of isomorphism \eqref{eq:bis}, the choice of $\tilde D,$ and the choice of
$\mathbf D_t .$ First, we prove the independence of the choice of
$\mathbf D_t.$ Indeed, since the space of Cauchy-Riemann operators on $E$ is
contractible, given any two families $\mathbf D_t$ and $\mathbf D_t',$ we
can construct a homotopy between them $\overline{\mathbf D}_{s,t}$, such
that
\[
\overline{\mathbf D}_{0,t} = \mathbf D_t, \qquad \overline{\mathbf
D}_{1,t} = \mathbf D_t'.
\]
Trivializing $\det(\overline{\mathbf D}_{t,s})$ over the homotopy proves
that $\mathbf D_t$ and $\mathbf D'_t$ give the same answer.

Now we turn to proving independence of the choice of isomorphism
\eqref{eq:bis} and the choice of $\tilde D.$ Another choice of isomorphism
\eqref{eq:bis} would induce a different operator $D_a'$ in place of $D_a.$
Also, let $\tilde D'$ be another Cauchy-Riemann operator on $\hat
E|_{\tilde\Sigma}.$ We prove that these new choices induce the same
orientation on $D.$ Choose homotopies $\mathbf D_{a,t}$ and $\mathbf{
\tilde D}_t$ such that
\[
\mathbf D_{a,1} = D_a, \qquad \mathbf D_{a,\frac{1}{2}} = D_a', \qquad
\mathbf{
\tilde D}_1 = \tilde D, \qquad \mathbf{
\tilde D}_{\frac{1}{2}} =
\tilde D'.
\]
We choose the family $\mathbf D_t$ so that, as before, $\mathbf
D_0 = D$ and $\mathbf D_1 = \#_a D_a \# \tilde D,$ but we require also that
\[
\mathbf D_{t} = \#_a \mathbf D_{a,t} \# \mathbf{
\tilde D}_t', \qquad
t\in \left[\frac{1}{2},1\right].
\]
Since this choice of $\mathbf D_{t}$ is as good as any other, it remains
only to show that the orientation on $\det(D_a')$ induced by the isomorphism
\[
\det(D_a') \stackrel{\sim}{\rightarrow} \det(D_{i,n}),\qquad \text{$i = -1$ or
$0,$}
\]
agrees with the orientation induced from $\det(D_a)$ by trivializing
$\det(\mathbf D_{a,t})$ over the interval $[\frac{1}{2},1].$ Similarly, we
must show that the complex orientation on $\det(\tilde D')$ agrees with
orientation induced from $\det(\tilde D)$ by trivializing $\det(\mathbf{
\tilde D}_t)$ over the interval $[\frac{1}{2},1].$ The latter agreement follows
from the compatibility of the topology of the determinant bundle over a family
with the canonical complex orientation.
To see the former agreement, note that the isomorphism \eqref{eq:bis} is
determined up to an automorphism preserving $Pin$ structure of the right hand
bundle pair. In the orientable case where $n \geq 3,$ by Lemma \ref{lem:aut},
all such automorphisms are homotopic to the identity. So, we may assume that
$\mathbf D_{a,t}$ is induced by a family of automorphisms. Then, it suffices
to note the that the determinant bundle is tautologically trivial over a family
of gauge-equivalent operators. The case $n=2$ may be reduced to the higher
dimensional case by stabilizing by a copy of the trivial bundle
pair. Indeed, $Q(2i)\oplus \id_{\underline\R}$ is homotopic to
$\id_{\underline\R^3}.$
In the non-orientable case, we need to consider the additional
possibility that the automorphism is homotopic to $-\id_\tau\oplus
\id_{\underline\C^{n-1}}.$ But $-\id_\tau\oplus \id_{\underline\C^{n-1}}$
clearly
preserves the orientation of $\det(D_{-1}),$ so this possibility does not
effect the argument. The remaining claims of the lemma follow immediately from
the construction.
\end{proof}
\begin{lemma}\label{lem:co}
If the boundary condition $F|_{(\partial\Sigma)_a}$ is orientable, then
reversing the orientation on $F|_{(\partial\Sigma)_a}$ will change the
canonical orientation on $\det(D)$ given in Proposition \ref{prop:orl}.
\end{lemma}
\begin{proof}
This is an immediate consequence of the proof of Proposition
\ref{prop:orl}.
\end{proof}

In the following lemmas, we use the fact that $H^1(B,\Z/2\Z)$ acts naturally
transitively on the set of $Pin$ structures on a vector bundle $V\rightarrow B.$
See \cite{KT}.
\begin{lemma}\label{lem:cp}
Changing the $Pin$ structure $\p$ of a Cauchy-Riemann $Pin$ boundary
value problem by the action of the generator of
\[
H^1((\partial\Sigma)_a,\Z/2\Z) \hookrightarrow H^1(\partial\Sigma,\Z/2\Z)
\]
reverses the canonical orientation of Proposition \ref{prop:orl}.
\end{lemma}
\begin{proof}
By the proof of Proposition \ref{prop:orl} it suffices to consider the
special cases $\underline D(-1,n)$ and $\underline D(0,n).$ If $n=1,$ Lemma
\ref{lem:cp} is tautological. So, we assume $n \geq 2.$ For the case $\underline
D(0,n),$ see \cite[Remark 21.6]{FOOO}. For the case $\underline D(-1,n),$ it
suffices to show that
\[
A := \id_\tau\oplus - \id_{\underline\C} \oplus \id_{\underline\C^{n-2}},
\]
which clearly reverses the orientation of $\det(D_{-1}),$ does not preserve
$Pin$ structure. We use the identification of Diagram \ref{eq:mdl} to show
that that $A|_{\tau_\R\oplus\underline\R^{n-1}}$ does not lift to an
automorphism of $P.$
Indeed, let $r_2 \in O(n)$ be the reflection that acts on $\R^n$ by the formula
\[
r_2(x_1,x_2,x_3,\ldots,x_n) = (x_1,-x_2,x_3,\ldots,x_n).
\]
and, as before, let $\pi:Pin(n)\rightarrow O(n)$ denote the canonical covering
map. The automorphism $A|_{\tau_\R\oplus\underline\R^{n-1}}$ acts on $\mathfrak
F(\tau_\R\oplus\underline\R^{n-1})$ by the explicit formula $(o,t)
\rightsquigarrow
(r_2o,t).$ If this automorphism were to lift to $P,$ it would be given by
left-multiplication by $a\in Pin(n)$ such that $\pi(a) = r_2.$ So,
thinking of $Pin(n)$ as the group generated by the unit vectors in the Clifford
algebra, we would have $a = \pm e_2.$ But this contradicts condition
\eqref{eq:pin}, since
\[
e_1 \cdot e_2 = - e_2 \cdot e_1. \tag*{\qedhere}
\]
\end{proof}

We now introduce a lemma that will play an important role in understanding
the significance of relative $Pin$ structures. Note that any real vector bundle
over a Riemann surface with non-empty boundary $\Sigma$ admits a $Pin$
structure because $\Sigma$ deformation retracts to a wedge of circles.
\begin{lemma}\label{lem:otb}
Let $V \rightarrow \Sigma$ be a real vector bundle over a Riemann
surface with boundary. Consider $\underline D =
(\Sigma,V\otimes\C,V|_{\partial\Sigma},\p,D).$ The canonical orientation of
$\det(D)$ is the same for any $\p$ that arises by restricting a $Pin$ structure
for $V$ over $\Sigma$ to $\partial\Sigma.$
\end{lemma}
\begin{proof}
Let $i:\partial\Sigma \rightarrow \Sigma$ denote the canonical inclusion.
By Lemma \ref{lem:cp}, it suffices to show that any change of $Pin$
structure over $\Sigma$ would change the $Pin$ structure over $\partial\Sigma$
by the action of the sum of the generators of $H^1((\partial\Sigma)_a)$ for an
even number of components $(\partial\Sigma)_a$ of $\partial\Sigma.$ Now, any two
$Pin$ structures of $V$ over $\Sigma$ may be
related by the action of $H^1(\Sigma,\Z/2\Z).$
So we may equivalently show that for all $\alpha \in H^1(\Sigma)$ we
have $i^*\alpha(\partial\Sigma) = 0 \pmod 2.$
But this follows immediately because with $\Z/2\Z$ coefficients $i^*$ is the
dual of $i_*,$ and tautologically $i_*([\partial\Sigma]) = 0.$
\end{proof}

Now, we will calculate the sign of conjugation on the canonical orientation of
the determinant line of a Cauchy-Riemann $Pin$ boundary value problem. More
precisely, given a Riemann surface $\Sigma,$ let $\overline\Sigma$ denote the
same topological surface with conjugate complex structure, and let
\[
t:\Sigma\rightarrow \overline\Sigma
\]
denote the tautological anti-holomorphic map.
Similarly, let $(\overline E, \overline F)$ denote the same real bundle pair
with the opposite complex structure on $E,$ and let
\[
T: E \rightarrow \overline E
\]
denote the tautological anti-complex-linear bundle map. Furthermore, a
Cauchy-Riemann operator $D$ on the bundle $E\rightarrow \Sigma$ is the same as a
Cauchy-Riemann operator $\overline D$ on the bundle $\overline E \rightarrow
\overline \Sigma.$ So, given any Cauchy-Riemann $Pin$ boundary problem
$\underline D,$ we may construct its conjugate $\overline{\underline D}.$
Clearly, we have a tautological map of Cauchy-Riemann $Pin$ boundary value
problems,
\[
\xymatrix{
\Gamma(\Omega^{0,1}(E))\ar[r]^{t^{-1*}\otimes T} & \Gamma(\Omega^{0,1}(\overline
E)) \\
\Gamma(E,F) \ar[r]^T \ar[u]^{D} & \Gamma(\overline E,\overline F)
\ar[u]^{\overline D}
}
\]
which we denote by
\[
\underline T : \underline D \rightarrow \overline{\underline D}.
\]
In the following proposition, we denote by $g_0$ the genus of
$\Sigma/\partial\Sigma$ and we write $n = \dim_\C E = \dim_\R F.$
\begin{prop}\label{prop:cl}
The sign of the induced isomorphism
\[
T: \det(D) \rightarrow \det(\overline D)
\]
relative to the canonical orientation is given by
\begin{align*}
\s^+_T(\D) &:= \frac{\mu(E,F) (\mu(E,F) + 1)}{2} + (1-g_0)n + mn \\
& \qquad + \sum_{a < b} w_1(F)((\partial\Sigma)_a) w_1(F)((\partial\Sigma)_b)
\mod 2,
\end{align*}
for a $Pin^+$ structure and
\begin{align*}
\s^-_T(\D) &:= \frac{\mu(E,F) (\mu(E,F) + 1)}{2} + (1-g_0)n + mn  \\
&\qquad + \sum_{a < b} w_1(F)((\partial\Sigma)_a) w_1(F)((\partial\Sigma)_b) +
w_1(F)(\partial\Sigma) \mod 2,
\end{align*}
for a $Pin^-$ structure.
\end{prop}
\begin{rem}
When $\Sigma = D^2,$ since $\mu(E,F) \cong w_1(F)(\partial\Sigma) \pmod 2,$ we
have the relatively simple formula
\begin{equation}\label{eq:asfd}
\s^\pm_T(\D) \cong \frac{\mu(E,F) (\mu(E,F) \pm 1)}{2} \pmod 2.
\end{equation}
\end{rem}
Before proving the Proposition, we will need the following lemma.
\begin{lemma}\label{lem:-}
The map $\underline T: \D(-1,n) \rightarrow \overline{\D(-1,n)}$
preserves orientation if $\p_{-1}$ is $Pin^+,$ but not if
$\p_{-1}$ is $Pin^-.$
\end{lemma}
\begin{proof}
In this proof, and later in this paper as well, we will need to make use of the
anti-holomorphic involution,
\[
c : \C P^1 \rightarrow \C P^1,  \qquad [z_0 : z_1] \rightsquigarrow [\bar z_0:
-\bar z_1]
\]
and the natural involution $\tilde c$ of the tautological bundle
$\tau$ covering $c.$ We note that $c$ preserves the two
hemispheres of $\C P^1$ which lie on either side of $\R P^1
\subset \C P^1.$ So we may restrict $c,\,\tilde c,$ to either of
the hemispheres $D^2 \subset \C P^1,$ and we denote the
restriction as well by $c,\,\tilde c.$ Furthermore, let $C$ denote
the bundle morphism of the trivial bundle $\underline\C \rightarrow \C P^1$
covering $c$ that acts on the fiber by complex conjugation. The
lemma will follow immediately from Proposition \ref{prop:orl} if
we show that
\[
\tilde c \oplus C^{\oplus n - 1} : \overline{\D(-1,n)} \rightarrow \D(-1,n)
\]
is an isomorphism of Cauchy-Riemann boundary value problems in the $Pin^+$ case
whereas
\[
\tilde c \oplus -C \oplus C^{\oplus n-2}: \overline{\D(-1,n)} \rightarrow
\D(-1,n)
\]
is an isomorphism in the $Pin^-$ case. We treat only the $Pin^-$ case since
the $Pin^+$ case is very similar and not as interesting. The only property of
being an isomorphism of Cauchy-Riemann boundary value problems which is not
immediately evident is the preservation of $Pin^-$ structure. To verify this, we
again work with the explicit model of Diagram \ref{eq:mdl} for $\mathfrak
F(\tau_\R\oplus \underline\R^{n-1})$ and $P.$
In this model, it is not hard to see that at the level of the frame bundle,
\[
\left(\tilde c \oplus -C \oplus C^{\oplus n-2}\right)(o,t) = (r_1r_2o,1-t).
\]
So, a lifting of this map to $P$ must act by
\[
(p,t) \rightsquigarrow (a \cdot p, 1-t)
\]
where $a \in Pin^-(n)$ such that $\pi(a)= r_1r_2.$ It remains to check that
this lifting respects the equivalence relation defining $P.$ Indeed, $(p,0)$
and $(e_1\cdot p, 1)$ represent the same point in $P,$ so we must have
\[
(a \cdot p, 1) \sim (a \cdot e_1 \cdot p, 0).
\]
But that is the same as
\[
a  = e_1 \cdot a \cdot e_1.
\]
Choosing, for example, $a = e_1 \cdot e_2,$ and using the Clifford
multiplication of $Pin^-,$ we verify
\[
e_1 \cdot e_1 \cdot e_2 \cdot e_1 = - e_2 \cdot e_1 = e_1 \cdot e_2.
\tag*{\qedhere}
\]
\end{proof}
\begin{proof}[Proof of Proposition \ref{prop:cl}]
The strategy for the proof is to degenerate $\Sigma$ and $E$ as in the proof of
Proposition \ref{prop:orl} and then treat each irreducible component of the
resulting nodal Riemann surface $\hat\Sigma$ separately. Let $\eta$ be the
number of boundary components $(\partial\Sigma)_a$ for which
$w_1(F)((\partial\Sigma)_a) \cong 1 \pmod 2.$ Then the first Chern class of the
bundle $\hat E|_{\tilde\Sigma}$ from the proof of Proposition \ref{prop:orl} is
given by
\[
c_1\left(\hat E|_{\tilde\Sigma}\right) = \frac{\mu(E,F) + \eta}{2}.
\]
Abbreviating $\mu = \mu(E,F),$ we calculate
\begin{align*}
\frac{\mu(\mu+1)}{2} - \frac{\mu + \eta}{2} &\cong \frac{\mu^2}{2}
- \frac{\eta}{2} \\
    & \cong \frac{\eta^2 - \eta}{2} \\
    & \cong \binom{\eta}{2} \\
    & \cong \sum_{a < b} w_1(F)((\partial\Sigma)_a)
w_1(F)((\partial\Sigma)_b) \pmod 2.
\end{align*}
Here, the second congruence uses the fact that $\mu^2 \simeq \eta^2 \pmod 4.$
So, with $\tilde D$ as in the proof of Proposition \ref{prop:orl}, the
Riemann-Roch theorem gives
\begin{align}
\ind_\C(\tilde D) &= c_1\left(\hat E|_{\tilde\Sigma}\right) +
n(1-g_0) \notag \\
& \cong \frac{\mu(\mu+1)}{2} + n(1-g_0) \notag \\
&\qquad + \sum_{a < b}
w_1(F)((\partial\Sigma)_a)
w_1(F)((\partial\Sigma)_b) \pmod 2. \label{eq:rr}
\end{align}
We now tally the sign of conjugation on each of the tensor factors
on the right-hand side of equation \eqref{eq:orl}. Note that conjugation on a
complex virtual vector space leads to a sign change which is exactly its
dimension $\!\!\pmod 2.$ So, conjugation changes the orientation of $\det(\tilde
D)$ in accordance with the formula \eqref{eq:rr}. Similarly, conjugation changes
the orientation of $\bigotimes_a\det(E_{\hat\gamma_a})$ by $hn.$ The orientation
change for $\det(D_a)$ was calculated in Lemma \ref{lem:-} and
accounts for the difference between $\s^+_T$ and $\s^-_T.$
\end{proof}

\begin{dfn}\label{dfn:sesf}
A \emph{short exact sequence of families of Fredholm operators}
\begin{equation*}
0 \rightarrow D' \rightarrow D \rightarrow D'' \rightarrow 0
\end{equation*}
consists of a
parameter space $B,$ short exact sequences of Banach space bundles over $B,$
\begin{align*}
0\rightarrow X' \rightarrow X \rightarrow X'' \rightarrow 0, \\
0 \rightarrow Y' \rightarrow Y \rightarrow Y'' \rightarrow 0,
\end{align*}
and Fredholm Banach space bundle morphisms
\begin{equation*}
D : X \rightarrow Y,\qquad D' : X' \rightarrow Y', \qquad D'' : X'' \rightarrow
Y'',
\end{equation*}
such that the diagram
\begin{equation*}
\xymatrix{
0 \ar[r] & Y' \ar[r] & Y \ar[r] & Y'' \ar[r] & 0 \\
0 \ar[r] & X' \ar[r]\ar[u]^{D'} & X \ar[r]\ar[u]^{D} & X'' \ar[r]\ar[u]^{D''}
&
0
}
\end{equation*}
commutes.
\end{dfn}

\begin{lemma}\label{lem:ses}
A short exact sequence of families of Fredholm operators
\begin{equation*}
0 \rightarrow D' \rightarrow D \rightarrow D'' \rightarrow 0
\end{equation*}
induces an isomorphism
\begin{equation*}
\det(D') \otimes \det(D'') \stackrel{\sim}{\rightarrow} \det(D).
\end{equation*}
\end{lemma}

\section{Orienting moduli of open stable maps}\label{sec:or}
First, we set the basic assumptions which will hold throughout this
section. Let $(X,\omega)$ be a symplectic manifold with $\dim X = 2n$ and let $L
\subset X$ be a Lagrangian submanifold. In the following, we assume $L$ is
relatively $Pin^\pm$ and fix a relative $Pin^\pm$ structure $\p$ on $L.$
Furthermore, if $L$ is orientable, fix an orientation on $L.$ Let
$(\Sigma,\partial \Sigma)$ denote a Riemann surface with boundary and assume
$\partial \Sigma = \coprod_{a=1}^m (\partial \Sigma)_a,$ where
$(\partial \Sigma)_a \simeq S^1.$ Now, for
\begin{equation*}
\d = (d,d_1,\ldots,d_m) \in H_2(X,L)\oplus H_1(L)^{\oplus m},
\end{equation*}
define $B^{1,p}(L,\Sigma,\d)$ to be the Banach manifold of $W^{1,p}$ maps
\[
u: (\Sigma,\partial \Sigma) \rightarrow (X,L)
\]
such that $u_*([\Sigma,\partial \Sigma]) = d$ and $u|_{(\partial \Sigma)_a
*}([(\partial \Sigma)_a])= d_a.$  Furthermore, define
\begin{equation*}
B^{1,p}_{\k,l}(L,\Sigma,\d) : = B^{1,p}(L,\Sigma,\d) \times \prod_a (\partial
\Sigma)_a^{k_a} \times \Sigma^l \setminus \Delta.
\end{equation*}
Each factor of $(\partial\Sigma)_a$ in the preceding product corresponds
to a marked point on $(\partial \Sigma)_a,$ and each factor of $\Sigma$
corresponds to a marked point on $\Sigma.$ $\Delta$ denotes the subset of
the product in which two marked points coincide. We will use $\vec z = (z_{ai})$
and $\vec w = (w_j)$ to denote marked points in $\partial \Sigma$ and $\Sigma$
respectively, and we use $\u = (u,\vec z, \vec w)$ to denote elements of
$B_{\k,l}^{1,p}(L,\Sigma,\d).$
Note that we may occasionally omit the $L,\,\Sigma, \d,$ from the preceding
notation when it is clear from the context.
There exist canonical evaluation maps
\begin{align*}
evb_{ai} &: \overline{\M_{\k,l}}(L,\Sigma,\d) \rightarrow L, \qquad i = 1
\ldots k_a, \; a = 1 \ldots m, \\
evi_j &: \overline{\M_{\k,l}}(L,\Sigma,\d) \rightarrow X, \qquad
j = 1 \ldots l,
\end{align*}
given by $evb_{ai}(\u) = u(z_{ai})$ and $evi_j(\u) = u(w_j).$
We note that the above notation is also used for the evaluation maps from the
moduli spaces of holomorphic curves, which are just restrictions of the
maps above.

Define the Banach space bundle $\E \rightarrow B^{1,p}_{\k,l}(L,\Sigma,\d)$
fiberwise by
\begin{equation*}
\E_\u := L^p(\Sigma, \Omega^{0,1}(u^*TX))
\end{equation*}
for $\u \in B^{1,p}_{\k,l}(L,\Sigma,\d).$ Now, fix $J \in \J_\omega,$ and $\nu
\in \PP.$ Let
\begin{equation*}
\overline\partial_{(J,\nu)}: B^{1,p}_{\k,l}(L,\Sigma,\d) \rightarrow \E
\end{equation*}
denote the section of $\E$ given by the $\nu$-perturbed Cauchy
Riemann operator. Using the canonical identification between the vertical
tangent spaces of $\E$ and $\E$ itself we define
\begin{equation*}
D := D\overline\partial_{(J,\nu)}: TB^{1,p}_{\k,l}(L,\Sigma,\d) \rightarrow \E
\end{equation*}
to be the vertical component of the linearization of
$\overline\partial_{(J,\nu)}.$ We will denote linearization at any given
$\u \in B^{1,p}_{\k,l}$ by $D_\u.$ Finally, define $\L \rightarrow
B^{1,p}_{\k,l}(L,\Sigma,\d)$ to be the determinant line bundle of the
family of Fredholm operators $D,$
\begin{equation*}
\L := \det(D).
\end{equation*}

The following proposition is a basic ingredient in the proof of Theorem
\ref{thm:or}. Suppose either $L$ is orientable and provided with an
orientation or $k_a \cong w_1(d_a) + 1 \mod 2.$
\begin{prop}\label{prop:or}
The combination of an orientation of $L$ if $L$ is orientable and the choice of
relative $Pin^\pm$ structure $\p$ on $L$ canonically determine an isomorphism
of line bundles
\begin{equation*}
\L \stackrel{\sim}{\longrightarrow} \bigotimes_{a,i} evb_{ai}^*\det(TL).
\end{equation*}
\end{prop}
\begin{proof}
Clearly, the proposition will follow immediately if we succeed in providing
\begin{equation*}
\L' : = \L\otimes\bigotimes_{a,i} evb_{ai}^*\det(TL)^*
\end{equation*}
with a canonical orientation depending only on the orientation of $L$ and $\p.$
We observe that it suffices to canonically orient the fiber $\L'_\u$ over each
$\u \in B_{\k,l}^{1,p}$ individually in a way that varies continuously with
$\u.$

Recall that the relative $Pin$ structure of $L$ specifies a triangulation
of the pair $(X,L).$ Using simplicial approximation, we homotope the map
$u: (\Sigma,\partial \Sigma) \rightarrow (X,L)$ to a map $\hat u:
(\Sigma,\partial\Sigma) \rightarrow (X^{(2)},L^{(2)}).$ Denote the homotopy by
\[
\Phi : [0,1] \times (\Sigma,\partial\Sigma) \rightarrow (X,L).
\]
We claim that the choice of $\Phi$ is unique up to homotopy. Indeed, suppose
$\Phi'$ is another such homotopy. Concatenating $\Phi$ and $\Phi',$ we obtain a
map
\begin{equation*}
\Phi \# \Phi': [-1,1] \times
(\Sigma,\partial\Sigma) \rightarrow (X,L).
\end{equation*}
By simplicial approximation we may homotope $\Phi \#\Phi'$ to map into
$(X^{(3)},L^{(3)}).$ Reparameterizing this homotopy, we obtain a homotopy from
$\Phi$ to $\Phi'.$ We denote the homotopy from $\Phi$ to $\Phi'$ by
\begin{equation*}
\Psi : [0,1]^2 \times (\Sigma,\partial\Sigma) \rightarrow (X,L),
\end{equation*}
such that,
\begin{equation*}
\Psi(0,t) = \Phi(t), \qquad \Psi(1,t) = \Phi', \qquad \Psi(s,0) = u.
\end{equation*}

Now, define
\begin{equation*}
B' : = \left \{ \left . (u,\vec z,\vec w) \in B_{\k,l}^{1,p} \right | u:
(\Sigma, \partial \Sigma) \rightarrow (X^{(3)},L^{(3)}) \right \}.
\end{equation*}
We now prove that the homotopy uniqueness of $\Phi$ implies that it suffices
to orient $\L'|_{B'}.$ Indeed, think of $\Phi$ (resp. $\Psi$) as a map from
$[0,1]$ (resp. $[0,1]^2$) to $B_{\k,l}^{1,p}.$ Given an orientation on
$\L'|_{B'},$ trivializing $\Phi^*\L'$ induces an orientation of $\L'_\u.$ This
orientation agrees with the orientation induced by any other homotopy $\Phi',$
because we may trivialize $\Psi^*\L'.$ We note that the orientation on $\L'_\u$
thus induced varies continuously with $\u.$ Indeed, given a one parameter family
$\u_t \in B_{\k,l}^{1,p}$ we may choose a homotopy of the one-parameter family
$\Phi_t$ and trivialize $\Phi_t^*\L'.$

We turn to orienting $\L'|_{B'}.$ The relative $Pin$ structure of $(X,L)$
provides a vector bundle $V \rightarrow X^{(3)}$ and a $Pin^\pm$ structure on
$V|_{L^{(3)}}\oplus TL|_{L^{(3)}}.$ We introduce the shorthand notation
\begin{equation}\label{eq:vrc}
V_\R : = V|_{L^{(3)}}, \qquad V_\C : = V \otimes \C.
\end{equation}
Again, it suffices to canonically orient each individual line $\L'_\u$ for $\u
\in B'$ in a way that varies continuously in families. Let $D_0$ be an
arbitrary Cauchy-Riemann operator on $u^*V\otimes\C.$  We consider the
operator $D_\u\oplus D_0,$
\[
\xymatrix@-20pt{   TB_{\k,l}^{1,p} \oplus
W^{1,p}(u^*V_\C,u^*V_\R) \ar[rr]^{D_\u \oplus D_0} \ar@{=}[d]&& \E
\oplus L^p(\Omega^{0,1}(u^*V_\C)) \ar@{=}[d].
\\
W^{1,p}(u^*(TX \oplus V_\C),u|_{\partial\Sigma}^*(TL \oplus V_\R))\oplus
\R^{|\k|}\oplus \C^l &&
L^p(\Omega^{0,1}(u^*(TX\oplus V_\C))).
}
\]
Clearly, there exists a short exact sequence of Fredholm operators
\[
0 \rightarrow D_\u \rightarrow D_\u\oplus D_0 \rightarrow D_0 \rightarrow 0.
\]
So, by Lemma \ref{lem:ses} there exists a natural isomorphism
\begin{equation*}
\det(D_\u) \stackrel{\sim}{\rightarrow} \det(D_\u \oplus D_0)\otimes\det(D_0)^*,
\end{equation*}
and, after tensoring on both sides by $\bigotimes_{a,i}
\left(evb_{ai}^*\det(TL)^*\right)_\u,$
\begin{align*}
\L'_\u &=
\det(D_\u)\otimes\bigotimes_{a,i}\left(evb_{ai}^*\det(TL)^*\right)_\u \\
&\simeq \det(D_\u \oplus D_0)\otimes\det(D_0)^*\otimes\bigotimes_{a,i}
\left(evb_{ai}^*\det(TL)^*\right)_\u.
\end{align*}
Choose a $Pin$ structure $\tilde\p_0$ on $u^*V \rightarrow \Sigma$ and define
$\p_0$ to be its restriction to $u^*V_\R \rightarrow \partial\Sigma.$
By Lemma \ref{lem:otb}, the canonical orientation that the Cauchy-Riemann $Pin$
boundary value problem
\[
\underline D_0 = (\Sigma, u^*V_\C,u^*V_\R,\p_0,D_0)
\]
induces on $\det(D_0)$ does not depend on the choice of $\tilde \p_0.$ So, it
suffices to orient
\begin{equation} \label{eq:so}
\L'_\u \otimes \det(D_0) \simeq \det(D_\u\oplus
D_0)\otimes\bigotimes_{a,i}\left(evb_{ai}^*\det(TL)^*\right)_\u.
\end{equation}
Note that by pull back, the relative $Pin$ structure on $L$ gives a $Pin$
structure on $u|_{\partial\Sigma}^*(TL\oplus V_\R),$ the boundary condition
for $D_\u \oplus D_0.$   If $L$ is orientable and given an
orientation, since by definition $V$ has an orientation, we have an induced
orientation on $u|_{\partial\Sigma}^*(TL\oplus V_\R).$ So, by Proposition
\ref{prop:orl}, we have a canonical orientation on $\det(D_\u
\oplus D_0).$ Since the orientation of $L$ is equivalent to an orientation of
$\det(TL),$ we have given everything on the right-hand side of equation
\eqref{eq:so} a canonical orientation.

If $L$ is not orientable, choose an arbitrary orientation on $(evb_{a1}^*TL)_\u$
for each $a$ such that $k_a \neq 0.$ The complex structure on $\Sigma$ induces a
natural orientation on $\Sigma$ and hence on $(\partial\Sigma)_a$ for each $a.$
For each $a$ and each $i \in [2,k_a],$ trivializing
$u|_{\partial\Sigma}^*TL$ along the oriented line segment in
$(\partial\Sigma)_a$ from $z_{a1}$ to $z_{ai}$ induces an orientation on
$(evb_{ai}^*TL)_\u.$ In the case that $u|_{(\partial\Sigma)_a}^*TL$ is
orientable, the choice of orientation on $(evb_{a1}^*TL)_\u$ induces an
orientation on $u|_{(\partial\Sigma)_a}^*TL.$  By Proposition \ref{prop:orl}, an
orientation on $u|_{(\partial\Sigma)_a}^*TL$ if orientable together with
the chosen orientation on $V_\R$ and the previously mentioned $Pin$
structure induces a canonical orientation on $\det(D_\u\oplus D_0)$ and hence
the whole right-hand side of \eqref{eq:so} is oriented after these choices.

Note that changing the orientation on $(evb_{a1}^*TL)_\u$ will change
all the orientations it induces. By Lemma \ref{lem:co}, the condition $k_a \cong
w_1(d_a) + 1 \mod 2$ now implies that changing the orientation on any
given $(evb_{a1}^*TL)_\u$ would make no difference because the total number of
ensuing orientation changes would be even. Finally, the choice of
$D_0$ is irrelevant because the space of real-linear Cauchy-Riemann operators
on $u^*V$ is contractible. Since the argument for orienting $\L'_\u$ applies
word for word for a one-parameter family, we have indeed canonically oriented
$\L'|_{B'}.$
\end{proof}

At this point, we will modify the canonical isomorphism
\[
\L \stackrel{\sim}{\longrightarrow} \bigotimes_{a,i} evb_{ai}^*\det(TL)
\]
we constructed in the proof of Proposition \ref{prop:or}. Indeed,
$B^{1,p}_{\k,l}(L,\Sigma,\d)$ consists of many connected components, at least
one for each ordering of the marked points on their respective boundary
components. Let
\[
{\bf\varpi} = (\varpi_1,\ldots,\varpi_m)
\]
where $\varpi_a$ is a permutation of the integers $1,\ldots,k_a.$ Define
\[
\sign({\bf \varpi}) : = \sum_a \sign(\varpi_a).
\]
Let
\[
B^{1,p}_{\k,l,{\bf \varpi}}(L,\Sigma,\d)
\]
denote the component of $B^{1,p}_{\k,l}(L,\Sigma,\d)$ where the boundary marked
points $(z_{ai})$ are ordered within $\partial\Sigma$ by the permutations ${\bf
\varpi}.$
\begin{dfn}\label{dfn:tci}
When $\dim L \cong 0 \pmod 2$ we define \emph{the canonical isomorphism}
\[
\L \stackrel{\sim}{\longrightarrow} \bigotimes_{a,i} evb_{ai}^*\det(TL)
\]
to be the isomorphism constructed in the proof of Proposition
\ref{prop:or} twisted by $(-1)^{\sign({\bf\varpi})}$ over the component of the
moduli space $B^{1,p}_{\k,l,{\bf \varpi}}(L,\Sigma,\d).$ If $\dim L \cong 1
\pmod 2,$ then we define the canonical isomorphism to be simply the isomorphism
constructed in the proof of Proposition \ref{prop:or}.
\end{dfn}

Now, we move on to orienting moduli spaces of stable maps. We will restrict
attention to stable maps of two components, one of which is the original
Riemann surface $\Sigma,$ and the other of which is a disk bubble. This
will suffice for the purposes of this paper. However, it is not hard to extend
the results below to stable maps of arbitrarily many components of arbitrary
topological type.

We consider the case that a disk bubbles off the boundary component
$(\partial\Sigma)_b$ along with $k''$ of the marked points on
$(\partial\Sigma)_b$
and $l''$ of the interior marked points. Let
\[
k' := k_b - k'',\qquad \k' := (k_1,\ldots,k',\ldots,k_m).
\]
Let $l' + l'' = l.$ It will turn out to be convenient to keep track of
exactly which marked points bubble off. So, let $\sigma \subset [1,k_b]$
denote the subset of boundary marked points that bubble off and let
$\sigma^c$ denote its complement. Furthermore, let $\varrho \subset [1,l]$
denote the set of interior marked points that bubble off and let $\varrho^c$
denote its complement. Let
\[
\d' = (d',d_1,\ldots,d_b',\ldots,d_m) \in H_2(X,L)\oplus H_1(L)^{\oplus m},
\qquad d'' \in H_2(X,L),
\]
satisfy
\[
d' + d'' = d, \qquad d_b' + \partial d'' = d_b.
\]
We will need to add an extra marked point to each of the two irreducible
components of the stable map in order to impose the condition that the two
components intersect. We denote by $z_0'$ the extra marked point on $\Sigma$ and
by $z_0''$ the extra marked point on $D^2.$ We will use
the notation
\[
\k' + e_b = (k_1,\ldots,k'+1,\ldots,k_m).
\]
We define the space of $W^{1,p}$ stable-maps with this combinatorial data to be
the fiber product
\[
B_{\k,\sigma,l,\varrho}^{1,p}(L,\Sigma,\d',d'') : =
B_{\k'+e_b,l'}^{1,p}(L,\Sigma,\d')\, _{evb_{0}'}\!\times_{evb_0''}
B_{k''+1,l''}^{1,p}(L,D^2,d'').
\]
Elements $\u \in B_{\k,\sigma,l,\varrho}^{1,p}(L,\Sigma,\d',d'')$ take the form
\begin{gather*}
\u = (\u',\u''),\qquad \u'\in \: B_{\k'+e_b,l'}^{1,p}(L,\Sigma,\d'),\qquad \u''
\in B_{k''+1,l''}^{1,p}(L,D^2,d''), \\
evb_{0}'(\u') = evb_0''(\u'').
\end{gather*}
Associated to each such $\u$ there is a nodal Riemann surface with boundary
\[
\hat \Sigma_\u := \Sigma \cup D^2 /z_{0}' \sim z_0''.
\]
and a continuous map
\[
u:  (\hat\Sigma_\u,\partial\hat\Sigma_\u) \rightarrow (X,L)
\]
given by $u'$ on $\Sigma$ and $u''$ on $D^2.$ We denote the node of $\hat
\Sigma_\u$ by $z_0.$ Let
\begin{align*}
p' &:  B_{\k,\sigma,l,\varrho}^{1,p}(L,\Sigma,\d',d'') \rightarrow
B_{\k'+e_b,l'}^{1,p}(L,\Sigma,\d'),  \\
p'' &:  B_{\k,\sigma,l,\varrho}^{1,p}(L,\Sigma,\d',d'')\rightarrow
B_{k''+1,l''}^{1,p}(L,D^2,d''),
\end{align*}
denote the natural projections. Note that when various indices are clear from
the context, we may abbreviate
\begin{gather*}
B' = B_{\k'+e_b,l'}^{1,p}(L,\Sigma,\d'), \qquad
B'' = B_{k''+1,l''}^{1,p}(L,D^2,d''), \\
B^\# : = B_{\k,\sigma,l,\varrho}^{1,p}(L,\Sigma,\d',d'').
\end{gather*}
Define the Banach space bundle $\E^\# \rightarrow
B^\#$ by
\begin{equation*}
\E^\# := p_1^*\E' \oplus p_2^* \E''.
\end{equation*}
Fiberwise, we have
\begin{equation*}
\E_\u^\# := L_p(\Sigma, \Omega^{0,1}(u'^*TX)\oplus\Omega^{0,1}(u''^*TX))
\end{equation*}
for $\u \in B_{\k,\sigma,l,\varrho}^{1,p}(L,\Sigma,\d',d'').$ If $J \in
\J_\omega,$ and $\nu
\in \PP,$ we let
\begin{equation*}
\overline\partial_{(J,\nu)}^\#: B^\# \rightarrow \E^\#
\end{equation*}
denote the section of $\E$ given by the $\nu$-perturbed Cauchy
Riemann operator. Here, the natural $\nu$-perturbed
Cauchy-Riemann operator has a vanishing inhomogeneous term on the disk bubble.
Using the canonical identification between the vertical
tangent spaces of $\E^\#$ and $\E^\#$ itself we define
\begin{equation*}
D^\# := D\overline\partial_{(J,\nu)}^\#:
TB_{\k,\sigma,l,\varrho}^{1,p}(L,\Sigma,\d',d'') \rightarrow \E^\#
\end{equation*}
to be the vertical component of the linearization of
$\overline\partial_{(J,\nu)}^\#.$ Finally, define $\L^\# \rightarrow
B^{1,p}_{\k,l}(L,\Sigma,\d)$ to be the determinant line bundle of the
family of Fredholm operators $D^\#,$
\begin{equation*}
\L^\# := \det(D^\#).
\end{equation*}
Again, suppose either $L$ is orientable and provided with an orientation or $k_a
\cong w_1(d_a) + 1 \mod 2.$
\begin{prop}\label{prop:ors}
The combination of an orientation of $L$ if $L$ is orientable and the choice of
relative $Pin^\pm$ structure $\p$ on $L$ canonically determine an isomorphism
of line bundles
\begin{equation*}
\L^\# \stackrel{\sim}{\longrightarrow} \bigotimes_{a,i} evb_{ai}^*\det(TL).
\end{equation*}
\end{prop}
\begin{proof}
The proof is the same as the proof of Proposition \ref{prop:or} except for one
extra subtlety, which is particularly important in the case that $L$ is not
orientable. That is, the Riemann surface underlying the Cauchy-Riemann $Pin$
boundary problem associated to $D^\#$ is singular. As in equation \eqref{eq:so}
of the proof of Proposition \ref{prop:or}, it
suffices to orient
\begin{equation}\label{eq:sos}
\L^{\#'}_\u \otimes\det(D^\#_0) \simeq \det(D^\#_\u \oplus D^\#_0)
\otimes \bigotimes_{a,i}\left(evb_{ai}^*\det(TL)^*\right)_\u
\end{equation}
for $\u$ such that
\[
u : (\hat\Sigma_\u,\partial\hat \Sigma_\u) \rightarrow (X^{(3)},L^{(3)}).
\]
To this end, we need to describe $D^\#_0$ and $D^\#$ in greater detail.
Using the notation \eqref{eq:vrc}, define
\[
d_{00}^0 : W^{1,p}(u'^*V_\C,u'^*V_\R) \oplus W^{1,p}(u''^*V_\C,u''^*V_\R)
\rightarrow (ev_{b0}^*V_\R)_\u
\]
by
\[
d_{00}^0(\xi',\xi'') = \xi'(z'_0) - \xi''(z''_0).
\]
Denote
\[
W^{1,p}\left(u^*V_\C,u|_{\partial\hat\Sigma_\u}^*V_\R\right) : = \ker(d_{00}^0).
\]
At this point, we introduce abbreviated notation
\begin{gather*}
W'_V := W^{1,p}(u'^*V_\C,u'^*V_\R),\qquad W_V'': =
W^{1,p}(u''^*V_\C,u''^*V_\R), \\
W_V: = W^{1,p}\left(u^*V_\C,u|_{\partial\hat\Sigma_\u}^*V_\R\right), \\
Y_V' : = L^p(u'^*V_\C),\qquad Y_V'' := L^p(u''^*V_\C), \\
Y_V : = L^p(u^*V_\C) = L^p(u'^*V_\C) + L^p(u'^*V_\C).
\end{gather*}
Then, choose arbitrary Cauchy-Riemann operators
\[
D_0' : W'_V \rightarrow Y'_V, \qquad D_0'' : W''_V \rightarrow Y''_V,
\]
and define
\[
D_0^\# := \left ( D_0'\oplus D_0''\right)|_{W_V} : W_V \rightarrow Y_V.
\]
Now, turning to $D^\#,$ define
\begin{equation*}
d_{00} : p'^*TB_{\k',l'}^{1,p}\oplus p''^*TB_{k'',l''}^{1,p} \rightarrow
evb_{0}^* TL
\end{equation*}
by
\[
(\xi',\xi'') \rightsquigarrow devb_{0}'(\xi') - dev_0''(\xi'').
\]
Note that
\begin{equation*}
TB_{\k,\sigma,l,\varrho}^{1,p}(L,\Sigma,\d',d'') = \ker(d_{00}).
\end{equation*}
So, we have a short exact sequence of families of Fredholm operators
\[
\xymatrix{\E_\u^\#\oplus Y_V \ar[r] & p'^*\E_\u'\oplus Y_V'\oplus
p''^*\E_\u''\oplus Y_V'' \ar[r] & 0 \\
TB_\u^\# \oplus W_V \ar[r] \ar[u]^{D_\u^\#\oplus D_0^\#} &
p'^*TB_\u'\oplus W_V'\oplus p''^*TB_\u''\oplus W_V''
\ar[r]^(.62){d_{00}\oplus d_{00}^0}\ar[u]^{D'_\u\oplus D_0'\oplus D_\u''\oplus
D_0''} &
evb_{0}^*\left(TL\oplus V_\R \right)_\u,\ar[u]
}
\]
and, hence, an isomorphism
\begin{equation}\label{eq:sosa}
\det(D_\u^\#\oplus D^\#_0) \stackrel{\sim}{\rightarrow} \det(D_\u'\oplus D_0')
\otimes \det(D_\u''\oplus D_0'')\otimes
evb_{0}^*\det(TL\oplus V_\R)_\u^*.
\end{equation}
Noting isomorphism \eqref{eq:sosa}, if $L$ is oriented, the whole right-hand
side of equation \eqref{eq:sos} is canonically oriented by arguing as in the
proof of Proposition \ref{prop:or}.

If $L$ is not orientable, choose an arbitrary orientation on
$(evb_{a1}^*TL)_\u$ for each $a$ such that $k_a \neq 0.$ The complex
structure on $\Sigma$ induces a natural orientation on $\Sigma$ and hence on
$(\partial\Sigma)_a$ for each $a.$ Similarly, the complex structure on $D^2$
induces a natural orientation on $\partial D^2.$ This said, any ordered pair of
points $(z,z'),\, z \neq z_0,$ in the same connected component of
$\partial\hat\Sigma$ can be connected by a unique oriented line segment from the
first to the second. For the non-singular boundary components of
$\partial\hat\Sigma,$ this is evident. For the singular boundary component,
\[
(\partial\hat\Sigma)_b : = (\partial \Sigma)_b \cup \partial D^2 / z_{0}' \sim
z_0'',
\]
we define the unique oriented line segment from $z$ to $z'$ as follows:
For concreteness, assume that $z \in (\partial\Sigma)_b.$ The same exact
definition applies if $z\in \partial D^2.$ Start from $z$ and proceed
in the direction of the orientation of $(\partial\Sigma)_b$ until reaching
either $z'$ or $z_0'.$ If $z'$ is reached first or if $z' = z_0',$ the path ends
there. If $z_0'$ is reached first, then continue the path starting from
$z_0''\in \partial D^2$ and proceeding along $\partial D^2$ in the direction of
the orientation. If $z'$ belongs to $\partial D^2$ then the path ends when it
reaches $z'.$ Otherwise it continues around $\partial D^2$ back to $z_0''$ and
then proceeds in the direction of the orientation from $z_0'$ along
$(\partial\Sigma)_b$ until it reaches $z'.$

For each $a$ and each $i \in [2,k_a],$ trivializing
$u|_{\partial\hat\Sigma}^*TL$ along the oriented line segment in
$(\partial\hat \Sigma)_a$ from $z_{a1}$ to $z_{ai}$ induces an orientation on
$(evb_{ai}^*TL)_\u$ as before. To orient the factor of $evb_{0}^*\det(TL\oplus
V_\R)_\u^*$ appearing on the right-hand side of isomorphism \eqref{eq:sosa},
since $V_\R$ is equipped with a chosen orientation by the definition of a
relative $Pin$ structure, it suffices to orient $(evb_{0}^*TL)_\u.$ For this
purpose, we proceed as follows: If $k_b \neq 0$ then trivializing
$u|_{\partial\hat\Sigma}^*TL$ along the oriented line segment from $z_{b1}$ to
$z_0$ induces an orientation on $(evb_{0}^*TL)_\u.$ If $k_b = 0$ then choose an
arbitrary orientation on $(evb_{0}^*TL)_\u.$ For all $a \neq b,$ if
$u|_{(\partial\hat\Sigma)_a}^*TL$ is orientable, the choice of orientation
on $(evb_{a1}^*TL)_\u$ induces an orientation on
$u|_{(\partial\Sigma)_a}^*TL$ since in this case $k_a \neq
0$ by the assumption $w_1(d_a) \cong k_a + 1 \pmod 2.$
Furthermore, we may induce an orientation on either or both of
$u|_{(\partial\Sigma)_b}^*TL$ and $u|_{\partial D^2}^*TL$ if orientable, from
the orientation on $(evb_{0}^*TL)_\u.$
By Proposition \ref{prop:orl}, the
orientation on $u|_{(\partial\Sigma)_a}^*TL$ if orientable together with
the chosen orientation on $V_\R$ and the previously mentioned $Pin$
structure induce a canonical orientation on $\det(D_\u'\oplus D_0').$
Similarly, the
orientation on $u|_{D^2}^*TL$ if orientable together with
the chosen orientation on $V_\R$ and the previously mentioned $Pin$
structure induce a canonical orientation on $\det(D_\u''\oplus D_0'').$
Hence the whole right-hand side of \eqref{eq:sosa} is oriented after these
choices. Since we have also chosen an orientation on $(evb_{ai}^*TL)_\u,$ it
follows that the entire right-hand side of \eqref{eq:sos} is oriented.

Note that changing the orientation on $(evb_{a1}^*TL)_\u$ will change
all the orientations it induces. Similarly, in the case $k_b = 0,$ changing
the chosen orientation on $(evb_{0}^*TL)_\u$ changes any orientation it induces.
By Lemma \ref{lem:co}, the condition $k_a \cong
w_1(d_a) + 1 \mod 2$ now implies that changing the orientation on any
given $(evb_{a1}^*TL)_\u$ or $(evb_{0}^*TL)_\u$ would make no difference because
the total number of ensuing orientation changes would be even.
\end{proof}

Now, we modify the canonical isomorphism
\[
\L^\# \stackrel{\sim}{\longrightarrow} \bigotimes_{a,i} evb_{ai}^*\det(TL)
\]
constructed in the proof of Proposition \ref{prop:ors} to be consistent with
the isomorphism over the moduli space of irreducible stable maps defined in
Definition \ref{dfn:tci}. For this purpose, we note that there is a canonical
ordering on the marked points in the boundary of the nodal curve
$\hat\Sigma_\u.$ Indeed, recall from the proof of Proposition \ref{prop:ors}
that given any pair of points $(z,z') \subset \partial\hat\Sigma$ such that $z
\neq z_0,$ there exists a unique oriented line segment from $z$ to $z'.$ We
define $z_{ai}$ to be ordered before $z_{ai'}$ if and only if the segment from
$z_{a1}$ to $z_{ai}$ lies within the segment from $z_{a1}$ to $z_{ai'}.$ Again,
we can divide the moduli space $B^\#$ into components $B^\#_{\varpi}.$
By analogy with the irreducible case, we make the following definition.
\begin{dfn}\label{dfn:tcis}
If $\dim L = 0 \pmod 2,$ we define \emph{the canonical isomorphism}
\[
\L^\# \stackrel{\sim}{\longrightarrow} \bigotimes_{a,i} evb_{ai}^*\det(TL)
\]
to be the isomorphism constructed in the proof of Proposition \ref{prop:ors}
twisted by $(-1)^{\sign(\varpi)}$ over the component $B^\#_{\varpi}.$ If $\dim
L = 1 \pmod 2,$ we define the canonical isomorphism to be the isomorphism
constructed in the proof of Proposition \ref{prop:ors}.
\end{dfn}

\section{The definition of the invariants revisited}\label{sec:def2}
At this point, we will rigorously define the integration carried out in
defining the invariants in Section \ref{sec:def} by generalizing the techniques
of Ruan and Tian \cite{RT}.

First, we will carefully define the inhomogeneous perturbation
to the Cauchy-Riemann equation relevant in the current situation. Let
$\mathcal C$ be a parameter space to be specified later. Let $\pi_i,\, i = 1,2,$
denote the projection from $\Sigma \times X \times \mathcal C$ to the $i^{th}$
factor and let $\pi_i'$ denote the restriction of $\pi_i$ to $\partial\Sigma
\times L \times \mathcal C.$ We define the set of inhomogeneous terms
$\PP$ to be the set of
sections
\[
\nu \in \Gamma\left(\Sigma\times X\times \mathcal
C,Hom\left(\pi_1^*T\Sigma,\pi_2^*T\Sigma\right)
\right)
\]
such that
\begin{enumerate}
\renewcommand{\theenumi}{(\roman{enumi})}
\renewcommand{\labelenumi}{\theenumi}
\item \label{it:1}
$\nu$ is $(j_\Sigma,J)$-anti-linear, i.e. $\nu \circ j_\Sigma = - J \circ \nu;$
\item \label{it:2}
$\nu|_{\partial\Sigma\times L\times \mathcal C}$ carries the sub-bundle
${\pi'_1}^*T\partial\Sigma\subset {\pi'_1}^* T\Sigma$ to the sub-bundle
${\pi'_2}^*(J\,T\!L)\subset {\pi'_2}^*TX.$
\end{enumerate}
For the time being, take $\mathcal C = B_{k,l}^{1,p}(L,\Sigma,\d)$ and define
the section $\bar\partial_{J,\nu}$ of $\E$ by
\[
 \bar\partial_{J,\nu}u := du\circ j_\Sigma + J \circ du - \nu(\cdot,u(\cdot),\u)
\in L^p\left(\Omega^{0,1}(u^*TX)\right).
\]
\begin{lemma}\label{lem:wp}
The operator $\bar\partial_{J,\nu}$ gives rise to an
elliptic boundary value problem.
\end{lemma}
\begin{proof}
Define $\widetilde \Sigma: = \Sigma\cup_{\partial \Sigma}\overline \Sigma.$
Clearly, we may extend any $\nu \in \PP$ over $\widetilde \Sigma$ to
obtain a $\tilde \nu$ that continues to satisfy condition \ref{it:1}. Let
$J_\nu$ be the automorphism of $T(\widetilde\Sigma\times X)$ given in matrix
form by
\[
J_{\tilde\nu} : = \begin{pmatrix}
j_\Sigma & 0 \\
\tilde \nu & J
\end{pmatrix}.
\]
It is not hard to check that condition \ref{it:1} implies that $J_{\tilde \nu}$
is an almost complex structure on $\widetilde \Sigma \times X.$ Condition
\ref{it:2} implies that $\partial \Sigma \times L\subset \widetilde \Sigma
\times X$ is a totally real submanifold. A map
\[
u: (\Sigma,\partial\Sigma) \rightarrow (X,L)
\]
satisfying $\bar\partial_{J,\nu} u = 0$ is equivalent to a
standard $J_{\tilde \nu}$-holomorphic map
\[
\tilde u :  (\Sigma,\partial\Sigma) \rightarrow (\widetilde\Sigma\times
X,\partial\Sigma\times L)
\]
satisfying $\pi_1\circ \tilde u = \id_\Sigma.$ We conclude that
$\bar\partial_{J,\nu}$ does indeed give rise to an elliptic boundary value
problem.
\end{proof}
In light of Lemma \ref{lem:wp}, we define
\[
\widetilde\M_{\k,l}(L,\Sigma,\d) : = \bar
\partial_{J,\nu}^{-1}(0) \subset B_{\k,l}^{1,p}(L,\Sigma,\d).
\]
Let $\phi$ be an anti-symplectic involution of $X$ such that $L = Fix(\phi).$ We
call a map
\[
u : (\Sigma,\partial\Sigma) \rightarrow (X,L)
\]
\emph{$\phi$-multiply covered} if there does not exist $z \in \Sigma$ such that
\[
du(z) \neq 0, \qquad u(z) \notin u(\Sigma\setminus \{z\}), \qquad u(z) \notin
\im \phi\circ u.
\]
Such maps are also commonly termed \emph{not $\phi$-somewhere injective.} A
standard argument shows that the moduli space of $\phi$-somewhere injective maps
has expected dimension for generic $J\in \J_{\omega,\phi}$ even when $\nu = 0.$
See \cite[Section 11]{FOOO}. However, if we take $\nu = 0,$ then the moduli
space $\widetilde\M_{\k,l}(L,\Sigma,\d)$ may be singular at $\phi$-multiply
covered maps even for a generic choice of $J.$ Assume for a moment that $\Sigma
= D^2$ and $\mu(d) > 2.$ Then by the following Lemma, the image of
$\phi$-multiply covered maps under the evaluation map has codimension greater
than or equal to $2.$
\begin{lemma}\label{lem:mc}
Suppose $u : (D^2,\partial D^2) \rightarrow (X,L)$ is $\phi$-multiply covered.
Then, we may factor $u$ as a composition $u = u'\circ \chi,$ where $u'$ is a
real
$J$-holomorphic map
\[
 u' : \C P^1 \rightarrow X, \qquad \phi\circ u' \circ c' = u',
\qquad \bar\partial_J u' =0,
\]
such that $u'|_{D^2}$ is $\phi$-somewhere injective,
and
\[
\chi : (D^2,\partial D^2) \rightarrow (\C P^1,\R P^1)
\]
is a holomorphic map of degree greater than or equal to $2.$
\end{lemma}
\begin{rem}
Since we may extend $\chi$ to a holomorphic map from $\C P^1$ to $\C P^1$
by the Schwarz reflection principle, $\chi$ is actually given by polynomial.
Thinking of $D^2$ as $H \cup \{\infty\}$ and thinking of $\C P^1$ as $\C\cup
\{\infty\},$ the boundary condition on $\chi$ implies that $\chi$ is given by a
real polynomial.
\end{rem}
\begin{proof}[Proof of Lemma \ref{lem:mc}]
Gluing together $u$ and $\phi\circ u\circ c'$ we obtain a real $J$-holomorphic
map
\[
\tilde u :\C P^1 \rightarrow X, \qquad \phi\circ \tilde u \circ c' = \tilde u.
\]
By a standard theorem \cite[Chapter 2]{MS}, we can factor $\tilde u = \tilde
\chi \circ \tilde u'$ where $\tilde u'$ is somewhere injective. Another
standard result \cite[Theorem E.1.2]{MS} says that $\tilde u'$ is injective
except at a finite number of points. Since the image
of $u'$ is clearly invariant under $\phi,$ $\phi$ induces an anti-holomorphic
involution of $\C P^1$ away from a finite number of points. So, removing
singularities, there exists an anti-holomorphic involution $c''$ on $\C P^1$
such that $\phi \circ \tilde u' \circ c'' = \tilde u'.$  Moreover, $c''$ fixes
the image of $\R P^1$ under $\tilde\chi.$ It follows that $c''$ is conjugate to
$c'$
by some biholomorphism $a \in PSL_2(\C).$ So, we may define
\[
u' = \tilde u' \circ a^{-1}, \qquad \chi = a \circ \tilde
\chi|_{D^2}.\tag*{\qedhere}
\]
\end{proof}
Consequently, under the previously mentioned conditions, $\phi$-multiply covered
maps are not important in the definition of intersection theory on the moduli
space. So, even in the case $\nu = 0,$ we may obtain a smooth moduli space by
defining
\begin{equation*}
\widetilde \M_{\k,l}^*(L,\Sigma,\d) : =
\widetilde\M_{\k,l}(L,\Sigma,\d)\setminus \{\text{$\phi$-multiply covered
maps}\}.
\end{equation*}
Then, to obtain an interesting intersection theory, we define
\begin{equation}\label{eq:nu0}
\M_{\k,l}(L,\Sigma,\d) : = \widetilde \M_{\k,l}^*(L,\Sigma,\d)/PSL_2(\R),
\end{equation}
where $PSL_2(\R)$ acts by
\[
(u,\vec z,\vec w) \rightsquigarrow (u\circ\varphi,(\varphi^{-1})^k(\vec
z),(\varphi^{-1})^l(\vec w)), \qquad \varphi \in PSL_2(\R).
\]
This choice of action ensures that the evaluation maps descend to the quotient.
Note that the action would not preserve the moduli space if we were to allow a
general inhomogeneous term $\nu.$

Equivalently, instead of quotienting by $PSL_2(\R),$ we could consider an
appropriate section of the group action. Fortunately, this approach, as observed
in \cite{RT}, generalizes to the situation where we allow a generic $\nu$ and so
provides a definition of the moduli space that works even when $(X,L)$ may admit
holomorphic disks of Maslov index $0$ or when $\Sigma \neq D^2.$

So, we construct a section of the reparameterization group action in the
following manner. Let
\[
\pi_j : \widetilde\M_{\k,l}(L,\Sigma,\d) \rightarrow \Sigma
\]
be the projection sending $(u,\vec z,\vec w) \rightsquigarrow w_j.$
First, suppose that $\Sigma \simeq D^2.$ Choose an interior point $s_0 \in
\Sigma$ and a line $\ell \subset D^2$ connecting $s_0$ to $\partial \Sigma$ such
that for any pair of points $(w,w') \in \Sigma$ there exists a unique
$\varphi \in
Aut(\Sigma) \simeq PSL_2(\R)$ that satisfies
\[
\varphi(w)  = s_0, \qquad \varphi(w') \in \ell.
\]
For the time being, assume $l \geq 2.$ We require that the dependence of $\nu$
on $\u \in B_{\k,l}^{1,p}$ factors through $\pi_1 \times \pi_2.$ Moreover,
letting $d(\cdot,\cdot)$ denote the distance function on $D^2,$ we
impose on $\nu$ the condition
\begin{equation}\label{eq:nu}
\nu(\cdot,\cdot,(w_1,w_2)) =
\frac{d(w_1,w_2)}{d(w_1',w_2')} \nu(\cdot,\cdot,(w_1',w_2')) .
\end{equation}
In particular, $\nu$ vanishes uniformly in the limit $w_2 \to w_1.$
We define
\begin{equation}\label{eq:M}
\M_{\k,l}(L,\Sigma,\d) : = (\pi_1 \times \pi_2)^{-1}(s_0 \times \ell) \subset
\widetilde \M_{\k,l}(L,\Sigma,\d).
\end{equation}
Standard arguments show that for a generic choice of $\nu$ satisfying
\eqref{eq:nu}, $\M_{\k,l}(L,\Sigma,\d)$ will be a smooth manifold of expected
dimension. We digress briefly to explain the significance of condition
\eqref{eq:nu}. To prove the invariance of $N_{\Sigma,\d,\k,l},$ we need to argue
that stable maps in the Gromov compactification of $\M_{\k,l}(L,\Sigma,\d)$
involving sphere bubbles occur only in codimension two. A sequence $\u^i \in
\M_{\k,l}(L,\Sigma,\d)$ such that $w_2^i \to w_1^i$ will Gromov converge to a
stable map consisting of one disk component and one sphere component with both
$w_1$ and $w_2$ on the sphere component. The nodal point where the sphere and
disk are attached is fixed at $s_0.$ However, there are no other fixed marked
points on the disk component to compensate for the remaining $T^1$ symmetry.
Condition \eqref{eq:nu} implies that the limit inhomogeneous term on the disk
component is zero. So, we can simply quotient by the residual $T^1$ action.
Although $\phi$-multiply covered maps can arise, we may disregard them because
by a standard argument, their image under the evaluation map has high
codimension.

If, $\Sigma \simeq S^1 \times I,$ we choose a line $\ell \subset \Sigma$
connecting the two boundary components of $\Sigma$ such that for any $w \in
\Sigma$ there exists a unique $\varphi \in Aut(\Sigma)\simeq T^1$ such that
$\varphi(w) \in \ell.$ We choose $\nu$ to be entirely independent of
$\u.$ Assuming for the time being that $l \geq 1,$ define
\[
\M_{\k,l}(L,\Sigma,\d) : = \pi_1^{-1}(\ell) \subset \widetilde
\M_{\k,l}(L,\Sigma,\d).
\]
Again,  for a generic choice of $\nu,$ standard arguments show that
$\M_{\k,l}(L,\Sigma,\d)$ will be a smooth manifold of expected dimension.

Now we turn to the case when $\Sigma \simeq D^2$ and $l = 0.$ This case is of
particular interest because it arises when $X$ is a Calabi-Yau manifold and
$L$ is a special Lagrangian submanifold. The cases $\Sigma \simeq D^2,\,l=1,$
and $\Sigma \simeq S^1 \times I,\, l = 0,$ use a very similar argument.

We start with the moduli space $\M_{\k,2}(L,\Sigma,\d)$ constructed above and
proceed as follows.
Choose smooth manifolds $A,B,$ and maps
\[
f: A \rightarrow X, \qquad g: B \rightarrow X
\]
such that $(A,f)$ and $(B,g)$ define pseudo-cycles representing the
Poincare-dual of any non-trivial $\phi$-anti-invariant
$2^{nd}$ cohomology class. The symplectic form $\omega$ provides at
least one example, and for simplicity, we will continue with this example.
\begin{lemma}\label{lem:AB}
We can choose $(A,f)$ and $(B,f)$ that are $\phi$-anti-invariant and do not
intersect $L.$
\end{lemma}
\begin{proof}
A straightforward transversality argument shows that we may assume $(A,f)$ is
transversal to $L.$ If necessary, replacing $A$ by $\frac{1}{2} A \coprod
-\frac{1}{2} A$ and $f$ by $f \coprod \phi\circ f,$ we may assume that the
pseudo-cycle $(A,f)$ is $\phi$-anti-invariant just like $\omega.$ Now, we
consider a local model for $(A,f)$ near $L$ and show how to modify $(A,f)$ near
$L$ to avoid intersecting $L.$ Choose $\phi$-invariant local symplectic
coordinates $\Theta : X \rightarrow \C^n$ near $L$ such that $\Theta(L) =
\R^n.$ We may assume that the image of $f$ is given by the union of the
vanishing sets of conjugate real-linear maps
\[
\ell,\,\bar\ell : \C^n \rightarrow \R^2.
\]
More precisely, choose complex coordinates $z = (z_1,\ldots, z_n)$ on $\C^n.$
Taking $x = (x_1,\ldots,x_n)$ and $y = (y_1,\ldots,y_n),$ we write $z = x
 + iy.$ So, we can decompose
\[
\ell(z) = \ell^x(x) + \ell^y(y), \qquad \bar\ell(z) = \ell^x(x) - \ell^y(y).
\]
Letting $\ell_i,\, i = 1,2,$ denote the $i^{th}$ component of $\ell,$ the
equations for the image of $f$ may be written as
\[
0 = \ell_i \cdot \bar\ell_i = \ell_i^{x^2} - \ell_i^{y^2}, \quad i = 1,2.
\]
So, choosing small constants $\epsilon_i > 0,\, i = 1,2,$ we modify $(A,f)$
so that locally the image of $f$ satisfies equations
\[
\ell_i^{x^2} - \ell_i^{y^2} = - \epsilon_i, \quad i = 1,2.
\]
Clearly, these equations have no real solutions. The same applies for $B.$
\end{proof}
Furthermore, we may assume that $(A\times B,f\times g)$ is
transversal to $evi_1 \times evi_2$ so that we may define
\[
\M_{\k,0}(L,\Sigma,\d) : = \frac{1}{\omega(d)^2} \M_{\k,2}(L,\Sigma,\d)
\times_{X\times X} \left ( A \times B \right).
\]
The factor $\frac{1}{\omega(d)^2}$ in front of the fiber-product means that
each point in the moduli space should be counted with weight
$\frac{1}{\omega(d)^2}.$ This correction is designed to cancel the contribution
of the divisors $(A,f)$ and $(B,g)$ as predicted by the divisor axiom of
formal Gromov-Witten theory. Indeed a map $u: (\Sigma,\partial\Sigma)
\rightarrow (X,L)$ representing the class $d \in H_2(X,L)$ should intersect a
pseudo-cycle Poincare dual to $\omega$ exactly $\omega(d)$ times.
\begin{rem}\label{rem:l>0}
Note that even when $l > 0,$ we are free to fix the group action by adding
divisor constraints just as in the case when $l = 0.$
\end{rem}

\begin{proof}[Proof of Theorem \ref{thm:or}]
Given this definition of $\M_{\k,l}(L,\Sigma,\d),$ Proposition \ref{prop:or}
immediately implies Theorem \ref{thm:or}. We choose the isomorphism according
to Definition \ref{dfn:tci}.
\end{proof}
Now, for a sufficiently generic choice of points $\vec x = (x_{ai}), x_{ai} \in
L,$
and pairs of points $\vec y  = (y_j),$
\[
y_j : \{0,1\} \rightarrow X, \qquad y_j(1) = \phi\circ y_j(0),
\]
the total evaluation map
\[
{\bf ev} : = \prod_{a,i}evb_{ai} \times \prod_j evi_j :
\M_{\k,l}(L,\Sigma,\d) \longrightarrow L^{|\k|} \times X^l
\]
will be transverse to
\[
\prod_{a,i} x_{ai} \times \prod_j y_j : \{0,1\}^l \rightarrow L^{|\k|} \times
X^l.
\]
So, assuming the dimension condition \eqref{eq:dim} is satisfied, we may define
\[
N_{\Sigma,\d,\k,l} : = \# \mathbf{ev}^{-1}(\vec x,\vec y).
\]
Here, $\#$ denotes the signed count with the sign of a given point $v \in
\mathbf{ev}^{-1}(\vec x,\vec y)$ depending on
whether or not the isomorphism
\[
d{\bf ev}_v : \det(T\M_{\k,l}(L,\Sigma,\d))_v \stackrel{\sim}{\rightarrow}
{\bf ev}^* \det\left(T\left(L^{|\k|}\times X^l\right)\right)_v
\]
agrees with the isomorphism of Theorem \ref{thm:or} up to the
action of the multiplicative group of positive real numbers. If the dimension
condition \eqref{eq:dim} is not satisfied, we define $N_{\Sigma,\d,\k,l} := 0.$

\section{The sign of conjugation on the moduli space}\label{sec:conj}
Now, suppose $\Sigma$ is biholomorphic to $\overline \Sigma.$ It follows that
there exists a complex conjugation $c: \Sigma \rightarrow \Sigma.$ Let
$\phi$ be an anti-symplectic involution of $X$ such that
$Fix(\phi) = L.$ Fix $J \in \J_{\omega,\phi},$ define $\PP_{\phi,c}$ to
be the set of $\nu \in \PP$ such that $d\phi\circ\nu\circ dc = \nu$ and let
$\nu \in \PP_{\phi,c}.$ We define an involution
\[
\phi_B: B_{\k,l}^{1,p}(L,\Sigma,\d) \rightarrow B_{\k,l}^{1,p}(L,\Sigma,\d)
\]
by
\[
(u,\vec z,\vec w) \rightsquigarrow \left(\phi \circ u \circ c,
(c|_{\partial D^2})^{|k|}(\vec z), c^l(\vec w)\right).
\]
Furthermore, define an involution $\phi_{\E} : \E
\rightarrow \E$ covering $\phi_B$ by sending $\eta \in \E_u$ to $d\phi \circ
\eta \circ dc \in \E_{\phi_B(u)}.$ It is easy to see that $\overline
\partial_{J,\nu}$ is $\phi_B-\phi_\E$ equivariant, and hence $\phi_\E$ induces
an involution $\phi_\L : \L \rightarrow \L$ covering $\phi_B.$ Now, the
trivial bundle morphism of
\[
\bigotimes_{a,i}evb_{ai}^*\det(TL)
\]
covers $\phi_B.$ So, the involution $\phi_\L$ induces an involution
$\phi_\L'$ of the bundle
\[
\L' : = Hom\left(\bigotimes_{a,i}evb_{ai}^*\det(TL),\L\right) \simeq
\L\otimes\bigotimes_{a,i}evb_{ai}^*\det(TL).
\]
covering $\phi_B.$ Now, the bundle $\L'$ has a canonical orientation
corresponding to the canonical isomorphism of Definition \ref{dfn:tci}. So,
the involution $\phi_\L'$ may either preserve the orientation component of the
complement of the zero section of $\L',$ exchange it with the opposite component
or some combination of the two over different connected components of the base.
If $\phi_\L'$ preserves the orientation of $\L',$ we say it has sign $0.$ If
$\phi_\L'$ does not preserve the orientation of $\L'$ anywhere, we say it has
sign $1.$

We would like to show the sign of $\phi_L$ is well defined and calculate it. To
properly formulate the result, we define a degree $0$ homomorphism
\[
\psi : H_*(X,L;\Z/2\Z) \rightarrow H_*(X;\Z/2\Z)
\]
on the level of singular chains as follows: We implicitly use $\Z/2\Z$
coefficients everywhere. Suppose
\[
\sigma\in C_*(X,L) :=C_*(X)/C_*(L)
\]
is a relative singular chain. Let $\widetilde \sigma
\in C_*(X)$ represent $\sigma.$ Define
\begin{equation*}
\widehat \psi(\sigma) := (\id + \phi_*) \widetilde \sigma.
\end{equation*}
$\widehat \psi$ is well defined because if $\sigma = 0$ then $\widetilde
\sigma \in C_*(L)$ and hence
\begin{equation*}
(\id + \phi_*)\widetilde \sigma = 2\widetilde \sigma \cong 0 \mod 2.
\end{equation*}
Since $\phi_*$ commutes with the boundary operator, so
does $\widehat \psi,$ so we can define $\psi := H(\widehat \psi).$ Now, let
$g_0$ denote the genus of $\Sigma/\partial\Sigma.$
\begin{prop}\label{prop:s}
Let $n = \dim L.$ If $\p$ is a relative $Pin^-$ structure then the involution
$\phi_\L'$ has sign
\begin{align*}
\s^-(\d,\k,l) &\cong \frac{\mu(d)(\mu(d)+1)}{2} + (1-g_0)n + mn + |\k| + l \\
&\qquad + w_2(V)(\psi(d))+ w_1(\partial d) + \sum_{a < b} w_1(d_a) w_1(d_b) \\
&\qquad + \sum_a w_1(d_a)(k_a - 1) + (n+1)\sum_a \frac{(k_a-1)(k_a-2)}{2} \mod
2.
\end{align*}
If $\p$ is a relative $Pin^+$ structure, $\phi_\L'$ has sign
\begin{align*}
\s^+(\d,\k,l) &\cong \frac{\mu(d)(\mu(d)+1)}{2} + (1-g_0)n + mn + |\k| + l\\
&\qquad + w_2(V)(\psi(d)) + \sum_{a < b} w_1(d_a) w_1(d_b) \\
&\qquad + \sum_a w_1(d_a)(k_a - 1) + (n+1)\sum_a \frac{(k_a-1)(k_a-2)}{2}\mod 2.
\end{align*}
\end{prop}
\begin{rem}
Suppose $\dim L \leq 3.$ Then $L$ is $Pin^-$ by the Wu relations, so we can
take $\p$ to be a standard $Pin^-$ structure. In particular, $w_2(V) = 0.$
Since $w_1(\partial d) \cong \mu(d) \mod 2,$ we have
\begin{align*}
\s^-(\d,\k,l) &\cong \frac{\mu(d)(\mu(d)-1)}{2} + (1-g_0)n + mn + |\k| + l \\
&\qquad + \sum_{a < b} w_1(d_a) w_1(d_b) + \sum_a w_1(d_a)(k_a - 1) \\
&\qquad + (n+1)\sum_a \frac{(k_a-1)(k_a-2)}{2} \mod 2.
\end{align*}
In the particularly simple case that $\Sigma \simeq D^2,$ we have
\begin{align*}
\s^-(d,\k,l) &\cong \frac{\mu(d)(\mu(d)-1)}{2} + k + l + \mu(d)(k-1) \\
 & \qquad \mbox{} + (n+1)\frac{(k-1)(k-2)}{2}\mod 2.
\end{align*}
\end{rem}
We omit the proof of the preceding proposition as it is very similar to the
proof of the next proposition and we do not actually use it.

Now, recall from Section \ref{sec:or} that
\[
B_{\k,\sigma,l,\varrho}^{1,p}(L,\Sigma,\d',d'') : =
B_{\k'+e_b,l'}^{1,p}(L,\Sigma,\d')\, _{evb_{0}'}\!\times_{evb_0''}
B_{k''+1,l''}^{1,p}(L,D^2,d'').
\]
Note that $D^2 \simeq \overline{D^2}$ as demonstrated by the standard
conjugation
\[
c: D^2 \rightarrow D^2.
\]
So, we have an involution $\phi_{B''}$
of the second factor of the fiber product. Since $L \subset Fix(\phi),$
the involution $\phi_{B''}$ of the second factor induces an involution
$\phi_{B^\#}$ of the whole fiber product. Similarly, the involution
$\phi_{\E''}$ of the bundle $\E''$ over the second factor of the fiber product
induces an involution $\phi_{\E^\#}$ on the bundle $\E^\# \rightarrow B^\#.$
Recall that the natural inhomogeneous perturbation $\nu$ for stable maps
vanishes on bubble components. In particular, it is $\phi$-invariant. So,
$\bar\partial_{J,\nu}^\#$ defines a $\phi_{B^\#}-\phi_{\E^\#}$ invariant
section of $\E^\#.$ Consequently, $\phi_{B^\#}$ and $\phi_{\E^\#}$ induce an
involution $\phi_{\L^\#}$ of the determinant bundle $\L^\# \rightarrow B^\#.$ As
before, $\phi_{\L^\#}$ induces an involution $\phi_{\L^\#}'$ of
\[
\L^{\#'} : = Hom\left(\bigotimes_{a,i}evb_{ai}^*\det(TL),\L^\#\right) \simeq
\L^\#\otimes\bigotimes_{a,i}evb_{ai}^*\det(TL).
\]
Proving that the sign of $\phi_{\L^\#}'$ is well defined and calculating it will
play a crucial role in the proof of the invariance of $N_{\Sigma,\d,\k,l}.$
Before writing down the formula for the sign of $\phi_{\L^\#},$ let us
introduce some new notation. We define
\begin{equation*}
\Upsilon'(d'',k'') : \cong \mu(d'')k'' \cong w_1(\partial d'') k'' \pmod 2
\end{equation*}
and
\begin{equation*}
\Upsilon''(d'_b,d'',k',k'') : \cong
\begin{cases}
0, & w_1(d'_b) = w_1(\partial d'') = 0 \\
k', & w_1(d'_b) = w_1(\partial d'') = 1 \\
k''-1, & w_1(d'_b) =  1,\; w_1(\partial d'') = 0 \\
k_a - 1 = k'' + k' -1, & w_1(d'_b) = 0,\; w_1(\partial d'') =  1.
\end{cases}
\end{equation*}
\begin{prop}\label{prop:sis}
Let $n = \dim L.$ Suppose the marked point $z_1$ does not bubble off, i.e. $1
\notin \sigma.$ Then the involution $\phi_{\L^\#}'$ of the line-bundle
\[
\L^{\#'} \rightarrow B^{1,p}_{\k,\sigma,l,\varrho}(L,\Sigma,\d',d'')
\]
has sign
\begin{align}
\s^{\#'}_\pm(d'',k'',l'') & : \cong  \frac{\mu(d'')(\mu(d'') \pm 1)}{2}
+ w_2(\psi(d''))+ k'' + 1 + l'' \notag\\
& \qquad  +  \Upsilon'(d'',k'') +
(n+1)\frac{k''(k''-1)}{2} \mod 2,\label{eq:sn'}
\end{align}
with $+$ in the $Pin^+$ and $-$ in the $Pin^-$ case. On the other hand, suppose
now that the marked point $z_1$ does bubble off, i.e. $1 \in \sigma.$ Then the
involution $\phi_{\L^\#}'$ has sign
\begin{align}
\s^{\#''}_\pm(d'_b,d'',k',k'',l'') & : \cong  \frac{\mu(d'')(\mu(d'') \pm 1)}{2}
+ w_2(\psi(d'')) \notag \\
& \quad + k'' + 1 + l'' +  \Upsilon''(d'_b,d'',k',k'') +
w_1(d'_b)w_1(\partial d'')
\notag\\
& \quad + (n+1)\left(\frac{(k''-1)(k''-2)}{2} + k_b k'' + k_b \right)
\mod 2, \label{eq:sn''}
\end{align}
with $+$ in the $Pin^+$ and $-$ in the $Pin^-$ case.
\end{prop}
\begin{rem}
Suppose $L$ is orientable and $\dim L$ is odd. Then, using the fact that
$\mu(d'')$ is even if $L$ is orientable, we have
\begin{equation}\label{eq:soo}
\s^{\#'} = \s^{\#''} = \frac{\mu(d'')}{2}+ w_2(\psi(d'')) + k'' + 1 +
l''.
\end{equation}
\end{rem}
\begin{proof}[Proof of Proposition \ref{prop:sis}:]
The first term in $\s^{\#}$ comes from the formula \eqref{eq:asfd}. This
accounts for the sign of conjugation on the moduli space of unmarked disks. The
terms $k'' + l'' + 1$ account for conjugation on the configuration space of the
marked points, adding one extra point for the incidence condition. Recall from
the proof of Proposition \ref{prop:ors} that the unique oriented path from $z
\neq z_0$ to $z'$ in the boundary of $\partial\hat \Sigma_\u$ played an
important role in determining the canonical orientation of $\L^{\#'}.$ This path
depended on the orientation of $\partial\hat\Sigma,$ which is reversed under
conjugation. The terms $\Upsilon'$ (resp. $\Upsilon'' +
w_1(d'_b)w_1(\partial d'')$) in $\s^{\#'}$ (resp. $\s^ {\#''}$) account for this
dependence. The remaining terms account for the reordering of the marked points
under conjugation, which plays a role only in even dimensions according to
Definition \ref{dfn:tcis}.

We now explain in more detail how the unique oriented path from $z \neq 0$ to
$z'$ changes under conjugation, and how that effects the orientation of
$\L'^\#.$ First, suppose $1 \notin \sigma.$ The path from $z_{b1}$ to $z_{bi}$
for $i \in \sigma$ will change when the orientation of the boundary of the
bubble $\partial D^2$ changes under conjugation. If $w_1(\partial d'') = 1,$
this change of path changes the orientation of $(evb_{ai}^*TL)_\u$ for each $i
\in \sigma.$ Since, $|\sigma| = k'',$ we obtain the total change of
orientation given by $\Upsilon'.$ The explanation of $\Upsilon''$ is similar.

The additional orientation change $w_1(d'_b)w_1(\partial d''_b)$ when $1 \in
\sigma$ enters because then the path from $z_{b1}$ to $z_0$ changes under
conjugation. This path effects the orientation of $(evb_0^*TL)_\u$ when
$w_1(\partial d''_b) = 1.$ The orientation of $(evb_0^*TL)_\u$ enters twice into
the orientation of $\L'^\#$ when $w_1(\partial d_b') = 0.$ Indeed, it determines
the orientation of $evb_{0}^*\det(TL\oplus V_\R)_\u^*,$  and it determines the
orientation of $u|_{(\partial\Sigma)_b}^*TL,$ which in turn determines the
orientation of $\det(D_\u'\oplus D_0').$ Both of these determinants appear on
the right-hand side of isomorphism \eqref{eq:sosa}. However, when $w_1(d'_b) =
1$ we cannot orient $u|_{(\partial\Sigma)_b}^*TL.$ So, the orientation of
$(evb_0^*TL)_\u$ enters only once into the orientation of $\L'^\#.$ So, there
is an extra contribution to the orientation change of $\L'^\#$ exactly when
\[
w_1(d'_b)w_1(\partial d''_b) = 1.\tag*{\qedhere}
\]
\end{proof}

\section{Proof of invariance}\label{sec:inv}
\begin{proof}[Proof of Theorem \ref{thm:inv}]
In order to prove independence from various choices, we construct
cobordisms from parameterized moduli spaces. Complications arise
in compactifying these cobordisms. For concreteness, we focus on
independence of a variation of the constraints on marked points. The proof of
independence of a variation of $J,\,\Sigma,$ or $A,B,$ is very similar.

Recall that the definition of $N_{\Sigma,\d,\k,l}$ depends on the
choice of points $\vec x = (x_{ai}),\, x_{ai} \in L,$ and pairs of
points, $\vec y = (y_j),$
\[
y_j : \{0,1\} \rightarrow X, \qquad y_j(1) = \phi(y_j(0)).
\]
Suppose we choose different points $\vec x'$ and $\vec y'$ satisfying the same
conditions.
This corresponds to changing the forms $\alpha_{ai}$ and $\gamma_j$ mentioned in
Theorem \ref{thm:inv}. Let
\begin{gather*}
\mathbf{x} : [0,1] \rightarrow L^{|\k|}, \qquad \mathbf{x}(0) = \vec x,
\qquad \mathbf{x}(1) = \vec x', \\
\mathbf{y} : [0,1] \times \{0,1\}^l \rightarrow X^l, \qquad \mathbf y(t,1) =
\phi(\mathbf y(t,0)), \\
\mathbf y(0,\star) = \vec y(\star), \qquad \mathbf y(1,\star) = \vec
y'(\star).
\end{gather*}
If we choose $\mathbf x$ and $\mathbf y$ generically, they we will be transverse
to the total evaluation map
\[
\mathbf{ev}: \M_{\k,l}(L,\Sigma,\d) \rightarrow L^{|\k|} \times X^l.
\]
So,
\begin{equation}\label{eq:W}
\W := \W(\mathbf x,\mathbf y) := \M_{\k,l}(\Sigma,L,\d)\:
_{\mathbf{ev}}\!\!\times_{(\mathbf{x}\times\mathbf{y})\circ\Delta}
\left([0,1]\times\{0,1\}^l \right)
\end{equation}
gives a smooth oriented cobordism between
\[
\mathbf{ev}^{-1}(\vec x,\vec y) \quad \text{and} \quad \mathbf{ev}^{-1}(\vec
x',\vec y').
\]
However, $\W$ is generally not compact. So, in order to show the invariance of
$N_{\Sigma,\d,\k,l}$ we must study the non-trivial stable
maps arising in the Gromov-compactification of $\W,$ which we denote by
$\partial_G\W.$

We now digress for a moment to describe $\partial_G\W$ more explicitly.
We define
\[
\widetilde\M_{\k,\sigma,l,\varrho}(L,\Sigma,\d',d'') : = \bar
\partial_{J,\nu}^{\#-1}(0) \subset
B_{\k,\sigma,l,\varrho}^{1,p}(L,\Sigma,\d',d'').
\]
Recall that the inhomogeneous perturbation $\nu$ vanishes on the bubble $D^2
\subset \hat\Sigma.$  This means that the moduli space
$\widetilde\M_{\k,\sigma,l,\varrho}(L,\Sigma,\d',d'')$ may be singular at
$\phi$-multiply covered maps. As long as $\mu(d'')>0,$ this does not present a
problem because Lemma \ref{lem:mc} then shows that the image under the
evaluation map of stable maps with $\phi$-multiply covered bubbled components
has codimension at least two. Also, constant holomorphic disks have expected
dimension. On the other hand, in the case that $(X,L)$ admits holomorphic disks
of positive energy of Maslov index zero, we must take $\phi$-multiply covered
maps into consideration. We postpone the argument in this case to Section
\ref{sec:km}.

We continue now with the description of $\partial_G\W.$ Since $\nu = 0$ on
bubble components, we have an action of $PSL_2(\R)$ on $\widetilde
\M_{\k,\sigma,l,\varrho}(L,\Sigma,\d',d'')$ given by
\[
(\u',(u'',\vec z'',\vec w'')) \rightsquigarrow
\left(\u',\left(u''\circ\varphi,(\varphi^{-1})^{k''}(\vec z''),
(\varphi^{-1})^l(\vec
w'')\right)\right), \quad\! \varphi \in PSL_2(\R).
\]
On the other hand, a generic perturbation term $\nu$ will break the
$Aut(\Sigma)$ invariance of $\bar\partial_{J,\nu}^\#.$
So, we construct a section of the $Aut(\Sigma)$ action that would exist if $\nu$
vanished in the following manner. Let
\[
\pi_j : \widetilde\M_{\k,\sigma,l,\varrho}(L,\Sigma,\d',d'') \rightarrow \Sigma
\]
be the projection sending $(u,\vec z,\vec w) \rightsquigarrow w_j.$ In
addition, we define
\[
\pi_0' : \widetilde\M_{\k,\sigma,l,\varrho}(L,\Sigma,\d',d'') \rightarrow \Sigma
\]
by $\pi_0'(\u) = z_0',$ the point where the bubble attaches.

First, suppose that $\Sigma \simeq D^2$ and $l\geq 2.$ Recall from the
construction of $\M_{\k,l}(L,\Sigma,\d)$ in Section \ref{sec:def2} that we chose
an interior point $s_0 \in \Sigma$ and a line $\ell \subset D^2$ connecting
$s_0$ to $\partial \Sigma.$ We imposed the conditions $w_1 = s_0$ and $w_2
\in \ell.$ Since $s_0$ is an interior point, $w_1$ cannot possibly bubble off
onto a disk bubble, i.e. $1 \notin \varrho.$ However, $w_2$ could bubble onto a
disk bubble that bubbles off at $\ell\cap \partial\Sigma.$ So, we consider the
following two cases. If $w_2$ does not bubble off, i.e. $2 \notin \varrho,$
define
\begin{align*}
\M_{\k,\sigma,l,\varrho}(L,\Sigma,\d',d'') &: = (\pi_1 \times \pi_2)^{-1}(s_0
\times \ell)/PSL_2(\R) \\
&\qquad \subset \widetilde \M_{\k,\sigma,l,\varrho}(L,\Sigma,\d',d'')/PSL_2(\R).
\end{align*}
If $w_2$ does bubble off, that is, $2 \in \varrho,$ we define
\[
\M_{\k,\sigma,l,\varrho}(L,\Sigma,\d',d'') : = (\pi_1 \times \pi_0')^{-1}(s_0
\times (\ell\cap\partial\Sigma))/PSL_2(\R).
\]

Now, suppose $\Sigma \simeq S^1 \times I$ and $l \geq 1.$ Recall from the
construction of $\M_{\k,l}(L,\Sigma,\d)$ that we chose a line $\ell \subset
\Sigma$ connecting the two boundary components of $\Sigma$ and imposed the
condition $w_1 \in \ell.$ So, $w_1$ could bubble off at a
disk connecting to $\ell \cap \partial\Sigma.$ So, we consider the
following two cases. If $w_1$ does not bubble off, i.e. $1 \notin\varrho,$
define
\[
\M_{\k,\sigma,l,\varrho}(L,\Sigma,\d',d'') : = \pi_1^{-1}(\ell)/PSL_2(\R)
\subset \widetilde
\M_{\k,\sigma,l,\varrho}(L,\Sigma,\d',d'')/PSL_2(\R).
\]
If $w_1$ does bubble off, i.e. $1 \in\varrho,$ define
\[
\M_{\k,\sigma,l,\varrho}(L,\Sigma,\d',d'') : =
\pi_0'^{-1}(\ell\cap\partial\Sigma)/PSL_2(\R).
\]

Now we turn to the case when $\Sigma \simeq D^2$ and $l = 0.$ The cases
$\Sigma \simeq D^2,\,l=1,$ and $\Sigma \simeq S^1 \times I,\, l = 0,$ use a very
similar argument, which we omit. Recall from the construction of
$\M_{\k,l}(L,\Sigma,\d)$ that we chose $\phi$-anti-invariant pseudo-cycles
$(A,f)$ and $(B,g)$ representing the Poincare dual of $\omega$ and
satisfying various transversality conditions. Taking
$\M_{\k,\sigma,2,\varrho}(L,\Sigma,\d',d'')$ as defined above, we may perturb
$(A,f)$ and $(B,g)$ slightly so that the evaluation map
\[
ev_1 \times ev_2 : \M_{\k,\sigma,2,\varrho}(L,\Sigma,\d',d'') \rightarrow
X\times X
\]
is transversal to $(A\times B, f\times g).$ So, we may define
\[
\M_{\k,\sigma,0,\varrho}(L,\Sigma,\d',d'') : = \frac{1}{\omega(d)^2}
\M_{\k,\sigma,2,\varrho}(L,\Sigma,\d',d'')
\times_{X\times X} \left ( A \times B \right).
\]
Perturbing $\mathbf x$ and $\mathbf y$ slightly assures they are transverse
to the total evaluation map
\[
\mathbf{ev} : \M_{\k,\sigma,l,\varrho}(L,\Sigma,\d',d'') \rightarrow X^l\times
L^{|\k|}.
\]
Each moduli space $\M_{\k,\sigma,l,\varrho}$ contributes a boundary stratum of
the cobordism $\W,$ which we denote by
\[
\partial_G \W_{\sigma,\varrho} : = \M_{\k,\sigma,l,\varrho}(\Sigma,L,\d',d'')\:
_{\mathbf{ev}}\!\!\times_{(\mathbf{x}\times\mathbf{y})\circ\Delta}
\left([0,1]\times\{0,1\}^l \right).
\]
In total, the boundary of the Gromov compactification of $\W$ takes the form
\[
\partial_G \W = \bigcup_{\substack{a \in [1,m],\: \sigma \subset [1,k_a] \\
\varrho \subset [1,l]}}  \partial_G \W_{\sigma,\varrho}.
\]

In general, one might expect an extra term in $\partial_G \W$ coming from
sphere bubbles attached to a constant disk. If there are no marked points on
the disk, this may happen in codimension $1.$ However, assumption
\eqref{eq:aod} precludes this possibility.

Recall that $\Z/2\Z$ acts on $B_{\k,\sigma,l,\varrho}(L,\Sigma,\d',d'')$ by the
involution $\phi_{B^\#}$ that exchanges a disk bubble with its conjugate.
Now, the boundary strata $\W_{\sigma,\varrho}$ are constructed from
$B_{\k,\sigma,l,\varrho}(L,\Sigma,\d',d''), X$ and $L,$ by considering the
vanishing set of a $\Z/2\Z$ equivariant section and then taking various fiber
products with respect to $\Z/2\Z$ equivariant maps. So, each stratum admits a
canonical $\Z/2\Z$ action by an involution that we denote by
$\phi_{\partial\W}.$ We claim
this action is fixed point free and orientation reversing. In other words,
\[
\# \partial_G \W = 0
\]
so that
\begin{align}
0 = \# \partial \W &= \# (\mathbf{ev}^{-1}(\vec x',\vec y') -
\mathbf{ev}^{-1}(\vec x, \vec y) + \partial_G \W)  \notag\\
&= \#\mathbf{ev}^{-1}(\vec x',\vec y') - \# \mathbf{ev}^{-1}(\vec
x, \vec y)\label{eq:pG}
\end{align}
and $N_{\Sigma,\d,\k,l}$ does not depend on the choice of $\vec x,\,\vec y.$

First, we show that $\phi_{\partial\W}$ is fixed point free. Indeed, as noted
above, we are presently considering the case where we may assume there are no
$\phi$-multiply covered disks of positive energy. By definition, a
$\phi$-somewhere injective disk cannot be a fixed point of $\phi_{\partial\W}.$
So, $\phi_{\partial\W}$ could only have a fixed point if a zero energy disk
bubbled off. That would correspond to an interior marked point moving to the
boundary. Clearly, marked points that are constrained away from $L$ cannot
possibly move to the boundary. This is where we use the fact that by Lemma
\ref{lem:AB} we have chosen $(A,f)$ and $(B,g)$ not to intersect $L.$

The following calculations show that $\phi_{\partial\W}$ reverses
orientation. First, we consider the case that $\dim X = 6$
and $L$ is orientable. Since $\dim L = 3,$ by Wu's relations, $w_2(TL) = 0.$ So,
by formula $\eqref{eq:soo},$ the sign of $\phi'_{\L^\#}$ is given by
\begin{equation}\label{eq:ofs}
\s^\# \cong \frac{\mu(d'')}{2} + k'' + l'' + 1.
\end{equation}
Furthermore, $\Z/2\Z$ acts on $\M_{\k,\sigma,l,\varrho}(L,\Sigma,\d',d'')$
with this same sign. Indeed, in the case $l \leq 2,$ we need to add
marked points and fiber product with the divisors along the corresponding
evaluation map. One of the extra marked points could be on the bubble
component, so that $\Z/2\Z$ acts on it non-trivially. However, since the
divisors are chosen to be $\phi$ anti-invariant, the total sign change induced
on $\M_{\k,\sigma,l,\varrho}$ from the extra marked point will be zero. Also,
we note at this point that the sign of the conjugation automorphism of
$PSL_2(\R)$ is zero.

Note that the sign of $\phi_{\partial\W}$ is independent of $l''.$ Indeed,
the action of $\phi$ on $X$ reverses orientation because $\phi^*\omega^3 =
-\omega^3.$ So, the sign of the $\Z/2\Z$ action on the fiber product of
$\M_{\k,\sigma,l,\varrho}(L,\Sigma,\d',d'')$ with $[0,1] \times \{0,1\}^l$ over
$X^l,$ where $\phi$ acts non-trivially on $l''$ of the factors of $X,$
does not depend on $l''.$ On the other hand, a straightforward virtual dimension
calculation shows that $\partial_G\W_{\sigma,\varrho}$ must be empty unless
\begin{equation*}
\mu(d'') = 2k'' + 4l''.
\end{equation*}
This in turn implies that $\mu(d'')/2 \cong k'' \pmod 2,$ or equivalently,
\[
\frac{\mu(d'')}{2} + k'' + 1 \cong 1 \pmod 2.
\]
So, by equation \eqref{eq:ofs}, $\phi_{\partial\W}$ reverses orientation.

Now we turn to the more difficult situation where $\dim X = 4$ and $L$ may not
be orientable. By the Wu relations, $L$ is $Pin^-.$ So, we assume that $\p$ is
given by a standard $Pin^-$ structure. By the same argument as above, we
conclude that $\Z/2\Z$ acts on $\M_{\k,\sigma,l,\varrho}(L,\Sigma,\d',d'')$ with
sign given by formulas \eqref{eq:sn'} or \eqref{eq:sn''} depending on whether
or not $1 \in \sigma.$ Note that $\phi$ preserves the orientation
of $X$ because $\phi^*\omega^2 = \omega^2.$ This implies that
\eqref{eq:sn'},\eqref{eq:sn''}, also give the right signs for the action of
$\phi_{\partial\W}.$ Next, observe that by virtual dimension counting, the
stratum $\partial_G\W_{\sigma,\varrho}$ will be empty unless
\begin{equation}\label{eq:vdr}
\mu(d'') + r = k'' + 2l''
\end{equation}
for $r = 0$ or $-1.$ We claim that this restriction implies that
the signs \eqref{eq:sn'} and \eqref{eq:sn''} always simplify to exactly $1.$

We first consider the case $1 \notin \sigma.$ Using the restriction
\eqref{eq:vdr}, we calculate
\begin{align}
\frac{k''(k'' - 1)}{2} &= \frac{(\mu(d'') + r - 2l'')(\mu(d'') + r -2l'' -
1)}{2} \notag \\
& = \frac{\mu(d'')^2\! +\! 2r\mu(d'')\! +\! r^2\! - \! 4l''(\mu(d'')\!+\!r)\!
+\! 4l''^2 \!- \!\mu(d'')
\!-\! r \!- \! 2l''}{2} \notag \\
& \cong \frac{\mu(d'')(\mu(d'') - 1)}{2} + \frac{r(r-1)}{2} + l'' + r\mu(d'').
\pmod 2 \label{eq:perm1}
\end{align}
Again using the restriction \eqref{eq:vdr}, we calculate
\begin{align}
\Upsilon'(d'',k'') \cong \mu(d'')k'' &\cong \mu(d'')^2 + r\mu(d'')
+2l''\mu(d'')\notag\\
&\cong \mu(d'') + r\mu(d'') \pmod 2. \label{eq:up1}
\end{align}
Substituting equations \eqref{eq:perm1} and \eqref{eq:up1} into \eqref{eq:sn'},
eliminating the remaining $k''$ by \eqref{eq:vdr} and canceling expressions
which occur in pairs yields
\begin{equation*}
\s^{\#'}_\pm(d'',k'',l'') \cong \frac{r(r+1)}{2} + 1 \pmod 2.
\end{equation*}
This is always exactly $1$ since $r = 0$ or $-1.$

We turn now to the case $1 \in \sigma.$ Using the restriction \eqref{eq:vdr},
we calculate
\begin{align}
\frac{(k''-1)(k''-2)}{2} &= \frac{(\mu(d'') + r - 2l'' - 1)(\mu(d'')+r-2l''
-2)}{2} \notag\\
&= \frac{1}{2}\left[\mu(d'')^2 + 2r\mu(d'')+ r^2 - 3\mu(d'')  - 3r
\right.\notag\\
& \qquad\qquad\qquad \left.+ 2 + 4l''^2 + 4l''(\mu(d'') + r) + 6l'' \right]
\notag\\
&\cong \frac{\mu(d'')(\mu(d'') + 1)}{2} + \frac{r(r+1)}{2} \notag\\
&\qquad \qquad \qquad + r\mu(d'') + l''+ 1 \pmod 2. \label{eq:perm2}
\end{align}
Furthermore, using the condition $w_1(d_b) \cong k_b + 1$ and \eqref{eq:vdr},
\begin{equation}
k_b k'' +  k_b  \cong (w_1(d_b) + 1)(\mu(d'') + r + 1) \pmod 2\label{eq:perm2'}.
\end{equation}
Substitute calculations \eqref{eq:perm2} and \eqref{eq:perm2'} in
sign formula \eqref{eq:sn''} and use restriction \eqref{eq:vdr}. Cancelling
pairs of similar terms, we obtain
\begin{align}
\s^{\#''}_- &\cong \frac{r(r+1)}{2} +  r\mu(d'') +  r \notag \\
& \qquad + (w_1(d_b) + 1)(\mu(d'') + r + 1) + \Upsilon'' + w_1(d'_b)w_1(d'')
\notag\\
& \cong r\mu(d'') +  r + (w_1(d_b) + 1)(\mu(d'') + r + 1) + \Upsilon''+
w_1(d'_b)w_1(d''). \label{eq:sn'a}
\end{align}
Here, the second congruence follows from the fact $r = 0$ or $-1.$

We now expand $\Upsilon''$ to further analyze $\s^{\#''}.$ Using the
fact that
\[
w_1(d_b) = w_1(d'_b) + w_1(\partial d''),
\]
it is easy to verify that
\begin{equation}\label{eq:up2}
\Upsilon'' = w_1(d_b)(k'' - 1) + w_1(\partial d'') k'.
\end{equation}
By repeatedly applying restriction \eqref{eq:vdr}, the condition that
$w_1(d_a) = 1 + k_a$ and the fact that $\mu(d'') \cong w_1(\partial d'') \pmod
2,$ we calculate,
\begin{align*}
w_1(\partial d'') k' & \cong w_1(\partial d'')(k_b - k'') \\
& \cong  w_1(\partial d'')(w_1(d_b) + 1 + \mu(d'') + r) \\
& \cong w_1(\partial d'')(w_1(d_b') + w_1(\partial d'') + 1 + \mu(d'') + r) \\
& \cong w_1(\partial d'')w_1(d_b') + \mu(d'')(1 + r) \pmod 2.
\end{align*}
Substituting this calculation in formula \eqref{eq:up2}, and using restriction
\eqref{eq:vdr} again, we obtain
\[
\Upsilon'' \cong w_1(d_b)(\mu(d'') + r  + 1) +  w_1(\partial d'')w_1(d_b') +
\mu(d'')(1 + r) \pmod 2.
\]
Substituting this expression for $\Upsilon''$ in formula \eqref{eq:sn'a} an
cancelling all repeated terms leaves
\begin{align*}
\s^{\#''}_- &\cong r \mu(d'') + r + (1)(\mu(d'') + r + 1) + \mu(d'')(1 +
r) \\
&\cong r \mu(d'') + r + \mu(d'') + r + 1 + \mu(d'') + \mu(d'') r \\
& \cong 1 \pmod 2,
\end{align*}
as desired.

\section{An equivariant Kuranishi structure}\label{sec:km}
In this section, we complete the proof of invariance of $N_{\Sigma,\d,\k,l}$ in
the case when $(X,L)$ may admit holomorphic disks of Maslov index $0.$
If $\dim L =2,$ the expected dimension of holomorphic disks with Maslov index
$0$ is negative. For generic $J,$ by Lemma \ref{lem:mc}, such disks don't
exist. So, we consider the case $\dim L = 3.$ By assumption, $L$ is
orientable. The main tool we use to prove invariance in this case is the notion
of a Kuranishi structure, introduced in \cite{FO} and extended in \cite{FOOO}.
For a summary of relevant information on Kuranishi structures, see Appendix
\ref{sec:aks}.

Suppose $(X,\K)$ is a space with Kuranishi structure
\[
\K = (V_p,E_p,\Gamma_p,s_p,\psi_p).
\]
Let $\iota$ be an involution of $X.$
\begin{dfn}
An \emph{extension} $\tilde \iota$ of an involution $\iota$ to an involution
of $\K$ consists of $\Gamma_p$-equivariant maps $\iota_p : V_p \rightarrow
V_{\iota(p)}$ and $\hat \iota_p : E_p \rightarrow E_{\iota(p)}$ covering
$\iota_p$ such that
\begin{enumerate}
\renewcommand{\theenumi}{(E\arabic{enumi})}
\renewcommand{\labelenumi}{\theenumi}
\item
$\iota_{\iota(p)}\circ \iota_p = \id_{V_p}.$
\item
$s_{\iota(p)} \circ \iota_p = \hat \iota_p \circ s_p.$
\item
$\psi_{\iota(p)} \circ \iota_p|_{s_p^{-1}(0)} = \iota \circ \psi_p.$
\item
$\iota_q$ maps $V_{pq} \subset V_q$ to $V_{\iota(p)\iota(q)} \subset
V_{\iota(q)}.$
\item
$\iota_p\circ\varphi_{pq} = \varphi_{\iota(p)\iota(q)}\circ\iota_q$ and
$\hat\iota_p\circ\hat\varphi_{pq} =
\hat\varphi_{\iota(p)\iota(q)}\circ\hat\iota_q.$
\end{enumerate}
\end{dfn}
Note that $\iota_p,\,\hat\iota_p,$ induce bundle morphisms
\[
 \iota^T_p:\det(TV_p) \otimes \det(E_p) \rightarrow \det(TV_{\iota(p)}) \otimes
\det(E_{\iota (p)})
\]
covering $\iota_p.$ Now, suppose that $(X,\K)$ has a tangent bundle given by
$\Phi_{pq}.$ We say that $\tilde \iota$ acts smoothly on $(X,\K)$ if
\begin{equation*}
\Phi_{\iota(p)\iota(q)}\circ \hat\iota_q =  \hat\iota_p \circ \Phi_{pq}.
\end{equation*}
If $\tilde \iota$ acts smoothly and $(X,\K)$ is oriented, the bundle morphisms
$\iota^T_p$ may either preserve or reverse the orientation of $\K$ over each
connected component of $X.$ Furthermore, let $\G = (P,V_p',E_p',s_p',\psi_p')$
be a good coordinate system for $\K.$ We say that $\G$ is compatible with
$\tilde \iota$ if $\iota(P) = P$ and $\iota_p(V_p') = V_{\iota(p)}'.$ Finally,
a collection $\{t_p\}_{p\in P}$ of multisections of $E_p'$ is said to be
$\iota_p - \tilde\iota_p$ equivariant if $t_{\iota(p)}\circ \iota_p =
\tilde \iota_p \circ t_p$ for all $p \in P.$

Now, we give the main idea of the proof. As we will explain, it is
possible to construct an oriented Kuranishi structure with tangent
bundle $\K = (V_p,E_p,\Gamma_p,s_p,\psi_p)$ on $\W.$ Moreover, we
may construct $\K$ so that the induced Kuranishi structure on
$\partial_G\W$ admits a smooth orientation reversing involution
$\tilde \phi_{\partial\W}$ extending $\phi_{\partial\W}.$ There
exists a good coordinate system $\G = (P,V_p',E_p',s_p',\psi_p')$ for
$\K$ such that its restriction to $\partial_G\W$ is compatible with
$\tilde\phi_{\partial\W}.$ Finally, we may choose
$\phi_{\partial\W^p} - \hat\phi_{\partial\W^p}$ equivariant
generic transverse multisections $s_{p,n}'$ satisfying conditions
\ref{c:1}-\ref{c:4} of Theorem \ref{thm:kp} and coinciding exactly with $s_{p}'$
away from $\psi_p^{-1'}(\partial_G\W).$

Note that the charts of the induced Kuranishi structure on
$\partial\W$ are just $\partial V_p.$ We define
\[
\partial_G V_p : = \psi_p^{-1}(\partial_G\W) \subset \partial V_p.
\]
The same applies for the charts of the good coordinate system and we
use the analogous notation. The vanishing sets of the $s_{p,n}'$
define a $1$-dimensional simplicial complex with boundary
contained in the $\partial V_p'.$ The boundary is simply a
collection of points with signed rational weights. By the same
reasoning as in equation \eqref{eq:pG}, it suffices to show that
the part of this boundary contained in $\partial_G V_p'$ has total
cardinality zero. Again, we use the sign reversing involutions
$\phi_{\partial \W^p}$ to cancel the points in pairs. The only
slight complication arises because $\phi_{\partial \W^p}$ may have
fixed points. However, this is easily resolved by the observation
that a point $x$ of the vanishing set of $s_{p,n}'$ fixed by
$\phi_{\partial \W^p}$ must have weight zero. Indeed, by
definition of transversality for multisections, each branch
$s_{p,n}^{i'}$ of $s_{p,n}'$ is transverse to zero. So, if
$s_{p,n}^{i'}$ vanishes at $x,$ the differential
\[
ds_{p,n}^{i'} : T_x\partial V_p' \stackrel{\sim}{\rightarrow}
(E_p')_x
\]
defines a non-zero element $\omega \in \det(TV_p') \otimes
\det(E_p').$ Since $s_{p,n}'$ is $\phi_{\partial
\W^p}-\hat\phi_{\partial \W^p}$ equivariant,
\[
\hat\phi_{\partial \W^p}\circ s_{p,n}^{i'} \circ \phi_{\partial
\W^p}
\]
must also be a branch of $s_{p,n}',$ defining an element $\omega'
\in \det(TV_p') \otimes \det(E_p').$ But since $\tilde
\phi_{\partial \W}$ is orientation reversing, we know that
$\omega$ and $\omega'$ belong to opposite components of
$\det(TV_p') \otimes \det(E_p')\setminus \{0\}.$ So, the branches
of $s_{p,n}'$ that vanish at $x$ come in pairs that induce
opposite orientations on $x.$ Therefore, the total weight of $x$
is zero.

We turn our attention to the construction of $\K.$ Recall from \eqref{eq:W} that
\[
\W = \M_{\k,l}(\Sigma,L,\d)\:
_{\mathbf{ev}}\!\!\times_{(\mathbf{x}\times\mathbf{y})\circ\Delta}
\left([0,1]\times\{0,1\}^l \right).
\]
Let $\partial\M_{\k,l}(\Sigma,L,\d)$ denote the union of all
the strata of the Gromov compactification of $\M_{\k,l}$ except $\M_{\k,l}$
itself. Let $\phi_{\partial\M}$ denote the involution of $\partial\M_{\k,l}$
induced by $\phi_{B^\#}.$ As explained in \cite[Appendix 2]{FOOO}, a weakly
submersive Kuranishi structure $\K_\M$ on $\overline\M_{\k,l}(\Sigma,L,\d)$
naturally induces a Kuranishi structure $\K$ on the fiber product $\W.$  If we
let $\Z/2\Z$ act on $X$ by $\phi$ and on $[0,1]\times\{0,1\}^l$ by exchanging
$0$ and $1,$ $\mathbf{x}\times\mathbf{y}$ is clearly $\Z/2\Z$ equivariant. So,
if $\K_\M|_{\partial\M_{\k,l}}$ admits an extension of $\phi_{\partial\M},$
then $\K|_{\partial_G\W}$ will admit an extension of $\phi_{\partial\W}.$
Consequently, we focus on the construction of $\K_\M.$

We assume $\Sigma = D^2$ and $l \geq 2.$ The other cases are similar.
Interpreting \eqref{eq:M} as a fiber product, we have
\[
 \M_{\k,l}(\Sigma,L,\d) = \widetilde \M_{\k,l}(\Sigma,L,\d) \times_{\Sigma}
\left(s_0 \times \ell\right).
\]
Let $\partial \widetilde\M_{\k,l}(\Sigma,L,\d)$ denote the union
of all strata of the Gromov compactification of
$\widetilde\M_{\k,l}$ except $\widetilde\M_{\k,l}$ itself. Let
$\phi_{\partial\widetilde\M}$ denote the involution of
$\partial\widetilde\M_{\k,l}$ induced by $\phi_{B^\#}.$  Again, we reduce to
constructing a weakly submersive Kuranishi structure
$\K_{\widetilde \M}$ on the Gromov compactification of
$\widetilde\M_{\k,l}(\Sigma,L,\d)$ such that its restriction to
$\partial\widetilde \M_{\k,l}(\Sigma,L,\d)$ admits an extension of
$\phi_{\partial\widetilde\M}.$

In \cite{FOOO}, Fukaya et al.\ constructed a Kuranishi structure on the moduli
space of $J$-holomorphic disks with Lagrangian boundary conditions. Their
construction generalizes immediately to the higher genus fixed
conformal structure situation considered in this paper. Away from
$\partial\widetilde\M_{\k,l},$ we use their construction without further
discussion. Near $\partial\widetilde\M_{\k,l},$ we must impose additional
conditions on several choices made in the construction to make sure we can find
an extension of $\phi_{\partial\widetilde\M}.$
So, we briefly recount the idea of the construction of the Kuranishi
neighborhood of a point $p \in \partial\widetilde\M_{\k,l}$ given in
\cite{FOOO}. We assume that $p$ is a stable map of two components. The
construction for more components is similar. By definition, $p$ is an
equivalence class of quadruples $\u = (\hat\Sigma,u,\vec z,\vec w) \in
B_{\k,\sigma,l,\varrho}^{1,p}(L,\Sigma,\d',d'')$ such that
$\bar\partial_{J,\nu}^\# u = 0.$ The equivalence relation in the case of two
components equates reparameterizations of the bubble component. When there are
more than two components, the equivalence relation also takes into consideration
automorphisms of the underlying tree of the stable map. We choose some $\u$
such that $[\u] = p.$ Locally trivializing $\E^\#$ by parallel translation
and projection to $\Lambda^{0,1}(TX),$ we define
\[
 D_\u^\# = D_\u\bar\partial_{J,\nu}^\#:
T_\u B_{\k,\sigma,l,\varrho}^{1,p}(L,\Sigma,\d',d'') \longrightarrow \E^\#_\u
\]
to be the linearization of $\bar\partial_{J,\nu}^\#$ at $\u.$ Since the bubble
component $u''$ of $u$ may be $\phi$-multiply covered, $D_\u^\#$ may not be
surjective even for generic $J$ and $\nu.$ However, since $D_\u^\#$ is Fredholm,
we may choose a finite dimensional subspace $E_\u \subset \E^\#_\u$ such that if
$\pi : \E^\#_\u \rightarrow \E^\#_\u/E_\u$ denotes the natural projection, then
$\pi\circ D_\u^\#$ is surjective. Possibly enlarging $E_\u,$ we may assume that
the evaluation maps from $\ker \pi \circ D_\u^\#$ to $T_{u(z_i)}L$ and
$T_{u(w_j)}X$ are surjective. This is necessary for the Kuranishi structure we
are describing to be weakly submersive. By elliptic regularity, we may choose
$E_\u$ to consist of smooth sections of $\Lambda^{0,1}(u^*TX).$ By the unique
continuation theorem, we may assume these sections are compactly supported away
from the singular point $z_0.$

Let $\u_\epsilon = (\hat\Sigma_\epsilon,u_\epsilon,\vec z_\epsilon,\vec
w_\epsilon)$ be sufficiently close to $\u.$ More precisely, choose a
small $\delta >0.$ We allow $\hat \Sigma_\epsilon$ to differ from $\hat \Sigma$
in a $\delta$ neighborhood $N_\delta$ of the singular point $z_0.$ Either $z_0$
may move slightly, or small neighborhoods of $z_0$ in each component of $\hat
\Sigma$ may be removed and their boundaries glued to smooth the singularity. In
particular, outside $N_\delta,$ there exists a canonical
identification of $\hat\Sigma$ with $\hat\Sigma_\epsilon.$ The pre-gluing
construction explained in detail in \cite{FO,MS} gives a smooth map $\tilde u :
(\hat\Sigma',\partial\hat\Sigma') \rightarrow (X,L)$ that agrees with $u$
outside $N_\delta$ and stays $\delta$-close to $u(z_0)$ within
$N_\delta.$ We assume that $u_\epsilon$ is $\delta$ close to $\tilde u$ in the
$W^{1,p}$ norm. Also, we assume that $\vec z_\epsilon,\vec w_\epsilon,$
are $\delta$-close to $\vec z, \vec w.$  Then, for $\delta$ sufficiently
small, there exist unique shortest geodesics from $u(z)$ to $u_\epsilon(z)$ for
each $z \in \hat\Sigma\setminus N_\delta.$ For $\delta$ sufficiently
small, we may assume $N_\delta$ is disjoint from the support of the sections of
$\Lambda^{0,1}(u^*TX)$ constituting $E_\u.$ So, we may parallel translate $E_\u$
along length minimizing geodesics and then project to $\Lambda^{0,1}(TX)$ to
obtain a subspace of $\E^{(\#)}_{\u_\epsilon}.$ Here, we parenthesize $\#$
because $\u_\epsilon$ may be an irreducible $W^{1,p}$ stable map. If $\delta$ is
sufficiently small this subspace has constant dimension. We let $\underline
E_\u$ denote the sub-bundle of $\E^{(\#)}$ so obtained. Similarly, we let
$\underline \pi : \E^{(\#)} \rightarrow \E^{(\#)}/\underline E_\u$ define the
projection to the quotient bundle.

According to \cite{FOOO,FO}, we may essentially define $V_p$ to be the set of
$\u_\epsilon$ as above such that $\underline \pi\circ
\bar\partial^{(\#)}_{J,\nu}u_\epsilon = 0.$ Then, we define $E_p =
\underline E_\u|_{V_p}$ and $s_p = \bar\partial^{(\#)}_{J,\nu}.$ The definition
of $\psi_p$ is tautological. A great deal of hard analysis then shows that $V_p$
is actually a smooth manifold with boundary modelled on
\[
 \ker\left(\pi\circ D_\u^\#\right) \times (R,\infty].
\]
However, we do not need to know the details of this analysis at all for our
purposes. Also, we note that it is necessary to further enlarge
$E_\u$ in order to construct the maps $\varphi_{pq},\hat\varphi_{pq}.$ However,
this step is essentially formal, and it is not hard to make it
$\phi$-equivariant. So, we do not discuss it here and instead refer the reader
to \cite[end of Section 18]{FOOO}, or for more detail \cite[Section 15]{FO}.
Finally, we note that canonical orientation of $\det(D^\#)$ induces an
orientation of $\det(TV_p) \otimes \det(E_p).$ Indeed, note that there exist
natural isomorphisms $\ker Ds_p \simeq \ker D\bar\partial_{J,\nu}$ and
$\coker Ds_p \simeq \coker D\bar\partial_{J,\nu}.$ On the other hand, the exact
sequence
\[
 0 \rightarrow \ker Ds_p \rightarrow TV_p \rightarrow E_p \rightarrow \coker
Ds_p \rightarrow 0
\]
induces a natural isomorphism
\[
 \det(TV_p) \otimes \det(E_p) \simeq \det(\ker Ds_p)\otimes\det(\coker Ds_p).
\]

We now detail additional conditions on the choices made in the above
construction, which ensure that $\K_{\widetilde\M}|_{\partial\widetilde\M}$
admits an extension of $\phi_{\partial\widetilde\M}.$ Clearly, we must
choose a $\phi$-invariant metric on $X$ for measuring all distances,
and constructing geodesics and parallel transport. We consider two cases: First
suppose $p_\phi : = \phi_{\partial\widetilde\M}(p) \neq p.$ We simultaneously
construct Kuranishi neighborhoods of $p$ and $p_\phi$ as
well as the extension of $\phi_{\partial\widetilde\M}.$ Indeed, given a
representative $\u$ of $p,$ we choose the
representative $\u_\phi : = \phi_{B^\#}(\u)$ of $p_\phi.$ Furthermore, we
choose $E_{\u_\phi} = \phi_{\E^{\#}}(E_\u).$ This is compatible with the
construction above because of the $\phi$-invariance of the metric. Again, by
$\phi$-invariance of the metric, it follows that $\phi_{\E^\#}$ maps
$\underline E_\u|_{B^\#}$ to $\underline E_{\u_\phi}|_{B^\#}.$
Since $\bar\partial_{J,\nu}^\#$ is $\phi_{B^\#}-\phi_{E^\#}$ equivariant, it
follows that $\phi_{B^\#}$ maps $\partial V_p$ to $\partial V_{p_\phi}.$ So, we
define the extension of $\phi_{\partial\widetilde\M}$ by
\[
\phi_{\partial\widetilde\M^p} : = \phi_{B^\#}|_{\partial V_p}, \qquad
\hat \phi_{\partial\widetilde\M^p} : = \phi_{E^\#}|_{E_p}.
\]
On the other hand, suppose that $\phi_{\partial\widetilde\M}(p) = p.$ This may
happen when $u''$ is a $\phi$-multiply covered disk of even multiplicity. It is
not hard to see that we may choose $\u$ representing $p$ such that
$\phi_{B^\#}(\u) = \u.$ Indeed, this follows from the fact that all
anti-holomorphic involutions of $D^2$ are conjugate under the action of
$PSL_2(\R).$ Furthermore, we choose $\E_{\u}$ to be
$\phi_{\E^\#}$-invariant. This said, we may define the extension of
$\phi_{\partial\widetilde\M}$ exactly as above. This completes the construction
of $\K_{\widetilde M}.$

Now, a minor adaptation of the proof of \cite[Lemma 6.3]{FO} gives
the good cover $\G.$ To obtain generic transverse multisections $s_{p,n}'$ such
that $s_{p,n}'|_{\partial_GV_p'}$ is $\tilde \phi_{\partial W}$ equivariant, we
use an argument from \cite[Section 11]{FOOO}. That is, since
$\K|_{\partial_G\W}$ admits an extension of $\phi_{\partial\W},$ it descends
naturally to a Kuranishi structure on $\W/(\Z/2\Z).$ Similarly,
$\G|_{\partial_G\W}$ descends to a good coordinate system $\hat\G=
(\hat P,\hat V_p',
\hat E_p',\hat s_p',\hat \psi_p')$ on $\partial_G\W/(\Z/2\Z).$ We denote by
\[
 \Pi : \partial\W \rightarrow  \partial\W/(\Z/2\Z)
\]
the quotient map. We can apply the standard machinery of Kuranishi
structures developed in \cite[Chapter 1]{FO}, reviewed in Theorem \ref{thm:kp},
to obtain generic transversal multisections $\hat s_{p,n}'$ perturbing $\hat
s_p.$ Pulling back $\hat s_{p,n}'$ under $\Pi,$ we obtain
$\phi_{\partial\W^p}-\hat\phi_{\partial\W^p}$ equivariant transverse
multisections over $\partial_G V'_p.$ Since $V_p'$ is a
manifold with corners, it is not hard to extend a transverse section from
$\partial_G V'_p$ to a neighborhood of $\partial_G V'_p.$  In fact, away from
corners, a neighborhood of $\partial_G V'_p$ is diffeomorphic to $[0,1)
\times \partial_G V_p'$ and we can extend sections transverse to zero as
constants over $[0,1).$ Since $\dim K = 1,$ a transverse section cannot have
zeros at corners of $V_p',$ so we can extend near corners
arbitrarily.
Note that by choosing a generic inhomogeneous perturbation $\nu$ we have
already made $s_p'$ transverse away from $\partial_G V_p'.$ So, we may patch the
$s_p'$ with the extensions we have just constructed to obtain the desired
$s_{p,n}',$ maintaining transversality everywhere. This completes the proof of
Theorem \ref{thm:inv}.
\end{proof}

\section{Calculations}\label{sec:calc}
In this section we prove Theorem \ref{thm:eq} and Example \ref{ex:q}. The main
tool of the proofs is the notion of a short exact sequence of Cauchy-Riemann
$Pin$ boundary value problems. The first step in understanding short exact
sequences of $Pin$ boundary problems is to understand short exact sequences of
bundles with $Pin$ structure. Suppose
\[
0\longrightarrow V' \longrightarrow V \longrightarrow V'' \longrightarrow 0
\]
is a short exact sequence of real vector bundles over a base $B.$
\begin{lemma}\label{lem:indp}
Assume that at least one of $V'$ and $V''$ is orientable. Then a $Pin$
structure on any two of $V,V',V,''$ naturally induces a $Pin$ structure on the
third.
\end{lemma}
\begin{proof}
For the proof of this Lemma, we write $\dim V' = n$ and $\dim V'' = m.$
Choosing a metric on $V,$ we may identify $V \simeq V' \oplus V''.$
By symmetry of the direct sum, we may assume that $V'$ is orientable. We use the
orientation of $V'$ to reduce its structure group to $SO(n).$ So, the Lemma
follows from the existence of the commutative square of group homomorphisms
\[
\xymatrix{ Spin(n) \times Pin(m) \ar[r]\ar[d]& Pin(n+m)\ar[d] \\
SO(n) \times O(n) \ar[r] & O(n+m).
}
\]
Indeed, we work on the level of transition functions which satisfy the co-cycle
condition. The commutativity of the two factors of the product of groups ensures
that the direct sum does not effect the cocyle condition.
\end{proof}
\begin{lemma}\label{lem:p1}
Let $V \rightarrow B$ be a real vector bundle. If $B$ is a one dimensional
manifold, then $V$ automatically carries a $Pin$ structure. On the other
hand, if $V$ is one dimensional and admits a $Pin$ structure, the $Pin$
structure on $V$ may be chosen canonically.
\end{lemma}
\begin{proof}
The first claim of the Lemma follows because the obstruction to the existence
of a $Pin$ structure is a second cohomology class. The second claim follows when
$B = \R P^1$ because, as noted in Remark \ref{rem:n=1}, all automorphisms of a
line bundle over $\R P^1$ preserve $Pin$ structure and so we can induce a $Pin$
structure canonically from a previously chosen one on $\tau_\R \rightarrow \R
P^1$ or $\underline\R \rightarrow \R P^1.$  This extends to general $B$ because
a $Pin$ structure on $V \rightarrow B,$ if it exists, is determined by its
restriction to loops in $B.$
\end{proof}

Now, let $E,E',E'',$ be complex vector bundles over a Riemann surface with
boundary $\Sigma$ and let $F,F',F'',$ be totally real subbundles of $E,E',E'',$
respectively, over $\partial \Sigma.$ Suppose further that we have an exact
sequence
\[
0 \longrightarrow E' \longrightarrow E \longrightarrow E'' \longrightarrow 0
\]
that restricts to an exact sequence
\begin{equation}\label{eq:Fex}
0 \longrightarrow F' \longrightarrow F \longrightarrow F'' \longrightarrow 0.
\end{equation}
We refer to such a short exact sequence as a short exact sequence of pairs of
vector bundles. Let $\p,\p',\p'',$ be $Pin$-structures on $F,F',F'',$
respectively. We say that $\p$ is \emph{compatible} with the short exact
sequence \eqref{eq:Fex} if $\p$ agrees with the $Pin$ structure naturally
induced on $F$ by $\p'$ and $\p''$ by Lemma \ref{lem:indp}.
If $F'$ or $F''$ is one dimensional, even if it does not come equipped with a
$Pin$ structure, we extend the notion of compatibility by equipping it with the
canonical $Pin$ structure of Lemma \ref{lem:p1}.

\begin{dfn}\label{dfn:sescr}
A \emph{short exact sequence} of Cauchy-Riemann $Pin$ boundary value
problems
\[
0 \longrightarrow \underline D' \longrightarrow \underline D \longrightarrow
\underline D'' \longrightarrow 0
\]
consists of
\begin{itemize}
\item
An exact sequence of pairs of vector bundles
\begin{equation}\label{eq:EFes}
 0 \longrightarrow (E',F') \longrightarrow (E,F) \longrightarrow (E'', F'')
\longrightarrow 0
\end{equation}
such that at least one of $F'$ and $F''$ is orientable.
\item
Orientations on each of $F,F',F'',$ that is orientable. If all three are
orientable, we assume the orientation of $F$ agrees with the orientation
induced from $F'$ and $F''.$
\item
$Pin$ structures $\p,\p',\p'',$ on $F,F',F'',$ respectively. If $F'$
(resp. $F''$) has dimension $1,$  we do not require $\p'$ (resp. $\p''$) as
part of the definition, since it may be chosen canonically by Lemma
\ref{lem:p1}.
\item
Cauchy-Riemann operators
\begin{gather*}
D : \Gamma(E,F) \longrightarrow \Gamma\left(\Omega^{0,1}(E)\right),\quad D':
\Gamma(E',F') \longrightarrow \Gamma\left(\Omega^{0,1}(E')\right),\\
D'': \Gamma(E'',F'') \longrightarrow \Gamma\left(\Omega^{0,1}(E'')\right),
\end{gather*}
such that the diagram
\begin{equation*}
\xymatrix{
0 \ar[r] & \Gamma\left(\Omega^{0,1}(E')\right) \ar[r] &
\Gamma\left(\Omega^{0,1}(E)\right) \ar[r] &
\Gamma\left(\Omega^{0,1}(E'')\right)\ar[r] & 0 \\
0 \ar[r] & \Gamma(E',F') \ar[r]\ar[u]^{D'} & \Gamma(E,F)\ar[r]\ar[u]^D &
\Gamma(E'',F'') \ar[r]\ar[u]^{D''} & 0
}
\end{equation*}
commutes.
\end{itemize}
\end{dfn}
Note that a short exact sequence of Cauchy-Riemann $Pin$ boundary value
problems is an example of a short exact sequence of Fredholm operators. See
Definition \ref{dfn:sesf}.
\begin{prop}\label{prop:sescr}
Let
\[
0 \longrightarrow D' \longrightarrow D''\longrightarrow D'\longrightarrow 0
\]
be a short exact sequence of Cauchy-Riemann $Pin$ boundary value problems.
Up to a universal sign, the isomorphism
\[
\det(D')\otimes\det(D'') \stackrel{\sim}{\longrightarrow} \det(D)
\]
given by Lemma \ref{lem:ses} respects the canonical orientations of Proposition
\ref{prop:orl} if and only if $\p$ is compatible with the short exact sequence.
The universal sign depends on the dimension of $E,E',E'',$ the topology of
$\Sigma$ and the orientability of $F,F',F'',$ restricted to each boundary
component of $\Sigma.$
\end{prop}
\begin{proof}
By a deformation argument, similar to the proof of Proposition \ref{prop:orl},
we would like to reduce the problem to a standard short exact sequence.
Throughout the proof, we assume that $\p$ is compatible with the
short exact sequence. The other case follows from Lemma \ref{lem:cp}. As in the
proof of Proposition \ref{prop:orl}, degenerate $\Sigma$ along curves $\gamma_a$
to a nodal Riemann surface $\hat \Sigma$ with nodal points $\hat \gamma_a,$ and
irreducible components $\Delta_a \simeq D^2$ and $\tilde \Sigma \simeq
\Sigma/\partial \Sigma.$ Simultaneously, degenerate $E,E',E'',$ to complex
vector bundles $\hat E, \hat E', \hat E'',$ over $\hat \Sigma$ that all satisfy
condition \eqref{eq:bis} for appropriate $n.$ Now, by degenerate, we mean
identify the fibers of $E|_{\gamma_a}$ (resp. $E'|_{\gamma_a},E''|_{\gamma_a}$)
with the single fiber $\hat E_{\hat\gamma_a}$ (resp.
$E'_{\hat\gamma_a},E''_{\hat\gamma_a}$). Such a degeneration satisfying
condition \eqref{eq:bis} is unique up to homotopy. Furthermore, we may choose
the degeneration of $E$ to extend the degeneration of $E'.$ These two
degenerations induce a degeneration of $E''$ via the short exact sequence
\eqref{eq:EFes}. So, we may assume that there exists a natural induced short
exact sequence
\begin{equation}\label{eq:hEFes}
0 \longrightarrow (\hat E',\hat F') \longrightarrow (\hat E,\hat F)
\longrightarrow (\hat E'',\hat F'') \longrightarrow 0.
\end{equation}

Choose a particular isomorphism \eqref{eq:bis} for $(\hat E'|_{\Delta_a},\hat
F'|_{\partial \Delta_a}).$ Extend it to an isomorphism \eqref{eq:bis} for $(\hat
E|_{\Delta_a},\hat F|_{\partial \Delta_a}).$ Denote the canonical bundle pairs
over $(D^2,\partial D^2)$ by
\[
(E_{i,n}, F_{i,n}) :=
\begin{cases}
(\tau\oplus \underline\C^{n-1}, \tau_\R \oplus \underline\R^{n-1}), & i = -1  \\
(\underline\C^n,\underline\R^n), & i = 0 .
\end{cases}
\]
The following diagram shows that we have a naturally induced isomorphism
\eqref{eq:bis} for $(\hat E''|_{\Delta_a},\hat
F''|_{\partial \Delta_a}):$
\begin{equation*}
\xymatrix{
0 \ar[r] & (E_{i,n},F_{i,n}) \ar[r] & (E_{i+j,n+m},F_{i+j,n+m}) \ar[r] &
(E_{j,m},F_{j,m}) \ar[r] & 0 \\
0 \ar[r] & (\hat E'|_{\Delta_a}, \hat F'|_{\partial\Delta_a}) \ar[r]\ar[u]
& (\hat E|_{\Delta_a}, \hat F|_{\partial\Delta_a}) \ar[r]\ar[u] & (\hat
E''|_{\Delta_a}, \hat F''|_{\partial\Delta_a}) \ar[r]\ar@{.>}[u] & 0.
}
\end{equation*}
Here, the top row makes sense because by assumption, either $i$ or $j$ or both
are $0.$ So, the top row is tautological. The bottom row is just a
restriction of short exact sequence \eqref{eq:hEFes}. Since the
morphisms in the top row commute with the canonical Cauchy-Riemann operators on
the bundles $E_{i,n},$ the isomorphisms \eqref{eq:bis} just chosen induce
Cauchy-Riemann operators $D_a,\,D_a'$ and $D_a'',$ on $\hat E|_{\Delta_a},\hat
E'|_{\Delta_a}$ and $\hat E''|_{\Delta_a},$ respectively, that commute with the
morphisms of the short exact sequence \eqref{eq:hEFes}. We claim that if $n +
m \geq 3,$ the preceding construction is unique up to homotopy. Indeed,
choosing a metric on $\hat E|_{\Delta_a}$ induces a splitting of the bottom
row of the above diagram. Since the space of metrics is contractible, this
choice is unique up to homotopy. Then, the middle vertical morphism determines
both of the other vertical morphisms. On the other hand, when $n+m \geq 3,$ the
middle vertical morphism is unique up to homotopy by Lemma \ref{lem:aut}. For
$n + m = 2,$ we use a stabilization argument as in the proof of Proposition
\ref{prop:orl} to reduce to the case $n + m = 3.$

Choose operators $\tilde D,\,\tilde D'$ and $\tilde D'',$ on $\hat E|_{\tilde
\Sigma},\, \hat E'|_{\tilde \Sigma}$ and $\hat E''|_{\tilde \Sigma},$ compatible
with the
short exact sequence \eqref{eq:hEFes}. Note that the induced isomorphism
\[
\det(\tilde D') \otimes \det(\tilde D'') \stackrel{\sim}{\longrightarrow}
\det(\tilde D),
\]
always preserves the canonical complex orientations of each side.
Finally, choose homotopies of operators $\mathbf D_t, \mathbf D'_t$ and
$\mathbf D''_t,$ on $E,\,E'$ and $E'',$ respectively, compatible with the
short exact sequence \eqref{eq:hEFes}, such that
\[
\mathbf D_0 = D, \qquad \mathbf D_1 = \#_a D_a \# \tilde D
\]
and similarly for $\mathbf D'_t$ and $\mathbf D''_t.$ Applying Lemma
\ref{lem:ses} to the short exact sequence of families of Fredholm operators,
\[
0 \longrightarrow \mathbf D'_t \longrightarrow \mathbf D_t\longrightarrow
\mathbf D''_t \longrightarrow 0
\]
proves that the sign is universal, as claimed.
\end{proof}

We now prove a technical lemma that will be useful in the proof of
Theorem \ref{thm:eq}. The idea of the proof is taken from
\cite[proof of Theorem C.1.10(iii)]{MS}, which in turn follows the work of
Hofer-Lizan-Sikorav \cite{HLS}.
\begin{lemma}\label{lem:surj}
Let $(E,F) \rightarrow (D^2,\partial D^2)$ be a vector bundle pair
with $\dim_\C E = 1,$ and denote its Maslov index by $\mu =
\mu(E,F) \geq -1.$ Let $D$ be a real Cauchy-Riemann operator on
$E.$ Let $z_1,\ldots,z_k \in \partial D^2,$ and $w_1,\ldots,w_l \in D^2$ be
distinct marked points. Assume $l + 2k = \mu + 1.$ Denote by
\[
ev: \ker(D) \rightarrow \R^k \oplus \C^l
\]
the evaluation map defined by
\[
\xi \rightsquigarrow (\xi(z_1),\ldots,\xi(z_k),\xi(w_1),\ldots,\xi(w_l)),
\quad \xi \in \ker D.
\]
Then $ev$ is always surjective.
\end{lemma}
\begin{proof}
The Fredholm index of $D$ is well known to be $\mu + 1.$ It
follows that the Fredholm index of $D \oplus ev$ is $0.$ So, if
$ev$ is not surjective, there must exist some non-zero $\xi \in \ker(D)$
such that $\xi(z_i) = 0,\, i = 1,\ldots,k,$ and $\xi(w_j) = 0,\, j = 1,
\ldots, l.$ According to \cite[proof of Theorem C.1.10(iii)]{MS}, there exists a
complex linear Cauchy-Riemann operator $D'$ on $E$ and a function $u \in
W^{1,p}(D^2,\C),$ such that $D'(u\xi) =0.$ This leads to a
contradiction since a holomorphic section of $(E,F)$ may have at
most $\mu$ zeros where interior zeros are counted twice.
\end{proof}
\begin{rem}
Note that it is crucial for this argument that the underlying Riemann surface is
a disk. Otherwise, the Fredholm index of $D$ is not $\mu + 1.$ This explains, at
least in part, why Welschinger's counting scheme does not immediately extend to
curves of higher genus.
\end{rem}

\begin{proof}[Proof of Theorem \ref{thm:eq}]
First, we treat the case $\dim L = 2.$ Welschinger's invariants are defined
only in the strongly semipositive case when $\Sigma = D^2.$ So, we take $\nu =
0.$ As explained in Section \ref{sec:def2}, we could consider the moduli space
defined by quotienting by the action of $PSL_2(\R)$, but for this proof it seems
more natural to take a section of the $PSL_2(\R)$ action. In particular, we add
two marked points constrained to divisors, even when $l > 0,$ as explained in
Remark \ref{rem:l>0}. As in Section \ref{sec:def2}, we denote by $(A,f)$ and
$(B,g)$ the divisor constraint pseudo-cycles. In this proof, we will refer to
the extra added marked points as $z_{-1}$ and $z_{-2}.$ As in Section
\ref{sec:def2}, we fix $z_{-1}$ to be at a point $s_0$ and we fix $z_{-2}$ to
lie on a line $\ell.$

Recall from Section \ref{sec:def2} that, by definition,
\[
N_{D^2,d,k,l} : = \# \mathbf{ev}^{-1}(\vec x,\vec y).
\]
So, we are counting holomorphic curves through a generic collection of marked
points, just like Welschinger. The sign of a given point $\u = (u,\vec z,\vec
w)\in \mathbf{ev}^{-1}(\vec x,\vec y)$ depends on whether or not the isomorphism
\begin{equation}\label{eq:omap}
d{\bf ev}_\u : \det(T\M_{\k,l}(L,\Sigma,\d))_\u \stackrel{\sim}{\rightarrow}
{\bf ev}^* \det\left(T\left(L^{|\k|}\times X^l\right)\right)_\u
\end{equation}
agrees with the isomorphism of Theorem \ref{thm:or} up to the
action of the multiplicative group of positive real numbers. We would like to
reduce this sign to the Welschinger sign associated with the curve $\u,$ up
to a universal correction factor. For this purpose, we apply Proposition
\ref{prop:sescr} to a particular short exact sequence of Cauchy-Riemann boundary
value problems. Indeed, we take the underlying short exact sequence of vector
bundle pairs
\begin{equation}\label{eq:sescru}
0 \longrightarrow (TD^2,T\partial D^2) \stackrel{du}{\longrightarrow}
(u^*TX,u^*TL) \longrightarrow (N^X_u,N^L_u) \longrightarrow 0.
\end{equation}
Pulling back the $Pin$ structure on $TL$ induces a $Pin$ structure $\p_u$ on
$u^*TL.$  We equip $(u^*TX,u^*TL)$ with the linearized $\bar\partial_J$ operator
$D_u = D\bar\partial_J|_u$ and we denote by $D'_u$ and $D''_u$ the natural
operators it induces on the other terms of the short exact sequence. Note that
$T\partial D^2$ is orientable and has a natural orientation, so we are in the
situation of Definition \ref{dfn:sescr}. The natural orientation on $T\partial
D^2$ also induces an orientation on $\ker D'_u.$

Before continuing, we introduce some notation for configuration space. Define
the configuration space of $k$ boundary points and $l$
interior points of the disk to be
\begin{equation*}
\mathcal C_{k,l} : = (\partial D^2)^k \times (D^2)^l\setminus \Delta.
\end{equation*}
Thinking of $\M_{k,l}$ as the fiber product
\begin{equation*}
\M_{k,l} : = \widetilde\M_{k,l+2} \times_{X^2 \times \mathcal C_{0,2}} (A
\times B \times s_0 \times \ell),
\end{equation*}
we obtain Diagram \ref{eq:qtdi}.
\begin{figure}[ht]
\centering
\includegraphics{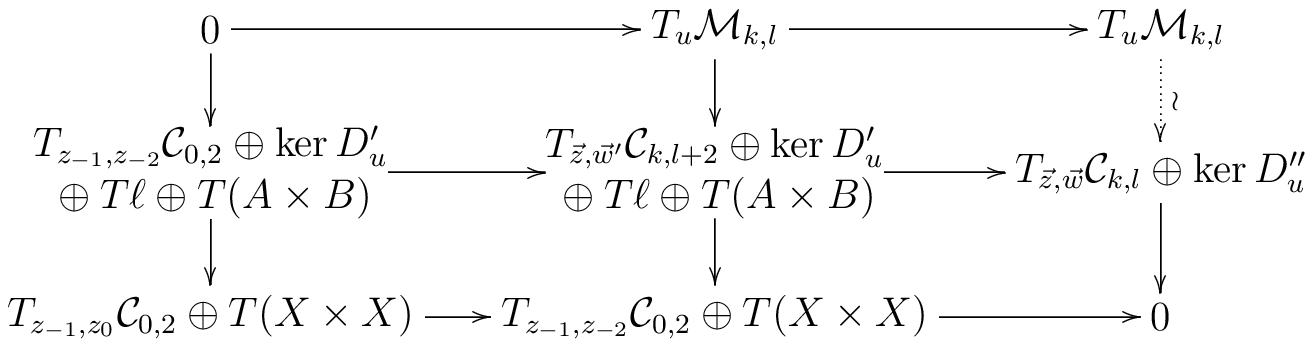}
\caption{}
\label{eq:qtdi}
\end{figure}
The central column of Diagram \ref{eq:qtdi} is the short-exact sequence for the
tangent space of the fiber product. The main content of the central row is the
short exact sequence of solutions of the the short exact sequence of
Cauchy-Riemann boundary value problems \eqref{eq:sescru}. This short exact
sequence exists because by assumption the Cauchy-Riemann operator at each term
of the sequence is surjective and we can apply the snake lemma from homological
algebra. Diagram \ref{eq:qtdi} shows how to construct a natural isomorphism
\begin{equation}\label{eq:mkl}
T\M_{k,l} \simeq T_{\vec z,\vec w}\mathcal C_{k,l}\oplus \ker
D_u''.
\end{equation}

First, suppose $u^*TL$ is orientable. Choose an orientation $\Omega_0$ on
$u^*TL.$ Then, $\Omega_0$ together with the complex orientation of
$T\partial D^2$ induces an orientation on $N_u^L.$ Furthermore, as in the proof
of Proposition \ref{prop:or}, $\Omega_0$ induces an orientation on
$\mathbf{ev}^* \det\left(T\left(L^{|\k|}\times X^l\right)\right)_v.$ If
$u^*TL$ is not orientable, as in the proof of Proposition \ref{prop:or},
$\mathbf{ev}^* \det\left(T\left(L^{|\k|}\times X^l\right)\right)_v$ admits
a canonical orientation. Since $D_u$ is surjective by assumption, by Proposition
\ref{prop:orl}, $\p_u$ and $\Omega_0$ induce an orientation $\Omega$ on
$\det(\ker(D_u)) \simeq \det(D_u).$  The central row
of Diagram \ref{eq:qtdi} shows that the $\Omega$ orientation on $\ker(D_u)$
induces an orientation of $T\M_{k,l} \simeq T_{\vec z,\vec
w}\mathcal C_{k,l}\oplus \ker D_u''$ as the quotient of an oriented vector
space by an oriented vector space. The sign of $\u$ is now the sign of map
\eqref{eq:omap} with respect to the orientations of the domain and range just
outlined.

On the other hand, the one dimension bundle $N_u^L$ carries a canonical
$Pin$ structure $\p_u''$ by Lemma \ref{lem:p1}.  Since $D_u''$ is surjective,
invoking Proposition \ref{prop:orl} again, $\p_u''$ and $\Omega_0$ induce an
orientation on $\ker(D_u'')$ and hence on $T\M_{k,l}$ by isomorphism
\eqref{eq:mkl}. By Proposition \ref{prop:sescr}, the $\Omega''$ orientation on
$T\M_{k,l}$ agrees with the $\Omega$ orientation on $T\M_{k,l}$ if and only if
$\p_u$ is compatible with the short exact sequence \eqref{eq:sescru}.

Recall that Welschinger's invariant counts curves with sign
determined by the parity of the isolated real double points. By
the adjunction formula, a rational curve $u$ of degree $d$ in a
symplectic four manifold has a topologically determined total
number of double points
\[
\delta(u) = \frac{d \cdot d - c_1(d) - 2}{2}.
\]
Complex double points come in pairs. So, the parity of the real
non-isolated double points is determined by the the parity of the
real isolated double points. On the other hand, the parity of the
real non-isolated double points determines the parity of the
number of twists of the real part of the curve about itself. This
exactly determines when the $Pin$ structure $\p_u$ is compatible
with the short exact sequence \eqref{eq:sescru}.

Finally, we claim that the map \eqref{eq:omap} always preserves
orientation if we consider the $\Omega''$ orientation on
$T\M_{k,l}$ together with the $\Omega_0$ orientation on
$\mathbf{ev}^* \det\left(T\left(L^{|\k|}\times
X^l\right)\right)_v.$ Indeed, by Lemma \ref{lem:surj}, if we vary
$D''$ to the standard complex linear Cauchy-Riemann operator on
$N^X_u$ and vary $N^L_u$ to a standard boundary condition keeping the marked
points distinct, the evaluation map will remain surjective during the whole
variation. So, the sign is standard for a given ordering of the marked
points. At this point, we use Definition \ref{dfn:tci}, which twists the
orientation of $T\M_{k,l}$ by $\sign(\varpi).$ Indeed, switching
the order of the marked points induces a change of sign of the evaluation
map \eqref{eq:omap} since each marked point leads to a codimension $n-1$
condition and $n = 2$ is even. The twisting
cancels this sign change.

The case $n=3$ is very similar. The same exact sequence
\eqref{eq:sescru} again plays a central role. Since $N^X_u$ is now
two dimensional, we need to define a $Pin$ structure on $N^L_u$ with which
$\p_u$ may or may not be compatible. In \cite{We2}, Welschinger does exactly
that using the splitting of $N^X_u$ into holomorphic line bundles. It is
important in that paper that $X$ be convex so that $J$ may be taken to be
integrable. Then the Cauchy-Riemann operator $D''_u$ is the standard one,
so the evaluation map is also standard and has  a standard sign. If $D''_u$
were not standard, we could not apply Lemma \ref{lem:surj} as before since
$N^X_u$ is no longer one dimensional. In \cite{We3}, Welschinger modifies the
definition of spinor states to take into account possible changes of
orientation arising from walls where the evaluation map is not surjective. This
allows him to extend the definition of his invariants to general
strongly-semipositive real symplectic manifolds $X.$
\end{proof}

Now we turn to the proof of the calculation in Example \ref{ex:q}. We
generalize Kontsevich's idea for calculating the closed Gromov-Witten
invariants of the quintic threefold \cite{K} to deal with the open case as well.
Extending previous notation in the case $n=1,$ we denote by $\tau$ the
tautological line bundle of $\C P^n$ equipped with its canonical complex
structure and we denote by $\tau_{\R}$ the tautological line bundle of $\R P^n.$
Furthermore, we let
\[
c' : \C P^n \rightarrow \C P^n, \qquad \tilde c': \tau \rightarrow \tau,
\]
denote complex conjugation and the bundle-map of $\tau$ covering complex
conjugation respectively. Throughout the following discussion, we use $\Gamma$
to denote holomorphic sections or holomorphic sections satisfying a real
boundary condition.

Let $s \in \Gamma(\tau_{\C P^4}^{*\otimes 5})$ be a real section, i.e.
\[
\tilde c' \circ s\circ  c' = s.
\]
We take
\[
X_s  = \{ s = 0 \} \subset \C P^4, \quad \omega_s  = \omega_{FS}|_X, \quad
\phi_s = c'|_{X_s}, \quad
L_s = Fix(\phi_s) \subset X_s.
\]
That is, $X_s$ is a quintic threefold equipped with the symplectic form
induced by the restriction of the Fubini-Study form of $\C P^4.$ Moreover, $X_s$
is equipped with an anti-symplectic involution $\phi_s$ and $L_s =
Fix(\phi_s)$ is the corresponding Lagrangian submanifold. Choosing $s$
generically, we may assume that $X_s$ is a smooth manifold. When it does not
lead to confusion, we may drop the subscript $s.$

We define a bundle $\mathcal F_d$ over the moduli space $\M_{0,0}(\R
P^4,D^2,d)$ of disks in $\C P^4$ with boundary in $\R P^4$ by
specifying its fibers,
\begin{equation*}
\mathcal F_{d_u} : = \Gamma(u^*\tau^{*\otimes5},u^*\tau_\R^{*\otimes5}), \quad u
\in \M_{0,0}(\R P^4,D^2,d).
\end{equation*}
By restricting to the image of each curve, $s$ induces a section $\hat s$ of
$\mathcal F_d$ that vanishes exactly on those curves entirely contained in
$X_s.$
\begin{lemma}
The total space of $\mathcal F_d$ is orientable for each $d.$ For $d$ odd,
$\mathcal F_d$ is an orientable vector bundle.
\end{lemma}
\begin{proof}
Let $u \in \M_{0,0}(\R P^4,D^2,d)$ and $\xi \in
\mathcal F_u.$ After
choosing a connection on $\mathcal F_d,$ there exists a canonical isomorphism
\begin{align*}
T_{(u,\xi)} \mathcal F_d  &\simeq
\Gamma\left(u^*\tau^{*\otimes5},u^*\tau_\R^{*\otimes5}\right) \oplus
\Gamma\left(u^*T\C P^4,u^*T \R P^4\right) \\
&\simeq \Gamma\left(u^*\left(\tau^{*\otimes5}\oplus T\C
P^4\right),u^*\left(\tau_\R^{*\otimes5}\oplus T \R P^4\right)\right).
\end{align*}
Since $\tau_\R^{*\otimes5}\oplus T\R P^4$ is orientable, after choosing an
orientation, Proposition \ref{prop:orl} gives a canonical orientation on each
tangent space. It is not hard to see that this orientation varies continuously
with $u$ and $\xi.$

When $d$ is odd, we have
\[
\mathcal F_{d_u} : =
\Gamma\left(u^*\tau^{*\otimes5},u^*\tau_\R^{*\otimes5}\right) \simeq
\Gamma\left(\tau^{*\otimes 5d}|_{D^2}, \tau_\R^{*\otimes5d}\right),
\]
where we think of $D^2$ as one hemisphere of $\C P^1$ with boundary $\partial
D^2 = \R P^1.$ Since $5d$ is odd, again Proposition \ref{prop:orl} gives each
$\mathcal F_u$ a canonical orientation that varies continuously with $u.$
\end{proof}

Let $V \rightarrow B$ be an orientable real vector bundle. We denote by $e(V)$
the Euler class of $V.$
\begin{prop}\label{prop:Keq}
Suppose $\hat s$ is transverse to the zero
section of $\mathcal F.$ Let
$d$ be odd. Then
\[
N_{D^2,d,0,0}  = e(\mathcal F_d)
\]
\end{prop}
\begin{rem}
This proposition should still hold true when $\hat s$ is not transverse to the
zero section of $\mathcal F.$ However, the proof will be slightly more
complicated. The author plans to address this issue in a future paper that will
calculate the invariants in higher degrees as well.
\end{rem}
\begin{rem}
If $d$ is even, when the details of the necessary corrections from real curves
with empty real part are worked out, an argument similar to the proof of
Proposition  \ref{prop:Keq} should show that $N_{D^2,d,0,0}$ is zero. Indeed,
$N_{D^2,d,0,0}$ should be given by the self-intersection number of the zero
section of $\mathcal F_d.$ Since $\dim \M_{0,0}(\R P^4,D^2,d) = 5d + 1,$ which
is odd when $d$ is even, the self intersection number should be zero.
\end{rem}
\begin{proof}[Proof of Proposition \ref{prop:Keq}]
Let $N_{L_s}$ denote the normal bundle of $L_s$ in $\R P^4.$ By the
adjunction formula, it is isomorphic to $\tau_\R^{*\otimes5}|_{L_s}.$ Since
$N_{L_s}$ is one dimensional, by Lemma \ref{lem:p1} we may choose its $Pin^+$
structure canonically. Equip $\tau_\R^{*\otimes 5}$ with a $Pin^+$ structure
corresponding to the $Pin^+$ structure of $N_{L_s}$ under the isomorphism of the
adjunction formula. Choose $Pin^+$ structures on $TL_s$ and $T\R P^4$ compatible
with the short exact sequence
\[
 0 \longrightarrow TL_s \longrightarrow T\R P^4 \longrightarrow N_{L_s}
\longrightarrow 0.
\]
As in isomorphism \eqref{eq:mkl} of the proof of Theorem \ref{thm:eq},
we may identify
\begin{equation}\label{eq:miso}
T_u M_{0,0}(L_s,D^2,d) \simeq \ker D''_u,
\end{equation}
where $D''_u$ is the operator that $D_u$ induces on the
bundle pair $(N^{X_s}_u,N^{L_s}_u).$ Equip $N^{L_s}_u$ with the $Pin^+$
structure induced by the short exact sequence
\[
0\longrightarrow T\partial D^2 \stackrel{du}{\longrightarrow} u^*TL_s
\longrightarrow N^{L_s}_u \longrightarrow 0.
\]
Then, by Diagram \ref{eq:qtdi} and Proposition \ref{prop:sescr}, we may assume
that isomorphism \eqref{eq:miso} is orientation preserving when $\ker D''_u$ is
given the canonical orientation of Proposition \ref{prop:orl}. From
Diagram \ref{eq:Ndi}
\begin{figure}[h]
\centering
\includegraphics{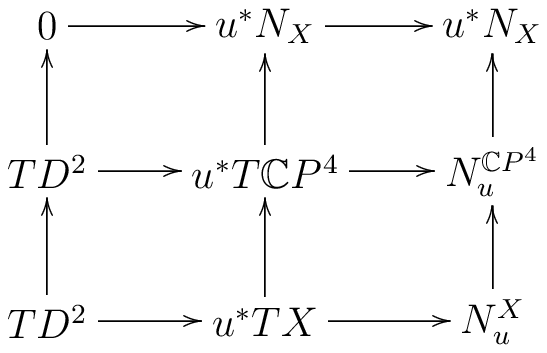}
\caption{}
\label{eq:Ndi}
\end{figure}
along with its conjugation invariant part, we deduce a short exact sequence of
Cauchy-Riemann boundary value problems
\begin{equation}\label{eq:sesNu}
0 \longrightarrow (N_u^X,N_u^L) \longrightarrow (N_u^{\C P^4},N_u^{\R P^4})
\longrightarrow (u^* N_X, u^* N_L) \longrightarrow 0.
\end{equation}
The Cauchy-Riemann operators at each term of the sequence are induced by the
rows of Diagram \ref{eq:Ndi}. The conjugation invariant parts of the rows of
Diagram \ref{eq:Ndi} induce $Pin^+$ structures on each of the boundary
conditions in short exact sequence \eqref{eq:sesNu}. The induced $Pin^+$
structures are compatible with the exact sequence because the already chosen
$Pin^+$ structures on the conjugation invariant parts of the first two columns
of Diagram \ref{eq:Ndi} are compatible.  Recall from the proof of the
adjunction formula that the isomorphism $N_{X_s} \simeq \tau^{*\otimes
5}|_{X_s}$ is given by the differential $ds.$ So, we have a diagram
\[
\xymatrix@-1pc{
    &&  &&          && \Gamma(u^* N_X, u^* N_L)
\ar[dd]^{\wr\:ds} \ar[drr]\\
0 \ar[rr] && \Gamma(N_u^X,N_u^L) \ar[rr] &&\Gamma(N_u^{\C P^4},N_u^{\R P^4})
\ar[rru]\ar[rrd]^{d\hat s} &&   &&  0.\\
    &&  &&          && \mathcal F_u \ar[urr] \ar@{=}[d]\\
    &&  &&          &&
\Gamma(u^*\tau^{*\otimes5},u^*\tau_\R^{*\otimes5}) \\
}
\]
Since $\C P^4$ is convex, $\Gamma(N_u^{\C P^4},N_u^{\R P^4})$ has expected
dimension. Since $\tau^{*\otimes5}$ is a line bundle, its sections always have
expected dimension. By assumption, $d\hat s$ is an isomorphism. So, by the
snake lemma, $\Gamma(N_u^X,N_u^L)$ must have expected dimension, i.e., $0.$ So,
its orientation is just a sign. Since all $Pin^+$ structures in the above
diagram have been chosen compatibly, the sign is given exactly by the sign of
$d\hat s,$ as claimed.
\end{proof}

\begin{proof}[Proof of Example \ref{ex:q}]
The section $s_F \in \Gamma(\tau^{*\otimes 5})$ defining the Fermat quintic does
not satisfy the assumptions of Proposition \ref{prop:Keq} even in degree $1.$
However, an elementary transversality argument shows that we may choose a nearby
section $s$ which does, in degree $1.$ Note that if $s$ is not sufficiently
close to $s_F,$ the topology of $L_s$ could be different
from that of $L_{s_F}.$ However, by Moser's method, if $s$ is sufficiently
close to $s_F,$ we know that $(X_s,\omega_s)$ is symplectomorphic to
$(X_{s_F},\omega_{s_F})$ by a $\phi_s - \phi_{s_F}$ equivariant
symplectomorphism. In particular, this symplectomorphism sends $L_s =
Fix(\phi_s)$ to $L_{s_F} = Fix(\phi_{s_F}).$ So, fixing some $s$ sufficiently
close to $s_F,$ we may think of the deformation
of $s_F$ to $s$ as a deformation of $\phi_{s_F}$-invariant complex structure on
$X_{s_F},$ which leaves the invariants unchanged. By Proposition \ref{prop:Keq},
it suffices to calculate $e(\mathcal F_1).$ Let $G(k,n)$ denote the Grassmannian
of real oriented $k$ planes in $n$ space and let $\tau_G$ denote its
tautological bundle. It is not hard to see that
\[
 \M_{0,0}(\R P^4,D^2,1) \simeq  G(2,5), \qquad \mathcal F_1 \simeq
Sym^5(\tau_G).
\]
Applying the splitting principle, we calculate the Pontryagin class
\begin{equation}\label{eq:p3}
p_3(Sym^5(\tau_G)) = 225 p_1(\tau_G)^3.
\end{equation}
Then, taking square roots, we have
\[
e(Sym^5(\tau_G)) = 15 e(\tau_G)^3 = 15 e(\tau_G^{\oplus 3}).
\]
Here, we have to include $G(2,5)$ into $G(2,n)$ for $n$ sufficiently large so
that both sides of equation \eqref{eq:p3} are not just zero. We use the unique
factorization property of the polynomial ring $H^*(G(2,\infty))$ to justify
taking square roots on both sides.

Finally, since there are two oriented $2$-planes in the intersection of three
generic hyperplanes in $\R^5,$ we know that
\begin{equation}\label{eq:int2}
\int_{G(2,5)} e(\tau_G^{\oplus 3}) = 2 \text{  or  } 0.
\end{equation}
To show the integral is actually $2,$ we proceed as follows.
Let $\hat G(k,n)$ denote the Grassmannian of unoriented $k$-planes in $n$-space
and let $\hat\tau_G$ denote its tautological bundle. Note that $\pi : G(2,5)
\rightarrow \hat G(2,5)$ is the orientation cover. Moreover, $\pi^*\hat \tau_G
\simeq \tau_G$ and $w_1(\hat \tau_G) = w_1(T\hat G(2,5)).$ So, both points
count the same and integral \eqref{eq:int2} is $2$ as desired.
\end{proof}

\appendix
\section{Kuranishi structures}\label{sec:aks}

In this appendix, we briefly review the definition of a space with Kuranishi
structure and perturbations thereof, as introduced in \cite{FO} and extended in
\cite{FOOO}. We essentially follow the conventions of \cite[Appendix 2]{FOOO}.
In the following discussion, we take $X$ to be a compact metrizable space.
\begin{dfn}\label{dfn:ks}
A \emph{Kuranishi structure with corners} on $X$ of dimension $d$ consists of
the following data:
\begin{enumerate}
\renewcommand{\theenumi}{(\arabic{enumi})}
\renewcommand{\labelenumi}{\theenumi}
\item
For each point $p \in X,$
\begin{enumerate}
\renewcommand{\theenumii}{(\arabic{enumi}.\arabic{enumii})}
\renewcommand{\labelenumii}{\theenumii}
\item
A smooth manifold with corners $V_p$ and a smooth vector bundle $E_p
\rightarrow V_p$ such that $\dim V_p - \rank E_p = d.$
\item
A finite group $\Gamma_p$ acting on $E_p \rightarrow V_p.$
\item
An $\Gamma_p$-equivariant smooth section $s_p$ of $E_p.$
\item
A homeomorphism $\psi_p$ from $s_p^{-1}(0)/\Gamma_p$ to a neighborhood of $p$ in
$X.$
\end{enumerate}
\item
For each $p \in X$ and for each $q \in \im \psi_p,$
\begin{enumerate}
\renewcommand{\theenumii}{(\arabic{enumi}.\arabic{enumii})}
\renewcommand{\labelenumii}{\theenumii}
\item
An open subset $V_{pq} \subset V_q$ containing $\psi_q^{-1}(q).$
\item
A homomorphism $h_{pq}:\Gamma_q \rightarrow \Gamma_p.$
\item
An $h_{pq}$-equivariant embedding $\varphi_{pq}: V_{pq} \rightarrow V_p$ and an
injective $h_{pq}$-equivariant bundle map $\hat\varphi_{pq}: E_q|_{V_{pq}}
\rightarrow E_p$ covering $\varphi_{pq}.$
\end{enumerate}
\end{enumerate}
Furthermore, the above data should satisfy the following compatibility
conditions:
\begin{enumerate}
\renewcommand{\theenumi}{(C\arabic{enumi})}
\renewcommand{\labelenumi}{\theenumi}
 \item
 $\hat\varphi_{pq}\circ s_q = s_p \circ \varphi_{pq}.$
 \item
 $\psi_q = \psi_p \circ \varphi_{pq}.$
 \item
 If $r \in \psi_q\left(s_q^{-1}(0)\cap V_{pq}/\Gamma_{q}\right),$ then in a
sufficiently small neighborhood of $r,$
 \[
 \varphi_{pq}\circ\varphi_{qr} = \varphi_{pr}, \qquad
\hat\varphi_{pq}\circ\hat\varphi_{qr} = \hat\varphi_{pr}.
 \]
\end{enumerate}
\end{dfn}
We note that occasionally we may regard $\psi_p$ as a map
from $s_p^{-1}(0)$ to $X$ by composing with the quotient map $V_p \rightarrow
V_p/\Gamma_p.$

A crucial ingredient in the construction of a perturbation of a space with
Kuranishi structure is the notion of the tangent bundle of a Kuranishi
structure. We take the following definition from \cite[Section 5]{FO}.
\begin{dfn}\label{dfn:ktb}
A \emph{tangent bundle} for a Kuranishi structure consists of a collection of
vector bundle isomorphisms
\[
 \Phi_{pq}: N_{U_p}U_q \stackrel{\sim}{\rightarrow} E_p|_{V_{pq}}/E_q|_{V_{pq}}
\]
covering the embeddings $\varphi_{pq}.$ Furthermore, if
\[
q \in \im \psi_p,\qquad
r
\in \psi_q\left(s_q^{-1}(0)\cap V_{pq}\right),
\]
then in a sufficiently small neighborhood of $r,$ we have a commutative diagram
\begin{equation*}
\xymatrix{
0 \ar[r] & N_{V_q}V_r \ar[r] \ar[d]^{\Phi_{qr}} & N_{V_p}{V_r} \ar[r]
\ar[d]^{\Phi_{pr}} & N_{V_p}V_q \ar[r]\ar[d]^{\Phi_{pq}} & 0 \\
0 \ar[r] & E_q/E_r \ar[r] & E_p/E_r \ar[r] & E_p/E_q \ar[r] & 0.
}
\end{equation*}
\end{dfn}

We also need the notion of an orientation for a Kuranishi structure, which
again comes from \cite[Section 5]{FO}.
\begin{dfn}
An \emph{orientation} of a Kuranishi structure with tangent bundle consists of
a family of trivializations of $\det(TV_p) \otimes \det(E_p)$ compatible with
the isomorphisms
\[
\det(TV_q) \otimes \det(E_q)|_{V_{pq}} \stackrel{\sim}{\rightarrow}\det(TV_p)
\otimes \det(E_p)|_{V_{pq}}
\]
induced by $\Phi_{pq}.$
\end{dfn}

Without providing full detail, we remind the reader of certain
definitions relating to multisections that are used in Section
\ref{sec:km} of this paper. For details, see \cite[Section 3]{FO}. In the
following, for $Z$ a space, we denote by $\mathcal S^\ell(Z)$ its
$\ell^{th}$ symmetric power. That is
\[
\mathcal S^\ell(Z) : = Z^\ell/S_\ell,
\]
where $S_\ell$ is the group of permutations of $\ell$ objects acting on
$Z^\ell$ by permuting the factors. Let $E \rightarrow V$ be a vector bundle. If
$U \subset V$ is sufficiently small, then $E|_U$ is trivial. So, if $\rank E =
r,$ a section of $E$ over $U$ may be specified by a map $U \rightarrow \R^r.$
Analogously, a multisection $s$ of $E$ over $U$ may be specified by a map $s_U:
U \rightarrow \mathcal S^\ell(\R^r).$ Note that globally, the multiplicity
$\ell$ can change. By definition, the multisection $s$ is said to be smooth if
after possibly shrinking $U,$ there exists a smooth lifting of $s_U$ to
the Cartesian product,
\[
\tilde s_U: U\rightarrow \left(\R^k\right)^\ell.
\]
The components of this lifting locally define $\ell$ sections $s_U^i$ of $E.$ We
call the $s^i_U$ branches of $s$ over $U.$ If $E \rightarrow V$ is a
$\Gamma$-equivariant vector bundle, then there is a natural notion of a
$\Gamma$-equivariant multisection coming from the induced action on the
symmetric power. Note that for $\Gamma$-equivariant smooth sections, we do not
require the local lifts $\tilde s_U$ to be $\Gamma$-equivariant. We call a
smooth multisection
transverse if each branch of each local lifting is transverse. The vanishing set
of a multisection $s$ is defined locally by
\[
 s^{-1}(0) \cap U = \bigcup_i \left(s^i_U\right)^{-1}(0).
\]
If $s$ is transverse and sufficiently generic, then $s^{-1}(0)$ admits a smooth
triangulation. If we fix a trivialization of $\det(E) \otimes \det(V),$
then the vanishing set of any smooth section is oriented. So, $s^{-1}(0)$
actually defines a rational singular chain by weighting each simplex of its
triangulation by the signed count of branches $s^i_U$ that vanish on it,
divided by $\ell.$

Finally, we need to recall the notion of a
\emph{good coordinate system} introduced in \cite[Section 6]{FO}. Fix a
Kuranishi structure on $X.$ We denote the various parts of the
Kuranishi structure by the same symbols as in Definition \ref{dfn:ks}. For $V_p'
\subset V_p,$ we denote by $E'_p,\psi'_p, s'_p,$ etc. the restrictions of all
the related parts of the Kuranishi structure. A good coordinate system specifies
a finite ordered set $P\subset X$ and $V'_p \subset
V_p$ for each $p \in P$ such that
\[
 X \subset \cup_{p\in P} \im \psi'_p.
\]
Furthermore, for $q,p \in P$ such that $q<p,$ it specifies a neighborhood
\[
V_{pq}'  \supset \psi_q^{-1'}(\im \psi'_p),
\]
an embedding $\varphi_{pq}': V_{pq}' \hookrightarrow V_p'$ and an injective
bundle map
\[
\hat
\varphi_{pq} : E_q'|_{V_{pq}'} \rightarrow E_p'
\]
covering $\varphi'_{pq}.$ Of course, we must require $s_p'\circ \varphi_{pq}' =
\hat \varphi_{pq}'\circ s_q'.$ Also, $\varphi_{pq}'$  (resp. $\hat
\varphi_{pq}'$) must respect the actions of $\Gamma_q$ and $\Gamma_p$ in such a
way as to define a map of the quotient orbifolds (resp. orbi-bundles). A few
additional technical conditions ensure requisite compatibility.

The following is a restatement of \cite[Theorem 6.4]{FO}.
\begin{thm}\label{thm:kp}
Let $(P,V_p',\psi_p',s_p',\varphi_{pq}',\hat\varphi_{pq}')$ be a good coordinate
system on a space  with Kuranishi structure $(X,\K).$ Suppose that $\K$ has a
tangent bundle in the sense of Definition $\ref{dfn:ktb}.$ Then, for each $p \in
P,$ there exists a sequence of smooth $\Gamma_p$-equivariant multisections
$s_{p,n}'$ of $E_p$ such that
\begin{enumerate}
\renewcommand{\theenumi}{(P\arabic{enumi})}
\renewcommand{\labelenumi}{\theenumi}
 \item \label{c:1}
 $s_{p,n}' \circ \varphi_{pq}' = \hat \varphi_{pq}'\circ s_{q,n}'.$
 \item
 $\lim_{n \to \infty} s_{p,n}' = s_p'.$
 \item
 $s_{p,n}'$ is transversal to $0.$
 \item \label{c:4}
 The restriction to $\im \varphi_{pq}'$ of the differential of the composition
of
any branch of $s_{p,n}'$ and the projection $E_p' \rightarrow E_p'/E_q'$
coincides
with the isomorphism $\Phi_{pq}' : N_{V_p'}V_q' \stackrel{\sim}{\rightarrow}
E_p'/E_q'.$
\end{enumerate}
\end{thm}
Fix $n$ sufficiently large. The rational chain given by the vanishing sets of
the $s_{p,n}'$ constructed in Theorem \ref{thm:kp} constitutes a transverse
perturbation of the space with Kuranishi structure $(X,\K).$

For the reader's convenience, we outline the proof of Theorem \ref{thm:kp}.
We use induction on the ordered set $P.$ Write $P = \{p_1,p_2,\ldots\}.$ Assume
the existence of perturbations $s_{p_i,n}'$ satisfying conditions
\ref{c:1}-\ref{c:4} for $i < j.$ The embeddings $\varphi_{p_jp_i}',$ the bundle
maps $\hat \varphi_{p_jp_i}'$ and the isomorphisms $\Phi_{p_jp_i},$ allow the
extension of $s_{p_i,n}'$ to a neighborhood of
\[
\bigcup_{i<j} \im \varphi_{p_jp_i}' \subset V_{p_j}.
\]
A small perturbation of the extension produces $s_{p_j,n}'$ as desired. Full
detail is given in \cite[Section 6]{FO}.

\end{document}